\documentclass{article}
\usepackage{graphicx}
\usepackage{float}
\usepackage{amsthm}
\usepackage{amsmath}
\usepackage{amssymb}
\usepackage{mathrsfs}
\usepackage{tikz-cd}
\tikzcdset{scale cd/.style={every label/.append style={scale=#1},
    cells={nodes={scale=#1}}}}
\usepackage{adjustbox}
\usepackage{stmaryrd}
\usepackage{enumitem}
\usepackage{placeins}
\usepackage{pdflscape}
\usepackage[colorlinks=true,citecolor=blue]{hyperref}
\usepackage{appendix}

\usepackage{comment}

\theoremstyle{definition}
\newtheorem{Definition}{Definition}[section]

\theoremstyle{plain}
\newtheorem{Theorem}[Definition]{Theorem}

\theoremstyle{plain}
\newtheorem{Proposition}[Definition]{Proposition}

\theoremstyle{definition}
\newtheorem{Notation}[Definition]{Notation}

\theoremstyle{plain}
\newtheorem{Lemma}[Definition]{Lemma}

\theoremstyle{plain}
\newtheorem{Corollary}[Definition]{Corollary}

\theoremstyle{plain}
\newtheorem{Conjecture}[Definition]{Conjecture}

\theoremstyle{definition}
\newtheorem{Example}[Definition]{Example}

\theoremstyle{remark}
\newtheorem{Remark}[Definition]{Remark}

\theoremstyle{plain}
\newcommand{\thistheoremname}{}
\newtheorem*{generic*}{\thistheoremname}
\newenvironment{namedthm*}[1]
  {\renewcommand{\thistheoremname}{#1}%
   \begin{generic*}}
  {\end{generic*}}

\newlength{\diagramwidth}

\title{On Étale Algebras and Bosonic Fusion 2-Categories}
\author{Hao Xu}
\date{November 2024}

\begin{document}

\bibliographystyle{alpha}

\maketitle

\begin{abstract}

We classify all connected and Lagrangian étale algebras in the Drinfeld center $\mathscr{Z}_1(\mathbf{2Vect}^\pi_G)$, where $G$ is a finite group and $\pi$ is a 4-cocycle on $G$. By Décoppet's result every bosonic fusion 2-category $\mathfrak{C}$ has its Drinfeld center equivalent to $\mathscr{Z}_1(\mathbf{2Vect}^\pi_G)$ for some $G$ and $\pi$. Combining this fact with classification of Lagrangian algebras in $\mathscr{Z}_1(\mathbf{2Vect}^\pi_G)$, we obtain a classification of bosonic fusion 2-categories. 

\end{abstract}

{\hypersetup{linkcolor=black}\tableofcontents}

\section{Introduction}

Étale algebras introduced in \cite{DMNO} are separable commutative algebras, which can be defined internal to any braided fusion category $\mathcal{B}$. An étale algebra $A$ is called connected if $\dim \mathbf{Hom}_{\mathcal{B}}(I,A) = 1$. Moreover, if its category of local modules $\mathbf{Mod}^{loc}_{\mathcal{B}}(A)$ is equivalent to $\mathbf{Vect}$, then $A$ is called a Lagrangian algebra.

Étale algebras are the essential ingredient to define the Witt equivalence between braided fusion categories. For any fusion category $\mathcal{C}$, its Drinfeld center $\mathcal{Z}_1(\mathcal{C})$ is a braided fusion category with Müger center $\mathcal{Z}_2(\mathcal{Z}_1(\mathcal{C})) \simeq \mathbf{Vect}$. A braided fusion category is called non-degenerate if its Müger center is equivalent to $\mathbf{Vect}$. Therefore, Drinfeld centers of fusion categories provides a natural family of non-degenerate braided fusion categories.

A non-degenerate braided fusion category $\mathcal{B}$ is called Witt trivial if it is equivalent to the Drinfeld center of some fusion category. Witt equivalence is an equivalence relation on non-degenerate braided fusion categories defined by modulo the Witt trivial class \cite{DMNO}. This equivalence relation drastically simplifies the classification of braided fusion categories, although it is still a challenging problem to classify all non-degenerate braided fusion categories up to Witt equivalence. Meanwhile, it also rises many interesting questions, such as generalizations to braided fusion categories with non-trivial Müger centers \cite{DNO}, Morita theoretic interpretation of Witt equivalences \cite{BJS,BJSS} and the group structure on the minimal modular extensions of symmetric fusion categories \cite{LKW}.

Étale algebras in braided fusion categories are also important from a physical perspective. They are used to describe the collections of condensable anyons in (2+1)D topological phases of matter \cite{Kon14}, so they are usually called \textit{condensable algebras} in the physical literature. Witt equivalence between braided fusion categories now corresponds to the equivalence of (2+1)D topological orders up to topological domain walls.

3D Dijkgraaf-Witten Theory, or finite gauge theory, is a well-studied example of (2+1)D topological order. As inputs there is a finite group $G$ and a 3-cocycle $\omega$ on $G$; mathematically, all line defects in this gauge theory form a braided fusion category $\mathcal{Z}_1(\mathbf{Vect}^\omega_G)$. Étale algebras condense to produce topological domain walls between $\mathcal{Z}_1(\mathbf{Vect}^\omega_G)$ and their categories of local modules. In particular, Lagrangian algebras classify all topological boundary conditions of $\mathcal{Z}_1(\mathbf{Vect}^\omega_G)$, i.e. fusion category $\mathcal{C}$ such that $\mathcal{Z}_1(\mathbf{Vect}^\omega_G) \simeq \mathcal{Z}_1(\mathcal{C})$. The classification of connected étale algebras in $\mathcal{Z}_1(\mathbf{Vect}^\omega_G)$ appeared first in \cite{DS16}.

\begin{Theorem}[{\cite[Theorem 2.19]{DS16}}]
    A connected étale algebra in $\mathcal{Z}_1(\mathbf{Vect}^\omega_G)$ is determined by:
    \begin{enumerate}
        \item A subgroup $H$ of $G$ up to conjugacy,
        
        \item A normal subgroup $N$ of $H$,
        
        \item A 2-cochain $\phi$ on $N$ such that $\mathrm{d} \phi = \omega|_N$,
        
        \item A 3-cocycle $\omega|_{H/N}$ on $H/N$ extending $\omega|_H$ along the canonical quotient homomorphism $H \twoheadrightarrow H/N$.
    \end{enumerate}
\end{Theorem}

In an unpublished exercise, Hao Zheng reformulated the data $\phi$ and $\omega|_{H/N}$ as monoidal functor structures on the canonical forgetful functor \[\mathbf{Vect}^{\pi|_H}_H \to \mathbf{Vect}^{\pi|_{N/N}}_{H/N}. \]

\begin{Theorem}[{\cite[Corollary 2.21]{DS16}}]
    A Lagrangian algebra in $\mathcal{Z}_1(\mathbf{Vect}^\omega_G)$ is determined by a subgroup $H$ of $G$ up to conjugacy and a A 2-cochain $\phi$ on $H$ such that $\mathrm{d} \phi = \omega|_H$.
\end{Theorem}

In this work, the author is motivated to categorify the above results. We will start with 4D Dijkgraaf-Witten Theory, where all surface and line defects form a 2-category $\mathscr{Z}_1(\mathbf{2Vect}^\pi_G)$ for some finite group $G$ and 4-cocycle $\pi \in \mathrm{H}^4(G;\Bbbk^\times)$. The braided fusion structure on this Drinfeld center has been examined in detail by Kong et al \cite{KTZ}. Étale algebras in $\mathscr{Z}_1(\mathbf{2Vect}^\pi_G)$ turn out to be $G$-crossed braided multifusion categories with associator twisted by $\pi$ (for which the precise definition will be given in Definition \ref{def:TwistedCrossedBraidedFusionCat}). We will classify all connected and Lagrangian étale algebras in $\mathscr{Z}_1(\mathbf{2Vect}^\pi_G)$ in Section \ref{sec:TwistedCrossedBraidedFusionCat}.

\renewcommand{\thistheoremname}{Theorem \ref{thm:ClassificationOfConnectedTwistedCrossedBraidedMultiFusionCategories}}

\begin{generic*}
    A connected {\'e}tale algebra in $\mathscr{Z}_1(\mathbf{2Vect}^\pi_G)$ is determined by
    \begin{enumerate}
        \item A subgroup $H$ of $G$ up to conjugacy,
        
        \item A normal subgroup $N$ of $H$,
        
        \item A braided fusion category $\mathcal{A}$ with an $H$-action $\gamma: H \to \mathbf{Aut}_{br}(\mathcal{A})$.
        
        \item A 4-group morphism $H/N \times \mathrm{B}^3 \Bbbk^\times \to \mathbf{BrPic}(\mathbf{Mod}(\mathcal{A}))$,
        
        \item An invertible 2-morphism: \[
            \begin{tikzcd}
                {H}
                    \arrow[d,"{(q,\pi|_H)}"']
                    \arrow[r,"\gamma"]
                & {\mathbf{Aut}_{br}(\mathcal{A})}
                    \arrow[d,"\mathbf{Z}"]
                \\ {H/N \times \mathrm{B}^3 \Bbbk^\times}
                    \arrow[r,"{(\Gamma,{-}\iota)}"']
                    \arrow[ur,Rightarrow,shorten >=3ex,shorten <=4ex]
                & {\mathbf{BrPic}(\mathbf{Mod}(\mathcal{A}))}
            \end{tikzcd}. \]
        The homotopy fiber of the above square produces a $\pi$-twisted $N$-crossed extension: \[\begin{tikzcd}
        {N \rtimes^\pi \mathrm{B}^2 \Bbbk^\times}
            \arrow[d]
            \arrow[r,"\theta"]
        & {\mathbf{Pic}(\mathcal{A})}
            \arrow[d,"\partial"]
        \\ {H}
            \arrow[r,"\gamma"']
            \arrow[ur,Rightarrow,shorten >=4ex,shorten <=4ex]
        & {\mathbf{Aut}_{br}(\mathcal{A})}
    \end{tikzcd}. \]
\end{enumerate}

    The corresponding $\pi$-twisted $G$-crossed braided multifusion category is given by \[\mathcal{C} := \mathbf{Fun}_H \left( G,\bigoplus_{g \in N} \theta(g) \right).\]
\end{generic*}

\renewcommand{\thistheoremname}{Corollary \ref{cor:ClassificationOfLagrangianTwistedCrossedBraidedFusionCategories}}

\begin{generic*}
    A connected {\'e}tale algebra in $\mathscr{Z}_1(\mathbf{2Vect}^\pi_G)$ (determined by the data stated in the theorem) is Lagrangian if and only if $\mathcal{A}$ is non-degenerate and $N=H$. In other word, a Lagrangian algebra in $\mathscr{Z}_1(\mathbf{2Vect}^\pi_G)$ is equivalent to a sequence $(H,\mathcal{A},\gamma,\theta,\phi)$ where \begin{enumerate}
        \item $H$ is a subgroup of $G$ determined up to conjugation,
        
        \item $\mathcal{A}$ is a non-degenerate braided fusion category with $H$-action $\gamma$,
        
        \item 3-group morphism $\theta: H \rtimes^\pi \mathrm{B}^2 \Bbbk^\times \to \mathbf{Pic}(\mathcal{A})$ corresponds to a $\pi|_H$-twisted $H$-crossed extension of $\mathcal{A}$,
        
        \item $\phi$ is an equivalence of $H$-crossed braided fusion 2-categories: \[\mathbf{2Vect}^{\pi|_H}_H \boxtimes_{\mathbf{Fun}(H,\mathbf{2Vect})} \mathscr{Z}^H_1(\mathbf{Mod}(\mathcal{A})) \simeq \mathbf{2Vect}_H \boxtimes \mathscr{Z}_1(\mathbf{Mod}(\mathcal{A}))\] by Theorem \ref{thm:LocalModulesOverTwistedCrossedBraidedFusionCategories}, or equivalently the cell filling the following diagram:
    \end{enumerate} \[\begin{tikzcd}
        {H}
            \arrow[d,"{\pi|_H}"']
            \arrow[r,"\gamma"]
        & {\mathbf{Aut}_{br}(\mathcal{A})}
            \arrow[d,"\mathbf{Z}"]
        \\ {\mathrm{B}^3 \Bbbk^\times}
            \arrow[r,"{-}\iota"']
            \arrow[ur,Rightarrow,shorten >=3ex,shorten <=4ex]
        & {\mathbf{BrPic}(\mathbf{Mod}(\mathcal{A}))} 
    \end{tikzcd}. \]
\end{generic*}

In the special case when the twist $\pi$ is trivial, we have an alternative description of the Drinfeld center $\mathscr{Z}_1(\mathbf{2Vect}_G)$. In \cite[Example 5.4.5]{D8}, Décoppet showed that $\mathbf{2Vect}_G$ and $\mathbf{2Rep}(G)$ are Morita equivalent. By \cite[Theorem 2.3.2]{D9}, Morita equivalent fusion 2-categories have braided equivalent Drinfeld centers.

On the other hand, Drinfeld et al essentially provided in \cite{DGNO} that $\mathbf{2Rep}(G)$ and $\Sigma \mathbf{Rep}(G)$ are equivalent as symmetric fusion 2-categories via equivariantization and de-equivariantization. Therefore their Drinfeld centers are also braided equivalent. Combining these two results, we have the following equivalence of three Drinfeld centers:

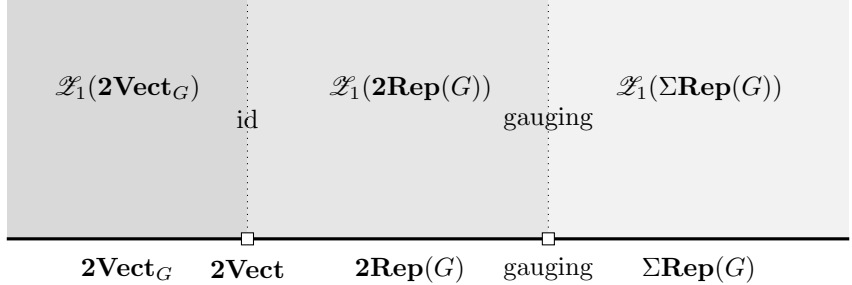
\begin{figure}[htbp]
\centering
\begin{tikzpicture}[scale=0.8]
    \fill[gray!30] (-4,0) rectangle (0,4) ;
    \fill[gray!20] (0,0) rectangle (5,4) ;
    \fill[gray!10] (5,0) rectangle (10,4) ;
    \draw[dotted] (0,4)--(0,0) node[midway] {$\mathrm{id}$} ;
    \draw[dotted] (5,4)--(5,0) node[midway] {\text{gauging}} ;
    \draw[very thick] (-4,0)--(10,0) ;
    \draw[fill=white] (-0.1,-0.1) rectangle (0.1,0.1) ;
    \draw[fill=white] (4.9,-0.1) rectangle (5.1,0.1) ;
    \node at (0,-0.5) {$\mathbf{2Vect}$} ;
    \node at (5,-0.5) {\text{gauging}} ;
    \node at (-2,2.5) {$\mathscr{Z}_1(\mathbf{2Vect}_G)$} ;
    \node at (2.7,2.5) {$\mathscr{Z}_1(\mathbf{2Rep}(G))$} ;
    \node at (7.5,2.5) {$\mathscr{Z}_1(\Sigma \mathbf{Rep}(G))$} ;
    \node at (-2,-0.5) {$\mathbf{2Vect}_G$} ;
    \node at (2.7,-0.5) {$\mathbf{2Rep}(G)$} ;
    \node at (7.5,-0.5) {$\Sigma \mathbf{Rep}(G)$} ;
\end{tikzpicture}
\caption{Braided equivalences of the centers}
\label{fig:EquivalentThreeCenters}
\end{figure}

Under these equivalences, one could classify all connected and Lagrangian étale algebras in $\mathscr{Z}_1(\mathbf{2Rep}(G))$ and $\mathscr{Z}_1(\Sigma \mathbf{Rep}(G))$.

\renewcommand{\thistheoremname}{Corollary \ref{cor:ConnectedAndLagrangianEtaleAlgebrasInZ(SigmaRepG)}}

\begin{generic*}
    An {\'e}tale algebra in $\mathscr{Z}_1(\Sigma \mathbf{Rep}(G))$ consisting of a braided multifusion category $\mathcal{M}$ and a braided functor $F:\mathbf{Rep}(G) \to \mathcal{M}$. This {\'e}tale algebra is connected if and only if $\mathcal{M}$ is fusion. 
        
    Its 2-category of local modules is equivalent to $\mathscr{Z}_1(\mathbf{Mod}(\mathcal{M}))$. Moreover, this {\'e}tale algebra is Lagrangian if and only if $\mathcal{M}$ is non-degenerate.
\end{generic*}

\renewcommand{\thistheoremname}{Theorem \ref{thm:ClassificationOfConnectedCrossedBraidedMultiFusionCategories}}

\begin{generic*}
    A connected {\'e}tale algebra in $\mathscr{Z}_1(\mathbf{2Rep}(G))$ is determined by
        \begin{enumerate}
            \item A subgroup $H$ of $G$ up to conjugacy, 
            
            \item A normal subgroup $N$ of $H$, 
            
            \item A braided fusion category $\mathcal{A}$ with an $H$-action $\rho: H \to \mathbf{Aut}_{br}(\mathcal{A})$,
            
            \item A 4-group morphism $\Gamma: H/N \to \mathbf{BrPic}(\mathbf{Mod}(\mathcal{A}))$.
            
            \item An invertible 2-morphism: \[\begin{tikzcd}
            {H}
                \arrow[d,twoheadrightarrow]
                \arrow[r,"\rho"]
            & {\mathbf{Aut}_{br}(\mathcal{A})}
                \arrow[d,"\mathbf{Z}"]
            \\ {H/N}
                \arrow[r,"\Gamma"']
                \arrow[ur,Rightarrow,shorten >=3ex,shorten <=4ex]
            & {\mathbf{BrPic}(\mathbf{Mod}(\mathcal{A}))}
        \end{tikzcd}, \] whose homotopy fiber is an $N$-crossed extension $\alpha: N \to \mathbf{Pic}(\mathcal{A})$. 
    \end{enumerate} 

    The corresponding $G$-crossed braided multifusion category is given by \[\mathcal{C} := \mathbf{Fun}_H \left( G,\bigoplus_{g \in N} \alpha(g) \right).\] Its 2-category of local modules is equivalent to $\mathscr{Z}_1(\mathbf{Mod}(\mathcal{C}^G))$. Moreover, this {\'e}tale algebra is Lagrangian if and only if $\mathcal{A}$ is non-degenerate and $N=H$.
\end{generic*}

Lastly, in Section \ref{sec:ClassificationOfBosonicFusion2Categories}, we use our classification of Lagrangian algebras in $\mathscr{Z}_1(\mathbf{2Vect}^\pi_G)$ to classify all bosonic fusion 2-categories. By \cite[Theorem 4.1.6]{D9}, every bosonic fusion 2-category $\mathfrak{C}$ has its Drinfeld center equivalent to $\mathscr{Z}_1(\mathbf{2Vect}^\pi_G)$ for some finite $G$ and $\pi$.

\renewcommand{\thistheoremname}{Theorem \ref{thm:ClassificationOfBosonicFusion2Categories}}

\begin{generic*}
    Any bosonic fusion 2-category $\mathfrak{C}$ is equivalent to the 2-category of right modules over its full center $[I,I]$ in $\mathscr{Z}_1(\mathfrak{C}) \simeq \mathscr{Z}_1(\mathbf{2Vect}^\pi_G)$. Since $[I,I]$ is a Lagrangian algebra in $\mathscr{Z}_1(\mathfrak{C})$, under this equivalence it is classified by $(H,\mathcal{A},\gamma,\phi)$ in Corollary \ref{cor:ClassificationOfLagrangianTwistedCrossedBraidedFusionCategories}.
\end{generic*}

In the proof we need a technical lemma on monoidal 2-categories enriched over braided monoidal 2-categories, which is added in the Appendix \ref{sec:M2CatEnrichedInB2Cat}.

\subsection*{Acknowledgements}
    The project was initiated by discussions with Liang Kong and Hao Zheng back in summer, 2021. The author was visiting Jamie Vicary's group at University of Cambridge as a Research Assistant by his generous support. Hao Zheng provided a reformulated proof of classification of connected étale algebras in $\mathcal{Z}_1(\mathbf{Vect}^\omega_G)$, which inspired the author to generalize the result to 2-categories later. The author obtained the first version of the classification theorem in 2022, and would like to thank Liang Kong and Hao Zheng for their valuable comments. 
    
    The author also would like to thank the hospitality of the Erwin Schrödinger International Institute for Mathematics and Physics (ESI) to invite the author to participate in the workshop ``\textit{Higher Structures and Field Theory}" in August 2022, during which the author met Alexei Davydov, Dmitri Nikshych, David Reutter and presented the first version of the classification theorem to them. The author would like to thank Alexei Davydov and Dmitri Nikshych particularly for their continuous support and encouragement.

    The main theorems of the current version requires a thorough understanding of the general theory of local modules in braided monoidal 2-categories. This extremely technical part has no chance to be completed without the help of Thibault Décoppet as a collaborator. The author would like to thank Thibault Décoppet for his patience and dedication to the project on local modules. He also pointed out the connection between the current work and the classification of fusion 2-categories. 

    The author first publicly announced the main theorems at the BIMSA-Tsinghua Quantum Symmetry Seminar on 13th June 2023 \cite{XuBIMSATalk1}, and later in another talk at BIMSA on 17th April 2024 \cite{XuBIMSATalk2}. The author would like to thank Hao Zheng and Yilong Wang for invitation and hosting of these talks.
    
    At the finishing stage of the project, the author was informed that Thibault Décoppet, Peter Huston, Theo Johnson-Freyd, Dmitri Nikshych, David Penneys, Julia Plavnik, David Reutter, Matthew Yu have just completed their classification of fusion 2-categories \cite{DHJFNPPRY} independently, and they will have a following paper on the equivariantization of higher categories in the future, which is also a key ingredient in the current work (see the end of Preliminaries). It would be more satisfying if this part of arguments could be made more explicit with the help of this coming paper.

    In October 2024, the author met David Reutter and Theo Johnson-Freyd at Hamburg. They provided valuable comments on formulations and proofs of the main theorems. In particular, they pointed out the importance of the twisted crossed extensions of braided fusion categories as a certain homotopy fiber, and reminded the author to include the base, i.e. a 4-group morphism playing the role of a crossed extensions of certain braided fusion 2-categories. The author would like to thank them for their help.
    
    The author was also informed by Liang Kong recently about Rui Wen's work on the string condensations in the (3+1)D finite gauge theory \cite{Wen24}. Rui Wen also studied certain étale algebras in $\mathscr{Z}_1(\mathbf{2Vect}_G)$ from a physical perspective. The author wishes to clarify any potential overlap with Rui Wen's work in the future and extends gratitude to Liang Kong for providing this information.

    The author was supported by DAAD Graduate School Scholarship Programme (57572629) and DFG Project 398436923 RTG 2491 ``\textit{Fourier Analysis and Spectral Theory}''.
\section{Preliminaries}

Fix a ground field $\Bbbk$ which is algebraically closed and of characteristic $0$. In the following context, $G$ is a finite group. Let us first recall the notion of fusion 1-categories \cite{Mu1,Mu2,Mu3,ENO02,EO03,Ost1,Ost2,ENO09,DGNO,EGNO}.

\begin{Definition}
    A \textbf{multifusion category} is a $\Bbbk$-linear monoidal category $\mathcal{C}$ such that:
    \begin{enumerate}
        \item $\mathcal{C}$ is additive and idempotent complete.
        
        \item $\mathcal{C}$ is finite semisimple, i.e. every object is a finite direct sum of simple objects.
        
        \item $\mathcal{C}$ is rigid, i.e. every object has a left and right dual.
    \end{enumerate}

    A multifusion category is called a \textbf{fusion category} if its unit object is simple.
\end{Definition}

\begin{Example}
    We denote the 1-category of finite dimensional vector spaces and linear maps over $\Bbbk$ as $\mathbf{Vect}$. It is equipped with a symmetric monoidal structure given by the tensor product of vector spaces and a monoidal unit given by $\Bbbk$. The associator, unitor and braiding are given by the canonical isomorphisms \[V \otimes (W \otimes U) \cong (V \otimes W) \otimes U,\] \[\Bbbk \otimes V \simeq V \simeq V \otimes \Bbbk, \] \[ V \otimes W \simeq W \otimes V, \] for any finite dimensional vector spaces $U,V,W$, respectively.
\end{Example}

\begin{Example}
    $\mathbf{Vect}_G$ is the fusion category of $G$-graded finite dimensional vector spaces, with monoidal product given by convolution on $G$. More explicitly, the objects are finite dimensional vector spaces $V = \bigoplus_{g \in G} V_g$ with $V_g$ the $g$-component of $V$, and the morphisms are linear maps $f: V \to W$ such that $f(V_g) \subset W_g$ for all $g \in G$. The monoidal product is given by $V \otimes W$, with the following grading: \[(V \otimes W)_g = \bigoplus_{h \in G} V_h \otimes W_{h^{-1}g}\] and the monoidal unit is $\Bbbk$ graded on unit $e$. The associator, unitor and braiding are induced directly from the associator, unitor and braiding of $\mathbf{Vect}$.
\end{Example}

\begin{Example}
    We can twist the associator of $\mathbf{Vect}_G$ with a (normalized) group cocycle $\omega \in \mathrm{H}^3(G;\Bbbk^\times)$: for any $g,h,k \in G$, the associator is replaced as follows \[(\Bbbk_g \otimes \Bbbk_h) \otimes \Bbbk_k \simeq \Bbbk_g \otimes (\Bbbk_h \otimes \Bbbk_k),\]
    \[(1 \otimes 1) \otimes 1 \mapsto \omega_{g,h,k} \, 1 \otimes (1 \otimes 1).\]
    We denote the new fusion category as $\mathbf{Vect}^\omega_G$.
\end{Example}

\begin{Example}
    $\mathbf{Rep}(G):= \mathbf{Fun}(\mathrm{B} G,\mathbf{Vect}) $ is the symmetric fusion category of finite dimensional representations of group $G$, with monoidal product given by tensor product of the underlying vector spaces. It is equivalent to the category of right modules over the group algebra $\Bbbk [G]$. The $G$-action on tensor product is induced by the coalgebra structure $\Delta: g \mapsto g \otimes g$, $\epsilon: g \mapsto 1$.
\end{Example}

\noindent \textit{Drinfeld Center of Multifusion Categories.} Drinfeld introduced the notion of Drinfeld center in his unpublished note on the quantum double of Hopf algebra. The first definition of Drinfeld center can be found independently in \cite{JS91,Maj91}. 

\begin{Definition}
    Let $\mathcal{C}$ be a monoidal category. The \textbf{Drinfeld center} of $\mathcal{C}$, denoted by $\mathcal{Z}_1(\mathcal{C})$, is a braided monoidal category consists of:
    \begin{enumerate}
        \item An object is a pair $(x,\beta^x)$, where $x$ is an object in $\mathcal{C}$ and $\beta^x_z: x \otimes z \to z \otimes x$ is a natural isomorphism, often called a \textit{half-braiding} on $x$, such that the following hexagon commutes naturally in $y,z$: \[\begin{tikzcd}
            {}
            & {x \otimes (y \otimes z)}
                \arrow[r,"\beta^x_{y \otimes z}"]
            & {(y \otimes z) \otimes x}
                \arrow[rd,"\alpha_{y,z,x}"]
            & {}
            \\ {(x \otimes y) \otimes z}
                \arrow[rd,"\beta^x_y \otimes z"']
                \arrow[ru,"\alpha_{x,y,z}"]
            & {}
            & {}
            & {y \otimes (z \otimes x)}
            \\ {}
            & {(y \otimes x) \otimes z}
                \arrow[r,"\alpha_{y,x,z}"']
            & {y \otimes (x \otimes z)}
                \arrow[ur,"y \otimes \beta^x_{z}"']
            & {}
            \end{tikzcd};\]

        \item A morphism from $(x,\beta^x)$ to $(y,\beta^y)$ is a morphism $f:x \to y$ in $\mathcal{C}$ such that the following diagram commutes for any $z$ in $\mathcal{C}$:
        \[\begin{tikzcd}[sep=large]
            {x \otimes z} 
                \arrow[r,"\beta^x_z"]
                \arrow[d,"f \otimes 1_z"']
            & {z \otimes x}
                \arrow[d,"1_z \otimes f"]
            \\ {y \otimes z}
                \arrow[r,"\beta^y_z"']
            & {z \otimes y}
        \end{tikzcd};\]

        \item The monoidal product is given by \[(x,b^x) \otimes (y,b^y) := (x \otimes y, b^{x \otimes y}),\] where the half-braiding $b^{x \otimes y}$ is given by the following hexagon on object $z$ in $\mathcal{C}$:
        \[\begin{tikzcd}
            {}
            & {(x \otimes y) \otimes z}
                \arrow[r,"\beta^{x \otimes y}_z"]
            & {z \otimes (x \otimes y)}
                \arrow[rd,"\alpha^{-1}_{z,x,y}"]
            & {}
            \\ {x \otimes (y \otimes z)}
                \arrow[rd,"x \otimes \beta^y_z"']
                \arrow[ru,"\alpha^{-1}_{x,y,z}"]
            & {}
            & {}
            & {(z \otimes x) \otimes y}
            \\ {}
            & {x \otimes (z \otimes y)}
                \arrow[r,"\alpha^{-1}_{x,z,y}"']
            & {(x \otimes z) \otimes y}
                \arrow[ur,"\beta^x_{z} \otimes y"']
            & {}
            \end{tikzcd};\]      

        \item The monoidal unit is given by $(I,\beta^I)$, where $I$ is the monoidal unit of $\mathcal{C}$ with the canonical half-braiding $\beta^I$ defined by \[\begin{tikzcd}
            {I \otimes x}
                \arrow[rr,"\beta^I_x"]
                \arrow[rd,"\lambda_x"']
            & {}
            & {x \otimes I}
                \arrow[ld,"\rho_x"]
            \\ {}
            & {x}
            & {}
        \end{tikzcd} \] given on $x$ in $\mathcal{C}$.

        \item The associator, left unitor and right unitor are inherited from $\mathcal{C}$.
        
        \item The braiding on objects $(x,\beta^x)$ and $(y,\beta^y)$ is given by the half-braiding $\beta^x_y:x \otimes y \to y \otimes x$.
    \end{enumerate}
\end{Definition}

In \cite{Ost2}, Ostrik proved that Drinfeld center is invariant under Morita equivalence. Morita equivalence of tensor categories, introduced by \cite{FRS} in the context of rational 2-dimensional conformal field theory and by \cite{Mu1} in the context of subfactors, is reformulated in \cite{EO03} using Morita dual category for a module category \cite{Ost1,Ost2}.

\begin{Proposition}
    The follows are equivalent for multifusion categories $\mathcal{C}$ and $\mathcal{D}$:
    \begin{enumerate}
        \item $\mathbf{Mod}(\mathcal{C}) \simeq \mathbf{Mod}(\mathcal{D})$ as finite semisimple 2-categories.
        
        \item There exists a finite semisimple left $\mathcal{C}$-module category $\mathcal{M}$ and an equivalence of multifusion categories \[\mathcal{D}^{mp} \simeq \mathbf{End}_\mathcal{C}(\mathcal{M}),\] where $\mathcal{D}^{mp}$ denotes the multifusion category with the same underlying finite semisimple category as $\mathcal{D}$ but with the direction of monoidal product reversed.
        
        \item There exists a separable algebra $A$ in $\mathcal{C}$ and an equivalence of multifusion categories \[\mathcal{D} \simeq \mathbf{Bimod}_\mathcal{C}(A),\] where $\mathbf{Bimod}_\mathcal{C}(A)$ denotes the 2-category of bimodules over $A$ in $\mathcal{C}$, with monoidal product given by the relative tensor product $\otimes_A$.
    \end{enumerate} If any of the above equivalent conditions hold, we say that $\mathcal{C}$ and $\mathcal{D}$ are \textbf{Morita equivalent}.
\end{Proposition}

\begin{Proposition}
    Let $\mathcal{C}$ and $\mathcal{D}$ be two multifusion categories. If $\mathcal{C}$ and $\mathcal{D}$ are Morita equivalent, then their Drinfeld centers $\mathcal{Z}_1(\mathcal{C})$ and $\mathcal{Z}_1(\mathcal{D})$ are equivalent as braided multifusion categories.

    Conversely, if $\mathcal{Z}_1(\mathcal{C})$ and $\mathcal{Z}_1(\mathcal{D})$ are equivalent as braided multifusion categories, then $\mathcal{C}$ and $\mathcal{D}$ are Morita equivalent.
\end{Proposition}

\begin{proof}
    The first part follows from the observation that \[\mathcal{Z}_1(\mathcal{C}) \simeq \mathbf{Fun}_{\mathcal{C} \boxtimes \mathcal{C}^{mp}}(\mathcal{C},\mathcal{C}).\] The second part is more subtle and was first prove in \cite{ENO09}.
\end{proof}

\begin{Example}
    By \cite{Wil}, Drinfeld centers of both $\mathbf{Vect}_G$ and $\mathbf{Rep}(G)$ are given by the category of $G$-equivariant $G$-graded finite dimensional vector spaces \[\mathcal{Z}_1(\mathbf{Vect}_G) = \mathcal{Z}_1(\mathbf{Rep}(G)) = \mathbf{Fun}(\mathrm{LB} G,\mathbf{Vect}) \simeq \bigoplus_{[g] \in \mathrm{Cl}(G)} \mathbf{Rep}(Z(g)),\] where $\mathrm{LB} G:=\mathbf{Map}(S^1,\mathrm{B} G)$ is the free looping groupoid for the one-point delooping groupoid $\mathrm{B} G$, $\mathrm{Cl}(G)$ is the set of conjugacy classes in $G$ and $Z(g)$ is the centralizer subgroup of $g \in G$. 
\end{Example}

\begin{Example}
    The Drinfeld center of $\mathbf{Vect}^\omega_G$, where the twist $\omega$ is a 3-cocycle on $G$, is also computed in \cite{Wil} as the category of $G$-equivariant $G$-graded finite dimensional vector spaces with the monoidal product twisted by $\omega$. Its underlying category can be decomposed as \[\mathcal{Z}_1(\mathbf{Vect}_G) = \mathbf{Fun}(\mathrm{L} \widetilde{\mathrm{B} G}^\omega,\mathbf{Vect}) \simeq \bigoplus_{[g] \in \mathrm{Cl}(G)} \mathbf{Rep}(Z(g),\tau_g (\omega)),\] where $\widetilde{\mathrm{B} G}^\omega$ is the extension of $\mathrm{B} G$ by $\omega$: \[\begin{tikzcd}
        {\mathrm{B}^2 \Bbbk^\times}
            \arrow[r]
            \arrow[d]
        & {\widetilde{\mathrm{B} G}^\omega}
            \arrow[d]
            \arrow[r]
        & {*}
            \arrow[d]
        \\ {*}
            \arrow[r]
        & {\mathrm{B} G}
            \arrow[r,"\omega"]
        & {\mathrm{B}^3 \Bbbk^\times}
    \end{tikzcd},\] equivalently, $\widetilde{\mathrm{B} G}^\omega$ is one-point delooping of a 2-group constructed from extension of $G$ by $\mathrm{B} \Bbbk^\times$. Hence its linearization is the twisted group algebra $\mathbf{Vect}_G^\omega$. Then $\mathrm{L} \widetilde{\mathrm{B} G}^\omega$ is the free looping 2-groupoid of $\widetilde{\mathrm{B} G}^\omega$. Any element $g \in G$ gives a transgression $\tau_g: \mathrm{H}^3(G;\Bbbk^\times) \to \mathrm{H}^2(Z(g);\Bbbk^\times)$; the twist $\tau_g(\omega)$ is the transgression of $\omega$ under $\tau_g$. If $g$ and $g'$ are conjugate in $G$, then the twists $\tau_g(\omega)$ and $\tau_{g'}(\omega)$ are equivalent.
\end{Example}

\medskip

\noindent {\it Picard Groups.} The theory of Brauer-Picard groups for fusion categories is proposed in \cite{ENO09} to study the extension theory of fusion categories. In this section, we will briefly recall the basic definitions and its development in \cite{DN13,DSPS14,DSPS13,BJS,DN21,JFR}.

\begin{Definition}
    The Brauer-Picard 3-group $\mathbf{BrPic}(\mathcal{C})$ for a multifusion category $\mathcal{C}$ is the maximal sub-3-group of invertible objects in the multifusion 2-category $\mathbf{Bimod}(\mathcal{C})$ of $\mathcal{C}$-bimodule categories, under the relative Deligne tensor product $\boxtimes_{\mathcal{C}}$.
\end{Definition}

\begin{Remark}
    The homotopy groups of the Brauer-Picard 3-group for a fusion category $\mathcal{C}$ have the following characterization:
    \begin{enumerate}
        \item $\pi_0(\mathbf{BrPic}(\mathcal{C}))$ is the Brauer-Picard group originally defined in \cite{ENO09}, with elements given by Morita classes of invertible $\mathcal{C}$-bimodule categories. By \cite[Theorem 1.1]{ENO09}, this group is isomorphic to the group of braided auto-equivalences of $\mathcal{Z}_1(\mathcal{C})$.

        \item $\pi_1(\mathbf{BrPic}(\mathcal{C}))$ is the group of invertible objects in Drinfeld center $\mathcal{Z}_1(\mathcal{C})$.

        \item $\pi_2(\mathbf{BrPic}(\mathcal{C})) \simeq \Bbbk^\times$.
    \end{enumerate}
\end{Remark}

\begin{Definition} \label{def:Picard3Group}
    The Picard 3-group $\mathbf{Pic}(\mathcal{A})$ for a braided multifusion category $\mathcal{A}$ is the maximal sub-3-group of invertible objects in the multifusion 2-category $\mathbf{Mod}(\mathcal{A})$ of right $\mathcal{A}$-module categories, under the relative Deligne tensor product $\boxtimes_{\mathcal{A}}$.
\end{Definition}

\begin{Remark}
    The homotopy groups of the Picard 3-group for a braided fusion category $\mathcal{A}$ have the following characterization:
    \begin{enumerate}
        \item $\pi_0(\mathbf{Pic}(\mathcal{A}))$ is the Picard group originally defined in \cite{ENO09}, with elements given by Morita classes of invertible right $\mathcal{A}$-module categories. By \cite[Theorem 5.2]{ENO09}, this group is isomorphic to the group of braided auto-equivalences of $\mathcal{A}$.

        \item $\pi_1(\mathbf{Pic}(\mathcal{A}))$ is the group of invertible objects in $\mathcal{A}$.

        \item $\pi_2(\mathbf{Pic}(\mathcal{A})) \simeq \Bbbk^\times$.
    \end{enumerate}
\end{Remark}

\begin{Remark}
    For any braided multifusion category $\mathcal{A}$, we have two choices of equivalences of 2-categories between left $\mathcal{A}$-module categories and right $\mathcal{A}$-module categories, \[ \mathrm{Ind}^{\pm}: \mathbf{Lmod}(\mathcal{A}) \simeq \mathbf{Mod}(\mathcal{A}),\] induced by the braiding on $\mathcal{A}$ and its reversed braiding \cite[Section 2.8]{DN13}. Therefore, the Picard 3-group $\mathbf{Pic}(\mathcal{A})$ can be equivalently defined using left $\mathcal{A}$-module categories.
\end{Remark}

\begin{Remark}
    When $\mathcal{A}$ is equivalent to the Drinfeld center $\mathcal{Z}_1(\mathcal{C})$ of some fusion category $\mathcal{C}$, the Picard 3-group $\mathbf{Pic}(\mathcal{Z}_1(\mathcal{C}))$ is equivalent to Brauer-Picard 3-group $\mathbf{BrPic}(\mathcal{C})$.
\end{Remark}

\begin{Notation}
    For a monoidal category $\mathcal{C}$, we denote the same underlying category $\mathcal{C}$ with the reversed monoidal product as $\mathcal{C}^{1mp}$. 
    
    For a braided category $\mathcal{A}$, we denote the same underlying monoidal category $\mathcal{A}$ with the reversed braiding as $\mathcal{A}^{2mp}$.
\end{Notation}

The following lemma is informed to the author by David Reutter.

\begin{Lemma} \label{lem:BraidedAutoEquivalencesAreTheSameAsMonoidalAuto2EquivalencesOfModularCategories}
    Given a braided fusion category $\mathcal{A}$, the 2-group of braided auto-equivalences of $\mathcal{A}$, denoted as $\mathbf{Aut}_{br}(\mathcal{A})$, is equivalent to $\mathbf{Aut}_{\otimes}(\mathbf{Mod}(\mathcal{A}))$, the 3-group of monoidal auto-equivalences of the fusion 2-category $\mathbf{Mod}(\mathcal{A})$. The latter 3-group turns out to have a trivial second homotopy group.
\end{Lemma}

\begin{proof}
    We first construct a monoidal 3-functor \[\Phi:\mathbf{Aut}_{\otimes}(\mathbf{Mod}(\mathcal{A})) \to \mathbf{Aut}_{br}(\mathcal{A})\] in the following way: 
    \begin{enumerate}
        \item Given any monoidal auto-equivalence $F$ of $\mathbf{Mod}(\mathcal{A})$, take its component on the monoidal unit $\mathcal{A}$: \[\Omega F: \mathbf{Hom}(\mathcal{A},\mathcal{A}) \to \mathbf{Hom}(\mathcal{A},\mathcal{A}).\] Note that the endo-hom category is braided, \[\mathbf{Hom}(\mathcal{A},\mathcal{A}) := \mathbf{Fun}_{\mathcal{A}^{1mp}}(\mathcal{A},\mathcal{A}) \simeq \mathcal{A}; \quad P \mapsto P(I).\] On the left side, given any two right $\mathcal{A}$-module functors $P,Q: \mathcal{A} \to \mathcal{A}$, their braiding $\beta_{P,Q}: P \circ Q \to Q \circ P$ is induced by the following diagram \[\begin{tikzcd}
            {\mathcal{A}}
                \arrow[rr,"Q"]
                \arrow[ddd,equal]
            & {}
            & {\mathcal{A}}
                \arrow[rr,"P"]
            & {}
            & {\mathcal{A}}
                \arrow[ddd,equal]
            \\ {}
            & {\mathcal{A} \boxtimes_\mathcal{A} \mathcal{A}}
                \arrow[ul,"\pmb{r}_\mathcal{A}"]
                \arrow[r,"Q \boxtimes_\mathcal{A} \mathcal{A}"]
                \arrow[d,equal]
            & {\mathcal{A} \boxtimes_\mathcal{A} \mathcal{A}}
                \arrow[r,"\mathcal{A} \boxtimes_\mathcal{A} P"]
                \arrow[u,"\pmb{r}_\mathcal{A}",bend left=20pt]
                \arrow[u,"\pmb{l}_\mathcal{A}"',bend right=20pt]
            & {\mathcal{A} \boxtimes_\mathcal{A} \mathcal{A}}
                \arrow[ur,"\pmb{l}_\mathcal{A}"']
                \arrow[d,equal]
            & {}
            \\ {}
            & {\mathcal{A} \boxtimes_\mathcal{A} \mathcal{A}}
                \arrow[r,"\mathcal{A} \boxtimes_\mathcal{A} P"']
                \arrow[dl,"\pmb{l}_\mathcal{A}"']
            & {\mathcal{A} \boxtimes_\mathcal{A} \mathcal{A}}
                \arrow[r,"Q \boxtimes_\mathcal{A} \mathcal{A}"']
                \arrow[d,"\pmb{r}_\mathcal{A}",bend left=20pt]
                \arrow[d,"\pmb{l}_\mathcal{A}"',bend right=20pt]
            & {\mathcal{A} \boxtimes_\mathcal{A} \mathcal{A}}
                \arrow[dr,"\pmb{r}_\mathcal{A}"]
            & {}
            \\ {\mathcal{A}}
                \arrow[rr,"P"']
            & {}
            & {\mathcal{A}}
                \arrow[rr,"Q"']
            & {}
            & {\mathcal{A}}
        \end{tikzcd},\] where the middle square is filled by the interchanger between $P$ and $Q$, and the remaining cells are filled by naturality of left unitor $\pmb{l}$ and right unitor $\pmb{r}$, and a canonical equivalence between $\pmb{l}_\mathcal{A}$ and $\pmb{r}_\mathcal{A}$, see \cite[Lemma 2.1]{GS}.

        Let us denote the inverse of the above equivalence as \[\mathcal{A} \simeq \mathbf{Fun}_{\mathcal{A}^{1mp}}(\mathcal{A},\mathcal{A}); \quad x \mapsto x \otimes - =: \Lambda_x. \]
        
        More explicitly, $\Phi(F)$ sends any object $x$ in $\mathcal{A}$ to $F(\Lambda_x)(I)$, which can be easily extended to an endo-functor on $\mathcal{A}$. For any objects $x$ and $y$ in $\mathcal{A}$, the monoidal functor structure on $\Phi(F)$ is defined via \[\begin{tikzcd}
            {\mathcal{A}}
                \arrow[rr,"F(y \otimes -)"]
                \arrow[ddrr,"F((x \otimes y) \otimes -)"',bend right=30pt]
            & {}
            & {\mathcal{A}}
                \arrow[dd,"F(x \otimes -)"]
            \\ {}
            & {}
                \arrow[Rightarrow,ur,"F(\alpha_{x,y,-})" {xshift=5pt,yshift=-20pt},shorten <= 5pt,shorten >= 10pt]
            & {}
            \\ {}
            & {}
            & {\mathcal{A}}
        \end{tikzcd}, \quad \begin{tikzcd}[column sep=25pt]
            {\mathcal{A}}
                \arrow[rr,"\mathrm{Id}_\mathcal{A}" {name=T}, bend left=30pt]
                \arrow[rr,"F(I \otimes -)"' {name=S},bend right=30pt]
                \arrow[Rightarrow,from=S,to=T,shorten <= 7pt,shorten >= 7pt,"F(\lambda_{-})"']
            & {}
            & {\mathcal{A}}
        \end{tikzcd},\]
        \[\Phi(F)(x) \otimes \Phi(F)(y) := F(\Lambda_x)(I) \otimes F(\Lambda_y)(I) \] 
        \[\xrightarrow[\sim]{F(x \otimes -) \text{ is right module functor}} F(x \otimes (I \otimes F(y \otimes I)))\] 
        \[\xrightarrow[\sim]{\text{left unitor}} F(x \otimes F(y \otimes I))\] 
        \[\xrightarrow[\sim]{F \text{ is functorial}} F(x \otimes (y \otimes I))\]
        \[\xrightarrow[\sim]{F(\alpha^{-1}_{x,y,I})} F((x \otimes y) \otimes I) =: \Phi(F)(x \otimes y),\]
        \[I \xrightarrow[\sim]{F \text{ is functorial}} F(I) \xrightarrow[\sim]{F(\lambda^{-1}_{I})} F(I \otimes I) := \Phi(F)(I).\]
        Moreover, the braiding is preserved by $\Phi(F)$ since the monoidal 2-functor structure on $F$ preserves the interchangers and unitors in $\mathbf{Mod}(\mathcal{A})$. Lastly, as $F$ is an 2-equivalence, its induces an equivalence on $\mathbf{Hom}(\mathcal{A},\mathcal{A})$. Therefore, we obtain a braided auto-equivalence $\Phi(F)$ on $\mathcal{A}$.
        
        \item Given two monoidal auto-equivalences $F,G$ of $\mathbf{Mod}(\mathcal{A})$ and a monoidal 2-natural isomorphism $\phi:F \to G$, we take its component on the monoidal unit $\mathcal{A}$: \[\begin{tikzcd}[sep=large]
            {\mathbf{Hom}(\mathcal{A},\mathcal{A})}
                \arrow[r,"\Omega G"]
                \arrow[d,"\Omega F"']
            & {\mathbf{Hom}(G(\mathcal{A}),G(\mathcal{A}))}
                \arrow[d,"\phi^*"]
            \\ {\mathbf{Hom}(F(\mathcal{A}),F(\mathcal{A}))}
                \arrow[r,"\phi_*"']
                \arrow[Rightarrow,ur,shorten <= 25pt,shorten >= 25pt,"\Omega \phi"]
            & {\mathbf{Hom}(F(\mathcal{A}),G(\mathcal{A}))}
        \end{tikzcd},\] or under the identification $\mathcal{A} \simeq \mathbf{Fun}_{\mathcal{A}^{1mp}}(\mathcal{A},\mathcal{A})$, and recall that $F,G,\phi$ preserve the monoidal unit, the above invertible 2-cell is equivalent to \[\begin{tikzcd}[sep=large]
            {\mathcal{A}}
                \arrow[r,"\Phi(G)"]
                \arrow[d,"\Phi(F)"']
            & {\mathcal{A}}
                \arrow[d,"\mathrm{Id}_\mathcal{A}"]
            \\ {\mathcal{A}}
                \arrow[r,"\mathrm{Id}_\mathcal{A}"']
                \arrow[Rightarrow,ur,shorten <= 15pt,shorten >= 15pt,"\Phi(\phi)"]
            & {\mathcal{A}}
        \end{tikzcd},\] i.e. a natural isomorphism $\Phi(\phi):\Phi(F) \to \Phi(G)$. Lastly, this induced natural isomorphism $\Phi(\phi)$ preserves the monoidal functor structures on $F$ and $G$ since $\phi$ is a monoidal 2-natural isomorphism.
        
        \item Given two monoidal auto-equivalences $F,G$ of $\mathbf{Mod}(\mathcal{A})$, two monoidal 2-natural isomorphisms $\phi_0,\phi_1:F \to G$, any invertible monoidal modification $\zeta:\phi_0 \to \phi_1$ induces an \textit{equality} of natural isomorphisms $\Phi(\phi_0)$ and $\Phi(\phi_1)$, since $F,G,\phi_0,\phi_1,\zeta$ all need to preserve the monoidal unit.
    \end{enumerate}

    Conversely, we construct a monoidal 3-functor \[\Psi:\mathbf{Aut}_{br}(\mathcal{A}) \to \mathbf{Aut}_{\otimes}(\mathbf{Mod}(\mathcal{A})) \] as follows: 
    \begin{enumerate}
        \item First, recall that by Ostrik's theorem \cite{Ost1}, the 2-category $\mathbf{Mod}(\mathcal{A})$ is equivalent to the Morita 2-category in $\mathcal{A}$, consists of:
        \begin{enumerate}
            \item An object is a separable algebra in $\mathcal{A}$.
            
            \item A 1-morphism from $A$ to $B$ is a finite projective $(A,B)$-bimodule, and the composition is given by the relative tensor product of bimodules.
            
            \item A 2-morphism between two finite projective $(A,B)$-bimodules $M$ and $N$ is an $(A,B)$-bimodule map $f:M \to N$.
        \end{enumerate}

        Moreover, since $\mathcal{A}$ is braided, this induces a monoidal structure on the Morita 2-category in $\mathcal{A}$ as the tensor product of two separable algebras $A$ and $B$ is equipped with a canonical separable algebra on $A \otimes B$. This promote the above equivalence to an equivalence of fusion 2-categories.
        
        Given any braided auto-equivalence $H$ of $\mathcal{A}$, the desired 2-functor $\Psi(H):\mathbf{Mod}(\mathcal{A}) \to \mathbf{Mod}(\mathcal{A})$ is defined equivalently on the Morita 2-category in $\mathcal{A}$ as follows:
        \begin{enumerate}
            \item A separable algebra $A$ in $\mathcal{A}$ is sent to $H(A)$, which is also a separable algebra using the monoidal functor structure of $H$.
            
            \item A finite projective $(A,B)$-bimodule $M$ is sent to $H(M)$, whose finite projective $(H(A),H(B))$-bimodule structure is also induced by the monoidal functor structure of $H$. Moreover, given an $(A,B)$-bimodule $M$ and a $(B,C)$-bimodule $N$, one has a canonical identification of $(H(A),H(C))$-bimodules, \[H(M) \otimes_{H(B)} H(N) \simeq H(M \otimes_B N).\]
            
            \item An $(A,B)$-bimodule map $f:M \to N$. is sent to $H(f):H(M) \to H(N)$, which is an $(H(A),H(B))$-bimodule map.
        \end{enumerate}
        
        Next, the monoidal 2-functor structure on $\Psi(H)$ is induced by the the braided functor structure on $H$. More explicitly, for any separable algebras $A$ and $B$ in $\mathcal{A}$, one has canonical algebra isomorphisms \[H(A \otimes B) \simeq H(A) \otimes H(B), \quad H(I) \simeq I.\]

        \item Given two braided auto-equivalences $H,K$ of $\mathcal{A}$ and a monoidal natural isomorphism $\eta:H \to K$, the desired monoidal 2-natural isomorphism $\Psi(\eta):\Psi(H) \to \Psi(K)$ is defined equivalently on the Morita 2-category in $\mathcal{A}$ as follows:
        \begin{enumerate}
            \item A separable algebra $A$ in $\mathcal{A}$ is sent to an algebra isomorphism $\eta_A:H(A) \to K(A)$, which induces an $(H(A),K(A))$-bimodule, whose underlying object is $K(A)$ with the left $H(A)$-action given by $\eta_A$, and the right $K(A)$-action given by the identity on $K(A)$. We denote it as ${}_{\eta_A} K(A)$.
            
            \item For any finite projective $(A,B)$-bimodule $M$, $\eta_M:H(M) \to K(M)$ gives rise to an $(H(A),K(B))$-bimodule morphism \[\begin{tikzcd}[sep=large]
                {H(A)}
                    \arrow[r,"{}_{\eta_A} K(A)"]
                    \arrow[d,"H(M)"']
                & {K(A)}
                    \arrow[d,"K(M)"]
                \\ {H(B)}
                    \arrow[r,"{}_{\eta_B} K(B)"']
                    \arrow[ur,Rightarrow,shorten <= 15pt,shorten >= 15pt,"\eta_M"]
                & {K(B)}
            \end{tikzcd}.\]

            \item Given any $(A,B)$-bimodule morphism $f:M \to N$, it is compatible with $\eta_M, \eta_N$ since
           \[\begin{tikzcd}[scale cd=0.8]
                {H(A)}
                    \arrow[rr,"{}_{\eta_A} K(A)"]
                    \arrow[dd,"H(N)",bend left=30pt]
                    \arrow[dd,"H(M)"',bend right=30pt]
                & {}
                & {K(A)}
                    \arrow[dd,"K(N)"]
                \\ {\begin{matrix} {}_{H(f)} \\ \Rightarrow \end{matrix}}
                & {}
                & {}
                \\ {H(B)}
                    \arrow[rr,"{}_{\eta_B} K(B)"']
                    \arrow[uurr,Rightarrow,shorten <= 40pt,shorten >= 40pt,"\eta_N"']
                & {}
                & {K(B)}
            \end{tikzcd} = \begin{tikzcd}[scale cd=0.8]
                {H(A)}
                    \arrow[rr,"{}_{\eta_A} K(A)"]
                    \arrow[dd,"H(M)"]
                & {}
                & {K(A)}
                    \arrow[dd,"K(N)",bend left=30pt]
                    \arrow[dd,"K(M)"',bend right=30pt]
                \\ {}
                & {}
                & {\begin{matrix} {}_{K(f)} \\ \Rightarrow \end{matrix}}
                \\ {H(B)}
                    \arrow[rr,"{}_{\eta_B} K(B)"']
                    \arrow[uurr,Rightarrow,shorten <= 40pt,shorten >= 40pt,"\eta_M"]
                & {}
                & {K(B)}
            \end{tikzcd}.\]

            \item For any two separable algebras $A$ and $B$, the algebra isomorphism $\eta_{A \otimes B}$ is compatible with $\eta_A \otimes \eta_B$. The trivial algebra $I$ is assigned with the identity $\eta_I = 1_I$.
        \end{enumerate}

        \item Finally, since the second homotopy group of $\mathbf{Aut}_{br}(\mathcal{A})$ is trivial, $\Psi$ is extended trivially on this level.
    \end{enumerate}

    It is straightforward to check that $\Phi \circ \Psi$ is equivalent to the identity on $\mathbf{Aut}_{br}(\mathcal{A})$. For the other direction, $\Psi \circ \Phi$ is equivalent to the identity on $\mathbf{Aut}_\otimes(\mathbf{Mod}(\mathcal{A}))$ by the universal property of Karoubi completion. In detail, Gaiotto and Johnson-Freyd \cite{GJF} showed that the Karoubi completion of the one-point delooping 2-category $\mathrm{B} \mathcal{A}$ is equivalent to the Morita 2-category in $\mathcal{A}$, which is also equivalent to $\mathbf{Mod}(\mathcal{A})$ by Ostrik's theorem. The monoidal structure on $\mathrm{B} \mathcal{A}$, which is witnessed by the existence of the twice delooping 3-category $\mathrm{B}^2 \mathcal{A}$, gives rise to the monoidal structure on $\mathbf{Mod}(\mathcal{A})$ after Karoubi completion. Therefore, any monoidal auto-equivalence of $\mathbf{Mod}(\mathcal{A})$ is determined by its restriction to $\mathrm{B} \mathcal{A} \to \mathbf{Mod}(\mathcal{A})$, which is essentially the same as the braided auto-equivalence it induces on the monoidal unit $\mathcal{A}$.
    
    So $\mathbf{Aut}_\otimes(\mathbf{Mod}(\mathcal{A}))$ and $\mathbf{Aut}_{br}(\mathcal{A})$ are equivalent 3-groups.
\end{proof}

\begin{Definition} \label{def:InnerAutomorphismsFromPicard3Group}
    For any braided fusion category $\mathcal{A}$, we can define a 3-group morphism by sending an invertible $\mathcal{A}$-module category to the inner automorphism generated by it, \[\partial: \mathbf{Pic}(\mathcal{A}) \to \mathbf{Aut}_\otimes(\mathbf{Mod}(\mathcal{A})); \quad \mathcal{M} \mapsto (\mathcal{P} \mapsto \mathcal{M} \boxtimes_\mathcal{A} \mathcal{P} \boxtimes_\mathcal{A} \mathcal{M}^{op}).\]
\end{Definition}

\begin{Remark} \label{rmk:PostnikovTowerOfPicard3Group}
    By Lemma \ref{lem:BraidedAutoEquivalencesAreTheSameAsMonoidalAuto2EquivalencesOfModularCategories}, the 3-group morphism $\partial$ is equivalent to a 3-group morphism from $\mathbf{Pic}(\mathcal{A})$ to $\mathbf{Aut}_{br}(\mathcal{A})$, so for simplicity we shall not distinguish these two 3-group morphisms and call both of them $\partial$. It has a 1-group truncation $\partial_*:\pi_0(\mathbf{Pic}(\mathcal{A})) \to \pi_0(\mathbf{Aut}_{br}(\mathcal{A}))$ and a 2-group truncation $\partial_*:\Pi_{\leq 1}(\mathbf{Pic}(\mathcal{A})) \to \Pi_{\leq 1}(\mathbf{Aut}_{br}(\mathcal{A}))$. 
    
    In \cite{ENO09}, the two mappings above are proven to be invertible when braided fusion category $\mathcal{A}$ is \textit{non-degenerate}, i.e. a braided fusion category whose \textit{M{\"u}ger center} is $\mathcal{Z}_2(\mathcal{A}) \simeq \mathbf{Vect}$. For a general braided fusion category $\mathcal{A}$, \cite{DN13} showed that the 1-group truncation of $\partial$ is part of a crossed module structure, called the \textit{Picard crossed module} of braided fusion category $\mathcal{A}$. This is the shadow of the 3-group crossed module structure on $\partial:\mathbf{Pic}(\mathcal{A}) \to \mathbf{Aut}_\otimes(\mathbf{Mod}(\mathcal{A}))$, where we let inner automorphisms on $\mathbf{Mod}(\mathcal{A})$ act on the 3-group $\mathbf{Pic}(\mathcal{A})$ via restriction.

    In \cite{DN21}, the fiber of $\partial_*:\Pi_{\leq 1}(\mathbf{Pic}(\mathcal{A})) \to \Pi_{\leq 1}(\mathbf{Aut}_{br}(\mathcal{A}))$ is characterized by $\Pi_{\leq 1}(\mathbf{Pic}_{br}(\mathcal{A}))$, where the 3-group $\mathbf{Pic}_{br}(\mathcal{A})$, called the \textit{braided Picard group} of braided fusion category $\mathcal{A}$, is defined by invertible braided module catagories over $\mathcal{A}$ \cite[Definition 4.1]{DN21}. Using Lemma \ref{lem:BraidedAutoEquivalencesAreTheSameAsMonoidalAuto2EquivalencesOfModularCategories}, we can immediately promote $\mathbf{Pic}_{br}(\mathcal{A})$ to be the homotopy fiber of the 3-group morphism $\partial:\mathbf{Pic}(\mathcal{A}) \to \mathbf{Aut}_{br}(\mathcal{A})$. Moreover, the existence of the 3-group crossed module structure on $\partial$ hints that $\mathbf{Pic}(\mathcal{A})$ itself is a homotopy fiber, as depicted in the following diagram:
    \[\begin{tikzcd}
        {\mathbf{Pic}_{br}(\mathcal{A})}
            \arrow[r]
            \arrow[d]
        & {\mathbf{Pic}(\mathcal{A})}
            \arrow[d,"\partial"]
            \arrow[r]
        & {*}
            \arrow[d]
        \\ {*}
            \arrow[r]
        & {\mathbf{Aut}_{br}(\mathcal{A})}
            \arrow[r,"\mathbf{Z}"]
        & {\mathbf{BrPic}(\mathbf{Mod}(\mathcal{A}))}
    \end{tikzcd}.\] The 4-group morphism $\mathbf{Z}$ maps a monoidal auto-equivalence of $\mathbf{Mod}(\mathcal{A})$ to the invertible bimodule 2-category it induces over $\mathbf{Mod}(\mathcal{A})$. The rigorous definition of the \textit{Brauer-Picard 4-group} for a fusion 2-category is given in \cite[Definition 2.6]{D11}. In our case, the Brauer-Picard 4-group of $\mathbf{Mod}(\mathcal{A})$ is in general not equivalent to the one-point delooping $\mathrm{B} \, \mathbf{Pic}_{br}(\mathcal{A})$. See \cite[Example 2.10]{D11} for the description of its homotopy groups for symmetric fusion category $\mathcal{A}$.
    
    For non-degenerate braided fusion category $\mathcal{A}$, this provides the top layer of the Postnikov tower for Picard 3-group $\mathbf{Pic}(\mathcal{A})$:
    \[\begin{tikzcd}
        {\mathrm{B}^2 \Bbbk^\times}
            \arrow[r]
            \arrow[d]
        & {\mathbf{Pic}(\mathcal{A})}
            \arrow[d,"\partial"]
            \arrow[r]
        & {*}
            \arrow[d]
        \\ {*}
            \arrow[r]
        & {\mathbf{Aut}_{br}(\mathcal{A})}
            \arrow[r,"\mathbf{Z}"]
        & {\mathrm{B}^3 \Bbbk^\times}
    \end{tikzcd}.\]
\end{Remark}

\medskip

\noindent \textit{Fusion 2-Categories.} Fusion 2-categories are introduced as a categorification of fusion 1-categories by Douglas and Reutter \cite{DR}, where they are originally used to construct state sum of 4-dimensional TQFTs. The essential observation is the correct categorification of Karoubi completeness, which is later generalized to weak $n$-categories by Gaiotto and Johnson-Freyd \cite{GJF}.

In this paper, we will mainly focus on several instances of fusion 2-categories to be introduced in the following. Therefore, for the sake of clarity, we will hide the explicit definition of fusion 2-categories, but readers can always refer to \cite[Section 2.1.1]{DR} if needed. In order to keep the diagrams simple and clear, we will use the following notations and abbreviations:
\begin{enumerate}
    \item Identity $n$-morphism of an $(n-1)$-morphism $f$ (here $n=1,2$) are denoted by $1_f$. We may also omit the subscript if there is no ambiguity.

    \item Usually, the monoidal product is denoted by $\Box$ and the monoidal unit is denoted by $I$. When needed to put the product in a more compact form, we will omit the symbol $\Box$.
    
    \item We also denote the associator by $\pmb{\alpha}$, left unitor by $\pmb{l}$, right unitor by $\pmb{r}$, pentagonator by $\pmb{\pi}$, three 2-unitors by $\pmb{\lambda}$, $\pmb{\mu}$ and $\pmb{\rho}$, respectively. If in addition this fusion 2-category admits a braiding, we denote the braiding by $\pmb{b}$, and the hexagonator by $\pmb{R}$ and $\pmb{S}$. See \cite[Appendix C]{SP} for a detailed explanation of these notations.
    
    \item Superscripts will be added if we need to distinguish monoidal structures on different 2-categories. 
\end{enumerate}

\begin{Example}
    $\mathbf{2Vect}$ is the fusion 2-category of finite semisimple categories, with monoidal product given by the Deligne tensor product $\boxtimes$. By Artin-Wedderburn, every finite semisimple category is equivalent to the category of modules over some finite semisimple algebras over $\Bbbk$. Eilenberg-Watts Theorem generalizes this identification, in the sense that $\mathbf{2Vect}$ is equivalent to Morita 2-category of finite semisimple algebras over $\Bbbk$, with monoidal product given by the tensor product of algebras.
    \begin{center}
        \begin{tabular}{|c|c|}
            \hline {} & {} \\[-1.5ex]
            finite semisimple algebra $A$ & finite semisimple category $\mathbf{Mod}(A)$ \\ [1ex]
            \hline {} & {} \\[-1.5ex]
            finite projective bimodule ${}_A M_B$ & functor $\mathbf{Mod}(A) \xrightarrow{- \otimes_A M} \mathbf{Mod}(B)$ \\ [1.5ex]
            \hline {} & {} \\[-1.5ex]
            module homomorphism $f:M \to N$ & natural transformation \\ [1ex]
            \hline {} & {} \\[-1.5ex]
            $A \otimes B$ & $\mathbf{Mod}(A) \boxtimes \mathbf{Mod}(B) \simeq \mathbf{Mod}(A \otimes B)$ \\ [1ex]
            \hline
        \end{tabular}
    \end{center}
\end{Example}

\begin{Example}
    $\mathbf{2Vect}_G$ is the fusion 2-category of $G$-graded finite semisimple categories, with monoidal product given by convolution on $G$. The underlying finite semisimple 2-category is given by $\mathbf{Fun}(G,\mathbf{2Vect})$, where $G$ is viewed as a 2-groupoid with trivial 2-morphisms. The monoidal product is given by the convolution on $G$, i.e. the grading on the Deligne tensor product is given by: \[(\mathcal{C} \boxtimes \mathcal{D})_g = \bigoplus_{h \in G} \mathcal{C}_h \boxtimes \mathcal{D}_{h^{-1}g}\] and the monoidal unit is $\mathbf{Vect}$ graded on unit $e$.
    
    Generalizing Eilenberg-Watts Theorem, we have an alternative Morita theoretic interpretation of $\mathbf{2Vect}_G$:
    \begin{enumerate}
        \item Objects are $G$-graded finite semisimple algebras $A = \bigoplus_{g \in G}A_g$, i.e. on top of the graded vector space there is a multiplication $m:A \otimes A \to A$ and unit $i: \Bbbk \to A$ such that $m(A_g,A_h) \subset A_{gh}$ and $i(\Bbbk) \subset A_e$ for any $g,h \in G$. Notice that we \textit{does not} require the grading to be faithful.
        
        \item 1-morphisms are $G$-graded finite projective bimodules ${}_A M_B = \bigoplus_{g \in G}M_g$, i.e. on top of the graded vector space there are left and right actions $l:A \otimes M \to M$ and $r:M \otimes B \to M$ such that $l(A_g,M_h) \subset M_{gh}$ and $r(M_g,B_h) \subset M_{gh}$ for any $g,h \in G$.
        
        \item 2-morphisms are $G$-graded module homomorphisms $f:M \to N$, i.e. linear maps $f:M \to N$ such that $f(M_g) \subset N_g$ for any $g \in G$.
        
        \item The monoidal product is given by the convolution of graded algebras, i.e. \[(A \otimes B)_g = \bigoplus_{h \in G} A_h \otimes B_{h^{-1}g}\] and the monoidal unit is the trivial algebra $\Bbbk$ graded trivially on unit $e$.
    \end{enumerate}
\end{Example}

\begin{Example} 
    We can twist the pentagonator of $\mathbf{2Vect}_G$ with a group cocycle $\pi \in \mathrm{H}^4(G;\Bbbk^\times)$: for any $g,h,k,l \in G$, the pentagonator is replaced as follows: \[\adjustbox{scale=0.9,center}{\begin{tikzcd}
        {((\mathbf{Vect}_g \boxtimes \mathbf{Vect}_h) \boxtimes \mathbf{Vect}_k) \boxtimes \mathbf{Vect}_l} 
            \arrow[r,"\omega_{gh,k,l}"] 
            \arrow[d,"\omega_{g,h,k}1"']
        & {(\mathbf{Vect}_g \boxtimes \mathbf{Vect}_h) \boxtimes (\mathbf{Vect}_k \boxtimes \mathbf{Vect}_l)} 
            \arrow[dd,"\omega_{g,h,kl}"] 
        \\ {(\mathbf{Vect}_g \boxtimes (\mathbf{Vect}_h \boxtimes \mathbf{Vect}_k)) \boxtimes \mathbf{Vect}_l}
            \arrow[d,"\omega_{g,hk,l}"']
        & {}
        \\ {\mathbf{Vect}_g \boxtimes((\mathbf{Vect}_h \boxtimes \mathbf{Vect}_k) \boxtimes \mathbf{Vect}_l)} 
        \arrow[r,"1\omega_{h,k,l}"']
        & {\mathbf{Vect}_g \boxtimes (\mathbf{Vect}_h \boxtimes (\mathbf{Vect}_k \boxtimes \mathbf{Vect}_l))}
    \end{tikzcd}}, \] 
    \[ (1 \boxtimes \omega_{h,k,l}) \circ \omega_{g,hk,l} \circ (\omega_{g,h,k} \boxtimes 1)(((\Bbbk_g \boxtimes \Bbbk_h) \boxtimes \Bbbk_k) \boxtimes \Bbbk_l) \] \[= \pi(g,h,k,l) \cdot (\omega_{g,h,kl} \circ \omega_{gh,k,l})(((\Bbbk_g \boxtimes \Bbbk_h) \boxtimes \Bbbk_k) \boxtimes \Bbbk_l). \]
    
    We denote the new fusion 2-category as $\mathbf{2Vect}^\pi_G$.
\end{Example}

\begin{Example}
    $\mathbf{2Rep}(G):= \mathbf{Fun}(\mathrm{B} G,\mathbf{2Vect})$ is the 2-category of fully dualizable 2-representations of $G$, e.g. a finite semisimple category $\mathcal{M}$ with $G$-action given by a monoidal functor \[\rho: G \to \mathbf{End}(\mathcal{M}).\] By the universal property of Karoubi completion, we have \[\mathbf{2Rep}(G) \simeq \mathbf{Fun}(\mathrm{B} \mathbf{Vect}_G,\mathbf{2Vect}) =: \mathbf{Lmod}(\mathbf{Vect}_G),\] where 2-functors are assumed to preserve $\mathbf{Vect}$-enriched structures. It is a symmetric fusion 2-category, with monoidal product given by the Deligne tensor product $\boxtimes$ and monoidal unit given by trivial 2-representation $\mathbf{Vect}$.
\end{Example}

\begin{Example}
    We denote $\Sigma \mathbf{Rep}(G) := \mathbf{Lmod}(\mathbf{Rep}(G))$ to be the symmetric fusion 2-category of finite semisimple left module 1-categories over $\mathbf{Rep}(G)$. The monoidal structure is given by the relative Deligne tensor product $\boxtimes_{\mathbf{Rep}(G)}$, see \cite{ENO09,DSPS14}.
\end{Example}

\begin{Remark}
    A monoidal category is \textit{strict} if all the associators and unitors are identities. Mac Lane \cite{Mac63} proved that any monoidal category is monoidally equivalent to a strict monoidal category.

    However, the situation is different for monoidal 2-categories. In general, a monoidal 2-category is not monoidally equivalent to a strict monoidal 2-category, namely a monoidal 2-category whose underlying 2-category is strict and the monoidal structure is strict. Nevertheless, there is a notion of \textit{semi-strict} monoidal 2-category, where the underlying 2-category is strict, the associators and unitors are identities, but there exists non-trivial interchangers for the monoidal product of 1-morphisms, see \cite{GPS} (under the name \textit{Gray monoid}) and \cite{BN,SP,BMS,DR} for elucidation on this topic.

    Semi-strictified monoidal 2-categories are more tractable, when we define structures and prove theorems using string diagrams, as it significantly reduce the number of coherence data and conditions to be checked. However, it is worth noting that natural examples of fusion 2-categories, such as those provided above, are often not semi-strict.
\end{Remark}

\begin{Remark}
    By the famous Schur's Lemma, if there exists some non-zero morphism between two simple objects in a finite semisimple category, then these two simple objects must be isomorphic. Unfortunately, this intuition does not generalize to finite semisimple 2-categories. Instead, one has a weaken version of Schur's Lemma: suppose $x,y,z$ are simple objects in a finite semisimple 2-category $\mathfrak{C}$, then any two of $\mathbf{Hom}(x,y)$, $\mathbf{Hom}(y,z)$, $\mathbf{Hom}(x,z)$ being non-zero implies the third is non-zero. 
    
    This property provides us with an equivalence relation on the set of isomorphism classes of simple objects in $\mathfrak{C}$: we say two simple objects are \textit{connected} if there exists some non-zero 1-morphism between them. We denote the set of connected components of $\mathfrak{C}$ as $\pi_0(\mathfrak{C})$.
\end{Remark}

\begin{Remark}
    A category is \textit{skeletal} if any two isomorphic objects are equal. Any category is equivalent to a skeletal one. This provides us with a computation-friendly approach to study multifusion categories (see \cite{EGNO} for summary): 
    \begin{enumerate}
        \item We start with the skeletal data of a (multi-)fusion category $\mathcal{C}$. Since its underlying category is finite semisimple, the set of isomorphism classes of simple objects $\mathcal{O}(\mathcal{C})$ provides a canonical $\mathbb{Z}_{\geq 0}$-basis for the set of isomorphism classes of objects. One can complete this set by formally adding inverses, and thus obtains the \textit{Grothendieck group} $\mathcal{K}(\mathcal{C})$. Then the monoidal product on $\mathcal{C}$ equipped this Abelian group with a multiplication, called the \textit{(multi-)fusion ring} of $\mathcal{C}$.
        
        \item The remaining monoidal structure: associator and unitors, can be encoded as additional data on the (multi-)fusion ring $\mathcal{K}(\mathcal{C})$. The \textit{Ocneanu Rigidity} asserts that any (multi-)fusion ring only admits finitely many monoidal structures. Special cases of these skeletal data have been studied in early literature of representation theory and physics with the name $6j$-\textit{symbols} or $F$-\textit{symbols}.
    \end{enumerate}

    By analogy, a 2-category is \textit{skeletal} if any two isomorphic objects are equal, and all hom categories are skeletal. Clearly, any 2-category is equivalent to a skeletal one. However, when we would like to apply the above strategy to study fusion 2-categories, we encounter two difficulties: 
    \begin{enumerate}
        \item Monoidal structure on $\mathfrak{C}$ does not induce a fusion ring structure on the Grothendieck group of $\mathfrak{C}$. A counterexample is the fusion 2-category $\mathbf{2Rep}(\mathbb{Z}/2)$, see \cite[Example 2.1.20]{DR}. As mentioned in the previous remark, there exists non-zero 1-morphisms between two non-isomorphic simple objects in a fusion 2-category $\mathfrak{C}$. We can classify them into connected components, which form a finite set $\pi_0(\mathfrak{C})$. The discovery of Karoubi completeness for 2-categories suggests that a good notion of structure in fusion 2-categories should depend on $\pi_0(\mathfrak{C})$ rather than $\mathcal{K}(\mathfrak{C})$. However, in general fusion rules of $\mathfrak{C}$ is not compatible with connectedness of simple objects, hence it does not descend to a multiplication on the Abelian group generated by $\pi_0(\mathfrak{C})$. It would be an interesting question to find the correct notion of \textit{fusion 2-ring for fusion 2-categories}.
        
        \item When we fix a representative in every isomorphism class of simple objects, the additional monoidal structure on $\mathfrak{C}$: associators, 1-unitors, pentagonators and 2-unitors, all descend to certain cohomology-type data on the set of isomorphism classes of simple objects and the set of isomorphism classes of simple 1-morphisms between simple objects. This data is called the $10j$-\textit{symbols} of $\mathfrak{C}$ in \cite[Definition 3.3.8]{DR}. $10j$-symbols are the key ingredients when we define state sum on 4-dimensional manifolds using fusion 2-categories.
    \end{enumerate}
\end{Remark}

\noindent {\it Drinfeld Center of Fusion 2-Categories.} The notion of Drinfeld center is generalized to monoidal 2-categories in \cite{BN,Cr}. Décoppet showed that the Drinfeld center of a fusion 2-category is a braided fusion 2-category, and Morita equivalent fusion 2-categories have braided equivalent Drinfeld centers \cite{D9}. In this paper, we will mainly focus on the following concrete examples of Drinfeld centers, hence we will omit the explicit definition of Drinfeld center of fusion 2-categories.

\begin{Example}
    Drinfeld center of $\Sigma \mathbf{Rep}(G)$, denoted by $\mathscr{Z}_1(\Sigma \mathbf{Rep}(G))$, has an explicit characterization in \cite{DN21} as the 2-category of \textit{braided module categories} over $\mathbf{Rep}(G)$, and later \cite{JFR} provides an alternative description as the 2-category of \textit{half-braided algebras and bimodules} in $\mathbf{Rep}(G)$. In a previous work of the author joint with Décoppet \cite{DX}, we showed that $\mathscr{Z}_1(\Sigma \mathbf{Rep}(G))$ is equivalent to the 2-category of \textit{local modules} over $\mathbf{Rep}(G)$, where $\mathbf{Rep}(G)$ is viewed as a braided algebra in $\mathbf{2Vect}$.
\end{Example}

\begin{Example} \label{exmp:DrinfeldCenterOf2VectG}
    Fusion 2-categories $\mathbf{2Vect}_G$ and $\mathbf{2Rep}(G)$ are Morita equivalent in the sense of \cite[Theorem 5.4.3]{D8}, thus their Drinfeld centers are equivalent braided fusion 2-categories. By the calculation from \cite{KTZ}, the two Drinfeld centers $\mathscr{Z}_1(\mathbf{2Vect}_G)$ and $\mathscr{Z}_1(\mathbf{2Rep}(G))$ are both identified with the 2-category of $G$-\textit{equivariant} $G$-\textit{graded finite semisimple categories}: \[\mathscr{Z}_1(\mathbf{2Vect}_G) = \mathscr{Z}_1(\mathbf{2Rep}(G)) = \mathbf{Fun}(\mathrm{LB} G, \mathbf{2Vect}) \simeq \bigoplus_{g \in \mathrm{Cl}(G)} \mathbf{2Rep}(Z(g)).\]
\end{Example}

\begin{Example} \label{exmp:DrinfeldCenterOfTwisted2VectG}
    After twisting the convolution product with a 4-cocycle $\pi$ on $G$, the Drinfeld center of fusion 2-category $\mathbf{2Vect}^\pi_G$ is identified with the 2-category of $G$-\textit{equivariant} $G$-\textit{graded finite semisimple categories} with the monoidal product twisted by $\pi$ in \cite{KTZ}. Its underlying 2-category has the decomposition \[\mathscr{Z}_1(\mathbf{2Vect}^\pi_G) = \mathbf{Fun}(\mathrm{L} \widetilde{\mathrm{B} G}^\pi, \mathbf{2Vect}) \simeq \bigoplus_{g \in \mathrm{Cl}(G)} \mathbf{2Rep}(Z(g),\tau_g(\pi)),\] where $\widetilde{\mathrm{B} G}^\pi$ is the extension of $\mathrm{B} G$ by $\pi$: \[\begin{tikzcd}
        {\mathrm{B}^3 \Bbbk^\times}
            \arrow[r]
            \arrow[d]
        & {\widetilde{\mathrm{B} G}^\pi}
            \arrow[d]
            \arrow[r]
        & {*}
            \arrow[d]
        \\ {*}
            \arrow[r]
        & {\mathrm{B} G}
            \arrow[r,"\pi"]
        & {\mathrm{B}^4 \Bbbk^\times}
    \end{tikzcd},\] equivalently, $\widetilde{\mathrm{B} G}^\pi$ is one-point delooping of a 3-group constructed from extension of $G$ by $\mathrm{B}^2 \Bbbk^\times$. Its linearization is $\mathbf{2Vect}_G^\pi$. Then $\mathrm{L} \widetilde{\mathrm{B} G}^\pi$ is the free looping 3-groupoid of $\widetilde{\mathrm{B} G}^\pi$. Any element $g \in G$ gives a transgression $\tau_g: \mathrm{H}^4(G;\Bbbk^\times) \to \mathrm{H}^3(Z(g);\Bbbk^\times)$; the twist $\tau_g(\pi)$ is the transgression of $\pi$ under $\tau_g$. If $g$ and $g'$ are conjugate in $G$, then the twists $\tau_g(\pi)$ and $\tau_{g'}(\pi)$ are equivalent.
\end{Example}

\medskip

\noindent {\it Algebras.} We are going to recall the notion of algebras in fusion 2-categories. It first appeared with the name \textit{pseudo-monoid} in the context of monoidal 2-categories in \cite{DS}. They can be viewed as a externalization of the notion of monoidal categories from $\mathbf{2Vect}$ to an arbitrary fusion 2-category.

\begin{Definition} \label{def:algebra}
    An algebra in a semi-strict monoidal 2-category $\mathfrak{C}$ consists of the following data:
    \begin{enumerate}
        \item An object $A$ in $\mathfrak{C}$;
        
        \item Two 1-morphisms $m:A \, \Box \, A \to A$ and $i:I \to A$;
        
        \item Invertible 2-morphisms \begin{center}
            \begin{tabular}{@{}c c c@{}}
            $\begin{tikzcd}[sep=small]
            A \arrow[rrrr, equal] \arrow[rrdd, "i1"'] &  & {} \arrow[dd, Rightarrow, "\lambda"', near start, shorten > = 1ex] &  & A \\
                                               &  &                           &  &   \\
                                               &  & AA \arrow[rruu, "m"']     &  &  
            \end{tikzcd},$
            
            &
            
            $\begin{tikzcd}[sep=small]
            AAA \arrow[dd, "1m"'] \arrow[rr, "m1"]    &  & AA \arrow[dd, "m"] \\
                                                        &  &                      \\
            AA \arrow[rr, "m"'] \arrow[rruu, Rightarrow, "\mu", shorten > = 2.5ex, shorten < = 2.5ex] &  & A                   
            \end{tikzcd},$
            
            &
            
            $\begin{tikzcd}[sep=small]
                                              &  & AA \arrow[rrdd, "m"] \arrow[dd, Rightarrow, "\rho", shorten > = 1ex, shorten < = 2ex] &  &   \\
                                              &  &                                             &  &   \\
            A \arrow[rruu, "1i"] \arrow[rrrr,equal] &  & {}                                          &  & A
            \end{tikzcd},$

            \end{tabular}
            \end{center} satisfying the coherence conditions (1)(2)(3)(4) in \cite[Definition 1.2.1]{D7}.
    \end{enumerate}
\end{Definition}

\begin{Definition}
    A braided algebra in semi-strict braided monoidal 2-category $\mathfrak{B}$ consists of the following data:
    \begin{enumerate}
        \item An algebra $(B,m,i,\mu,\lambda,\rho)$ in $\mathfrak{B}$;
        
        \item An invertible 2-morphism \[\begin{tikzcd}[sep=small]
            {AA} \arrow[rrrr, "m"] \arrow[rrdd, "\pmb{b}_{A,A}"'] &  & {} &  & {A} \\
                                               &  &                           &  &   \\
                                               &  & AA \arrow[rruu, "m"'] \arrow[uu, Rightarrow, "\beta"', shorten < = 1ex]     &  &  
            \end{tikzcd}\] satisfying the coherence conditions (15)(16)(17) in \cite[Definition 3.1]{DY22}.
    \end{enumerate}
\end{Definition}

\begin{Definition}
    A rigid algebra $A$ in a monoidal 2-category $\mathfrak{C}$ is an algebra whose multiplication $m:A \, \Box \, A \to A$ admits an $(A,A)$-bimodule right adjoint $m^*:A \to A \, \Box \, A$. For a detailed list of the data and coherence conditions of a rigid algebra, see \cite[Definition 2.1.1]{D7}.
\end{Definition}

\begin{Definition}
    A separable algebra $A$ in a monoidal 2-category $\mathfrak{C}$ is a rigid algebra whose counit $\epsilon^m: m \circ m^* \to 1_A$ admits an $(A,A)$-bimodule section $\gamma^m: 1_A \to m \circ m^*$. Coherence conditions of a separable algebra can be found in \cite[Definition 2.1.7]{D7}.
\end{Definition}

\begin{Proposition}[{\cite[Corollary 5.1.2]{D9}}]
    Any rigid algebra in a multifusion 2-category is automatically separable.
\end{Proposition}

Following \cite{DMNO}, we call a separable braided algebra in a braided 2-category $\mathfrak{B}$ as an \textbf{{\'e}tale algebra}. Here are some examples of {\'e}tale algebras.

\begin{Example}
    In $\mathfrak{B} = \mathbf{2Vect}$, an {\'e}tale algebra is just a braided multifusion category.
\end{Example}

Let $G$ be a finite group.

\begin{Example}
    In $\mathfrak{B} = \mathbf{2Rep}(G)$, an {\'e}tale algebra is a braided multifusion category $\mathcal{M}$ with a $G$-action, i.e. a monoidal functor $G \to \mathbf{Aut}_{br}(\mathcal{M})$.
\end{Example}

\begin{Example}
    In $\mathfrak{B} = \Sigma \mathbf{Rep}(G)$, an {\'e}tale algebra is a braided multifusion category $\mathcal{M}$ with a $\mathbf{Rep}(G)$-action, i.e. a braided functor $\mathbf{Rep}(G) \to \mathcal{Z}_2(\mathcal{M})$.
\end{Example}

\begin{Example} \label{exmp:BraidedMultiFusionCatOverRepG}
    In $\mathfrak{B} = \mathscr{Z}_1(\Sigma \mathbf{Rep}(G))$, by \cite[Lemma 3.2.1]{DX}, an {\'e}tale algebra consists of a braided multifusion category $\mathcal{M}$ with a braided functor $\mathbf{Rep}(G) \to \mathcal{M}$.
\end{Example}

\begin{Example} \label{exmp:BraidedCrossedMultiFusionCat}
    In $\mathfrak{B} = \mathscr{Z}_1(\mathbf{2Rep}(G))$, an {\'e}tale algebra is a $G$-\textit{crossed braided multifusion category} $\mathcal{M}$, see \cite{Tur00,Tur08,Tur10,Mu4,ENO09,DGNO,Cui16,JPR}.
\end{Example}

\begin{Example} \label{exmp:BraidedTwistedCrossedMultiFusionCat}
    In $\mathfrak{B} = \mathscr{Z}_1(\mathbf{2Vect}^\pi_G)$ where twisting $\pi$ is given as a group 4-cocycle on $G$, an {\'e}tale algebra is a $G$-crossed braided multifusion category $\mathcal{M}$ which trivializes twisting $\pi$. We will provide more details in Section \ref{sec:TwistedCrossedBraidedFusionCat}.
\end{Example}

Finally, we introduce the notion of connected algebras in fusion 2-categories, which generalizes the condition for multifusion categories to be fusion.
\begin{Definition}
    An algebra $A$ in a (linear) monoidal 2-category $\mathfrak{C}$ is connected if its unit $i: I \to A$ is simple in the linear category $\mathbf{Hom}(I,A)$.
\end{Definition}
In Section \ref{sec:BraidedFusionCatOverRepG}, Section \ref{sec:CrossedBraidedFusionCat} and Section \ref{sec:TwistedCrossedBraidedFusionCat}, we will describe connectedness for {\'e}tale algebras in $\mathscr{Z}_1(\Sigma \mathbf{Rep}(G))$, $\mathscr{Z}_1(\mathbf{2Rep}(G))$ and $\mathscr{Z}_1(\mathbf{2Vect}^\pi_G)$, respectively.

\medskip

\noindent {\it Modules.} Now we recall the notion of modules in monoidal 2-categories, following \cite{D4,D7}. They are externalizations of the notion of module category over a monoidal category to an arbitrary monoidal 2-category.

\begin{Definition}
    Given an algebra $(A,m,i,\mu,\lambda,\rho)$ in a semi-strict monoidal 2-category $\mathfrak{C}$, a right module over $A$ consists of the following data:
    \begin{enumerate}
        \item An object $M$ in $\mathfrak{C}$;
        
        \item A 1-morphism $n^M:M \, \Box \, A \to M$ in $\mathfrak{C}$;
        
        \item Invertible 2-morphisms \begin{center}
            \begin{tabular}{@{}c c@{}}
            $\begin{tikzcd}[sep=small]
            MAA \arrow[dd, "1m"'] \arrow[rr, "n^M1"]    &  & MA \arrow[dd, "n^M"] \\
                                                        &  &                      \\
            MA \arrow[rr, "n^M"'] \arrow[rruu, Rightarrow, "\nu^M", shorten > = 2.5ex, shorten < = 2.5ex] &  & M                   
            \end{tikzcd},$
            
            &
            
            $\begin{tikzcd}[sep=small]
                                              &  & MA \arrow[rrdd, "n^M"] \arrow[dd, Rightarrow, "\rho^M", shorten > = 1ex, shorten < = 2ex] &  &   \\
                                              &  &                                             &  &   \\
            M \arrow[rruu, "1i"] \arrow[rrrr,equal] &  & {}                                          &  & M
            \end{tikzcd},$

            \end{tabular}
            \end{center} satisfying the coherence conditions (5)(6) in \cite[Definition 1.2.3]{D7}.
    \end{enumerate}

    Similarly, a left module over $A$ consists of:
    \begin{enumerate}
        \item An object $N$ in $\mathfrak{C}$;
        
        \item A 1-morphism $l^N:A \, \Box \, N \to N$ in $\mathfrak{C}$;
        
        \item Invertible 2-morphisms \begin{center}
            \begin{tabular}{@{}c c@{}}
            $\begin{tikzcd}[sep=small]
            AAN \arrow[dd, "m1"'] \arrow[rr, "1l^N"]    &  & AN \arrow[dd, "l^N"] \\
                                                        &  &                      \\
            AN \arrow[rr, "l^N"'] \arrow[rruu, Rightarrow, "\kappa^N", shorten > = 2.5ex, shorten < = 2.5ex] &  & N                   
            \end{tikzcd},$
            
            &
            
            $\begin{tikzcd}[sep=small]
                N \arrow[rrrr, equal] \arrow[rrdd, "i1"'] &  & {} \arrow[dd, Rightarrow, "\lambda^N"', near start, shorten > = 1ex] &  & N \\
                                                   &  &                           &  &   \\
                                                   &  & AN \arrow[rruu, "l^N"']     &  &  
                \end{tikzcd},$

            \end{tabular}
            \end{center} satisfying the coherence conditions (10)(11) in \cite[Definition 1.3.1]{D7}.
    \end{enumerate}
\end{Definition}

\begin{Definition}
    Take an algebra $(A,m,i,\mu,\lambda,\rho)$ in a semi-strict monoidal 2-category $\mathfrak{C}$. Given two right $A$-modules $(M,n^M,\nu^M,\rho^M)$ and $(N,n^N,\nu^N,\rho^N)$, a right $A$-module 1-morphism between them consists of a 1-morphism $f:M \to N$ in $\mathfrak{C}$ and an invertible 2-morphism \[\begin{tikzcd}[sep=small]
        {MA} 
            \arrow[dd, "f1"'] 
            \arrow[rr, "n^M"] 
        & {}
        & {M} 
            \arrow[dd, "f"] 
            \
        \\ {}
        & {}  
        & {}
        \\ {NA} 
            \arrow[rr, "n^N"']
            \arrow[rruu, Rightarrow, "\psi^f", shorten > = 3ex, shorten < = 2.7ex] 
        & {}
        & {N}
    \end{tikzcd}\] satisfying the coherence conditions (8)(9) in \cite[Definition 1.2.6]{D7}.

    Given two right $A$-module 1-morphisms $(f,\psi^f)$ and $(g,\psi^g)$ between $M$ and $N$, a right $A$-module 2-morphism between them consists of a 2-morphism $\eta:f \to g$ in $\mathfrak{C}$ such that the coherence condition in \cite[Definition 1.2.7]{D7} holds.

    Similarly, given two left $A$-modules $(M,l^M,\kappa^M,\lambda^M)$ and $(N,l^N,\kappa^N,\lambda^N)$, a left $A$-module 1-morphism between them consists of a 1-morphism $f:M \to N$ in $\mathfrak{C}$ and an invertible 2-morphism \[\begin{tikzcd}[sep=small]
        {A M} 
            \arrow[dd, "1f"'] 
            \arrow[rr, "l^M"] 
        & {}
        & {M} 
            \arrow[dd, "f"] 
            \
        \\ {}
        & {}  
        & {}
        \\ {A N} 
            \arrow[rr, "l^N"']
            \arrow[rruu, Rightarrow, "\chi^f", shorten > = 3ex, shorten < = 2.7ex] 
        & {}
        & {N}
    \end{tikzcd}\] satisfying the coherence conditions (12)(13) in \cite[Definition 1.3.2]{D7}.

    Given two left $A$-module 1-morphisms $(f,\chi^f)$ and $(g,\chi^g)$ between $M$ and $N$, a left $A$-module 2-morphism between them consists of a 2-morphism $\eta:f \to g$ in $\mathfrak{C}$ such that the coherence condition in \cite[Definition 1.3.3]{D7} holds.
\end{Definition}

\begin{Lemma}
    For any algebra $A$ in a semi-strict monoidal 2-category $\mathfrak{C}$, right $A$-modules, right $A$-module 1-morphisms and right $A$-module 2-morphisms in $\mathfrak{C}$ form a 2-category, which we denote as $\mathbf{Mod}_\mathfrak{C}(A)$. Similarly, left $A$-modules, left $A$-module 1-morphisms and left $A$-module 2-morphisms in $\mathfrak{C}$ form a 2-category, which we denote as $\mathbf{Lmod}_\mathfrak{C}(A)$.
\end{Lemma}

\begin{Proposition}[{\cite[Proposition 3.1.2]{D7}}]
    Let $A$ be a separable algebra in a multifusion 2-category $\mathfrak{C}$. Then the 2-category of right $A$-modules $\mathbf{Mod}_\mathfrak{C}(A)$ and the 2-category of left $A$-modules $\mathbf{Lmod}_\mathfrak{C}(A)$ are finite semisimple.
\end{Proposition}

\begin{Definition}
    Given two algebras $A$ and $B$ in a semi-strict monoidal 2-category $\mathfrak{C}$, an $(A,B)$-bimodule consists of the following data:
    \begin{enumerate}
        \item An object $M$ in $\mathfrak{C}$;
        
        \item A left $A$-module structure $(M,l^M,\kappa^M,\lambda^M)$;
        
        \item A right $B$-module structure $(M,n^M,\nu^M,\rho^M)$;
        
        \item An invertible 2-morphism \[\begin{tikzcd}[sep=small]
            {AMB} 
                \arrow[rr, "l^M1"] 
                \arrow[dd, "1n^M"'] 
            & {}
            & {MB} 
                \arrow[dd, "n^M"] 
                \
            \\ {}
            & {}  
            & {}
            \\ {AM} 
                \arrow[rr, "l^M"']
                \arrow[rruu, Rightarrow, "\mu^M", shorten > = 3ex, shorten < = 2.7ex] 
            & {}
            & {M}
        \end{tikzcd}\] satisfying the coherence conditions (38)(39)(40)(41) in \cite[Definition 2.4.6]{Xu24}.
    \end{enumerate}

    Given two $(A,B)$-bimodules $M$ and $N$, an $(A,B)$-bimodule 1-morphism between them consists of:
    \begin{enumerate}
        \item A 1-morphism $f:M \to N$ in $\mathfrak{C}$;
        
        \item A left $A$-module 1-morphism structure $(f,\chi^f)$ and a right $B$-module 1-morphism structure $(f,\psi^f)$, satisfying the coherence condition (42) in \cite[Definition 2.4.7]{Xu24}.
    \end{enumerate}

    Given two $(A,B)$-bimodule 1-morphisms $f,g:M \to N$, an $(A,B)$-bimodule 2-morphism between them consists of a 2-morphism $\eta:f \to g$ in $\mathfrak{C}$ such that $\eta$ satisfies the coherence conditions for left $A$-module 2-morphisms and right $B$-module 2-morphisms.
\end{Definition}

\begin{Lemma}
    For any algebras $A$ and $B$ in a semi-strict monoidal 2-category $\mathfrak{C}$, $(A,B)$-bimodules, $(A,B)$-bimodule 1-morphisms and 2-morphisms in $\mathfrak{C}$ form a 2-category, which we denote as $\mathbf{Bimod}_\mathfrak{C}(A,B)$. Moreover, when $A=B$, we can simply denote $\mathbf{Bimod}_\mathfrak{C}(A,B)$ as $\mathbf{Bimod}_\mathfrak{C}(A)$.
\end{Lemma}

\begin{Proposition}[{\cite[Proposition 3.1.3]{D7}, \cite[Theorem 3.2.8]{D8}}]
    Let $A$ and $B$ be separable algebras in a multifusion 2-category $\mathfrak{C}$. Then the 2-category of $(A,B)$-bimodules $\mathbf{Bimod}_\mathfrak{C}(A,B)$ is finite semisimple. Moreover, the 2-category of $(A,A)$-bimodules $\mathbf{Bimod}_\mathfrak{C}(A)$ is equipped with a multifusion structure, induced by the \textit{relative tensor product}, see \cite[Definition 3.1.3]{D8}.
\end{Proposition}

\begin{Definition} \label{def:LocalModule}
    Let $B$ be a braided algebra in a semi-strict braided monoidal 2-category $\mathfrak{B}$. A local $B$-module consists of the following data:
    \begin{enumerate}
        \item An object $M$ in $\mathfrak{B}$;
        
        \item A right $B$-module structure $(M,n^M,\nu^M,\rho^M)$;
        
        \item An invertible 2-morphism $h^M$:
        \[\begin{tikzcd}
            {BM} 
                \arrow[rr, "\pmb{b}_{B,M}"] 
            & {}
                \arrow[dd, Rightarrow, "h^M", shorten > = 3ex, shorten < = 3ex] 
            & {MB} 
                \arrow[dd, "n^M"] 
                \
            \\ {}
            & {}  
            & {}
            \\ {MB} 
                \arrow[rr, "n^M"']
                \arrow[uu, "\pmb{b}_{M,B}"] 
            & {}
            & {M}
        \end{tikzcd},\] which satisfies the coherence conditions (21)(22)(23) in \cite[Definition 2.1.1]{DX}.
    \end{enumerate}
    
    Given two local $B$-modules $M$ and $N$, a local $B$-module 1-morphism between them consists of a 1-morphism $f:M \to N$ in $\mathfrak{B}$ satisfying the coherence condition (24) in \cite[Definition 2.1.2]{DX}.

    Given two local $B$-module 1-morphisms $f,g:M \to N$, a local $B$-module 2-morphism between them is simply a right $B$-module 2-morphism $\eta:f \to g$.
\end{Definition}

\begin{Lemma}
    For any braided algebra $B$ in a semi-strict braided monoidal 2-category $\mathfrak{B}$, local $B$-modules, local $B$-module 1-morphisms and 2-morphisms in $\mathfrak{B}$ form a 2-category, which we denote as $\mathbf{Mod}^{loc}_\mathfrak{B}(B)$.
\end{Lemma}

\begin{Proposition}[{\cite[Theorem 2.3.4]{DX}}]
    Let $B$ be an {\'e}tale algebra in a braided multifusion 2-category $\mathfrak{B}$. Then the 2-category of local $B$-modules $\mathbf{Mod}^{loc}_\mathfrak{B}(B)$ is finite semisimple, and it is equipped with a braided multifusion structure, induced by the relative tensor product over $B$.
\end{Proposition}

Following \cite{DMNO}, we categorify the notion of Lagrangian algebra via the characterization of vanishing local modules.

\begin{Definition}
    A Lagrangian algebra $L$ in a braided fusion 2-category $\mathfrak{B}$ is a connected {\'e}tale algebra whose 2-category of local modules is trivial: \[\mathbf{Mod}^{loc}_\mathfrak{B}(L) \simeq \mathbf{2Vect}.\]
\end{Definition}

\begin{Example}[{\cite[Example 3.3.2]{DX}}]
    In $\mathfrak{B} = \mathbf{2Vect}$, a Lagrangian algebra is a non-degenerate braided fusion category, i.e. a braided fusion category whose \textit{M{\"u}ger center} is equivalent to $\mathbf{Vect}$.
\end{Example}

In Section \ref{sec:BraidedFusionCatOverRepG}, Section \ref{sec:CrossedBraidedFusionCat} and Section \ref{sec:TwistedCrossedBraidedFusionCat}, we will describe Lagrangian algebras in $\mathscr{Z}_1(\Sigma \mathbf{Rep}(G))$, $\mathscr{Z}_1(\mathbf{2Rep}(G))$ and $\mathscr{Z}_1(\mathbf{2Vect}^\pi_G)$, respectively.

\medskip

\noindent {\it Equivariantization.} Fusion categories $\mathbf{Vect}_G$ and $\mathbf{Rep}(G)$ are Morita equivalent, with $\mathbf{Vect}$ as an invertible bimodule \cite{Ost2}. In other word, $\mathbf{2Rep}(G)$ and $\Sigma \mathbf{Rep}(G)$ are equivalent as finite semisimple 2-categories. Moreover, this equivalence preserves symmetric monoidal structures on both sides via a concrete construction named \textbf{equivariantization} and \textbf{de-equivariantization}, see \cite{AG,Bru,Gai05,Kir,Mu0} for references during the development of this subject. In the following, we quote results from \cite[Chapter 4]{DGNO} for a complete overview.

\begin{Definition}
    Let $\mathcal{C}$ be a finite semisimple category with a $G$-action, i.e. a monoidal functor $\rho: G \to \mathbf{Aut}(\mathcal{C})$. Then its equivariantization is defined as \[\mathcal{C}^G := \mathbf{Vect} \boxtimes_{\mathbf{Vect}_G} \mathcal{C}, \] which inherits a canonical left $\mathbf{Rep}(G)$-action from $\mathbf{Vect}$.
\end{Definition}

\begin{Remark}
    The above definition can be unpacked and rewritten as follows: 
    \begin{enumerate}
        \item An object of $\mathcal{C}^G$ is a pair $(V,\{\omega_g\}_{g \in G})$, where $V$ is an object in $\mathcal{C}$ and $\omega_g:V \to \rho(g)(V)$ is an isomorphism for each $g \in G$, such that for any $g,h \in G$, the following diagram commutes:
        \[\begin{tikzcd}[sep=large]
            {V}
                \arrow[r,"\omega_g"]
                \arrow[d,"\omega_{gh}"']
            & {\rho(g)(V)}
                \arrow[d,"\rho(g)(\omega_h)"]
            \\ {\rho(gh)(V)}
            & {\rho(g)(\rho(h)(V))}
                \arrow[l,"\rho_{g,h}(V)"]
        \end{tikzcd}.\]
        
        \item A morphism from $(V,\{\omega^V_g\}_{g \in G})$ to $(W,\{\omega^W_g\}_{g \in G})$ consists of a morphism $f:V \to W$ in $\mathcal{C}$ such that for each $g \in G$, the following diagram commutes:
        \[\begin{tikzcd}[sep=large]
            {V}
                \arrow[r,"f"]
                \arrow[d,"\omega^V_g"']
            & {W}
                \arrow[d,"\omega^W_g"]
            \\ {\rho(g)(V)}
                \arrow[r,"\rho(g)(f)"']
            & {\rho(g)(W)}
        \end{tikzcd}.\]
    \end{enumerate}
\end{Remark}

\begin{Definition}
    Let $\mathcal{D}$ be a finite semisimple category with a $\mathbf{Rep}(G)$-action, i.e. a monoidal functor $F: \mathbf{Rep}(G) \to \mathbf{End}(\mathcal{D})$. Then its de-equivariantization is defined as \[\mathcal{D}_G := \mathbf{Vect} \boxtimes_{\mathbf{Rep}(G)} \mathcal{D}, \] which inherits a canonical $G$-action from $\mathbf{Vect}$.
\end{Definition}

\begin{Remark}
    The above definition is equivalent to the following:
    \begin{enumerate}
        \item Recall that $\mathbf{Fun}(G)$ is a separable algebra in $\mathbf{Rep}(G)$. An object of $\mathcal{D}_G$ is a pair $(V,\theta)$ where $V$ is an object in $\mathcal{D}$ and $\theta:\mathbf{Fun}(G) \odot V \to V$ equips $V$ with a left $\mathbf{Fun}(G)$-module structure, i.e. the following diagrams commute:
        \[\begin{tikzcd}[sep=large]
            {(\mathbf{Fun}(G) \otimes \mathbf{Fun}(G)) \odot V}
                \arrow[d,"m \odot V"']
                \arrow[r,"\cong"]
            & {\mathbf{Fun}(G) \odot (\mathbf{Fun}(G) \odot V)}
                \arrow[d,"\mathbf{Fun}(G) \odot \theta"]
            \\ {\mathbf{Fun}(G) \odot V}
                \arrow[d,"\theta"']
            & {\mathbf{Fun}(G) \odot V}
                \arrow[d,"\theta"]
            \\ {V}
                \arrow[r,equal]
            & {V}
        \end{tikzcd},\]
        \[\begin{tikzcd}[sep=large]
            {V}
                \arrow[d,equal]
                \arrow[r,"\cong"]
            & {\Bbbk \odot V}
                \arrow[d,"i \odot V"]
            \\ {V}
            & {\mathbf{Fun}(G) \odot V}
                \arrow[l,"\theta"]
        \end{tikzcd}.\]
        \item A morphism from $(V,\theta^V)$ to $(W,\theta^W)$ consists of a morphism $f:V \to W$ in $\mathcal{D}$ such that the following diagram commutes:
        \[\begin{tikzcd}[column sep=50pt,row sep=30pt]
            {\mathbf{Fun}(G) \odot V}
                \arrow[r,"\mathbf{Fun}(G) \odot f"]
                \arrow[d,"\theta^V"']
            & {\mathbf{Fun}(G) \odot W}
                \arrow[d,"\theta^W"]
            \\ {V}
                \arrow[r,"f"']
            & {W}
        \end{tikzcd}.\]
    \end{enumerate}
\end{Remark}

\begin{Example}
    Consider $\mathbf{Vect}$ with the trivial $G$-action. We have $\mathbf{Vect}^G \simeq \mathbf{Rep}(G)$ as regular module over $\mathbf{Rep}(G)$. Conversely, $\mathbf{Rep}(G)_G \simeq \mathbf{Vect}$.
\end{Example}

\begin{Example}
    Consider $\mathbf{Vect}_G$ with the left $G$-action via translation. Then we have $(\mathbf{Vect}_G)^G \simeq \mathbf{Vect}$ as $\mathbf{Rep}(G)$-module categories. Conversely, $\mathbf{Vect}_G$ can be identified with the de-equivariantization of $\mathbf{Vect}$ with trivial $\mathbf{Rep}(G)$-action.
\end{Example}

\begin{Example}
    Take any subgroup $H$ in $G$. Then $\mathbf{Vect}_{G/H}$ has the left $G$-action via translation. We have $(\mathbf{Vect}_{G/H})^G \simeq \mathbf{Rep}(H)$ as $\mathbf{Rep}(G)$-module categories.
\end{Example}

\begin{Example}
    Take any 2-cocycle $\mu \in \mathrm{H}^2(G;\Bbbk^\times)$, and consider the tensor functor $\mathbf{Vect}_G \to \mathbf{Vect}$ twisted by $\mu$. Then the equivariantization of $\mathbf{Vect}$ with this twisted $G$-action is equivalent to the category of $\mu$-projective $G$-representations $\mathbf{Rep}(G,\mu)$ with the canonical $\mathbf{Rep}(G)$-action.
\end{Example}

\begin{Example}
    Consider $\mathbf{Vect}_G$ with $G$ acting via conjugation. Then we have $(\mathbf{Vect}_G)^G \simeq \mathcal{Z}_1(\mathbf{Vect}_G) \simeq \mathcal{Z}_1(\mathbf{Rep}(G))$, with $\mathbf{Rep}(G)$ embedding into its Drinfeld center. More generally, we can twist $\mathbf{Vect}_G$ by a 3-cocycle $\omega \in \mathrm{H}^3(G;\Bbbk^\times)$, and the equivariantization would be $(\mathbf{Vect}^\omega_G)^G \simeq \mathcal{Z}_1(\mathbf{Vect}^\omega_G)$, which still has $\mathbf{Rep}(G)$ embedding into it.\footnote{This is also famously known as the minimal modular extension of $\mathbf{Rep}(G)$, see \cite{LKW}.}
\end{Example}

By the universal property of the relative tensor product, we can extend the equivariantization and de-equivariantization to an equivalence of 2-categories. Moreover, this equivalence preserves symmetric monoidal structures on both sides, namely for any finite semisimple categories $\mathcal{M}$ and $\mathcal{N}$ with $G$-actions, we have a canonical equivalence of $\mathbf{Rep}(G)$-module categories: \[\mathcal{M}^G \boxtimes_{\mathbf{Rep}(G)} \mathcal{N}^G \simeq (\mathcal{M} \boxtimes \mathcal{N})^G.\]

In summary, there is an equivalence of \textit{symmetric fusion $2$-categories}:
$$\begin{tikzcd}
{\{\text{module categories of }\mathbf{Vect}_G\text{ with }\boxtimes\text{ as monoidal product}\}}
    \arrow[d,"\text{equivariantization}",shift left=5pt]
\\ {\{\text{module categories of }\mathbf{Rep}(G)\text{ with }\boxtimes_{\mathbf{Rep}(G)}\text{ as monoidal product}\}}
    \arrow[u,"\text{de-equivariantization}",shift left=5pt]
\end{tikzcd}$$
As a consequence, there are one-to-one correspondences between various types of algebras in both 2-categories.

\medskip

Separable algebras in $\mathbf{2Rep}(G)$ are multifusion categories with $G$-actions, e.g. a multifusion category $\mathcal{C}$ with a monoidal functor $G \to \mathbf{Aut}_\otimes(\mathcal{C})$.

On the other hand, separable algebras in $\Sigma \mathbf{Rep}(G)$ are $\mathbf{Rep}(G)$-module multifusion categories, e.g. a multifusion category $\mathcal{D}$ with a braided functor $\mathbf{Rep}(G) \to \mathcal{Z}_1(\mathcal{D})$. It is connected if and only if the underlying multifusion category $\mathcal{D}$ is fusion.

\begin{Proposition}
    Equivariantization and de-equivariantization give a one-to-one correspondence between multifusion categories with $G$-actions and $\mathbf{Rep}(G)$-module multifusion categories. 
    
    In particular, this provides a one-to-one correspondence between fusion categories with $G$-actions and fusion categories with embeddings of $\mathbf{Rep}(G)$ into their Drinfeld centers.
\end{Proposition}

Similarly, Étale algebras in $\mathbf{2Rep}(G)$ are braided multifusion categories with $G$-actions, e.g. a braided multifusion category $\mathcal{A}$ with a monoidal functor $G \to \mathbf{Aut}_{br}(\mathcal{A})$.

On the other hand, étale algebras in $\Sigma \mathbf{Rep}(G)$ are $\mathbf{Rep}(G)$-module braided multifusion categories, e.g. a braided multifusion category $\mathcal{B}$ with a braided functor $\mathbf{Rep}(G) \to \mathcal{Z}_2(\mathcal{B})$. It is connected if and only if the underlying braided multifusion category $\mathcal{B}$ is fusion.

\begin{Proposition}
    Equivariantization and de-equivariantization give a one-to-one correspondence between braided multifusion categories with $G$-actions and $\mathbf{Rep}(G)$-module braided multifusion categories. 
    
    In particular, this provides a one-to-one correspondence between braided fusion categories with $G$-actions and braided fusion categories with embeddings of $\mathbf{Rep}(G)$ into their M{\"u}ger centers.
\end{Proposition}

Furthermore, any equivalence between $\mathbf{2Rep}(G)$ and $\Sigma \mathbf{Rep}(G)$ as symmetric fusion 2-categories induces an equivalence of braided fusion 2-categories between $\mathscr{Z}_1(\mathbf{2Rep}(G))$ and $\mathscr{Z}_1(\Sigma \mathbf{Rep}(G))$.

\begin{Proposition}
    Equivariantization and de-equivariantization give a one-to-one correspondence between $G$-crossed finite semisimple categories and braided $\mathbf{Rep}(G)$-module categories \cite{DN21}.
\end{Proposition}

\begin{Proposition} \label{prop:EquivariantizationOfCrossedBraidedFusionCat}
    Equivariantization and de-equivariantization give a one-to-one correspondence between $G$-crossed braided multifusion categories and braided multifusion categories over $\mathbf{Rep}(G)$, e.g. a braided multifusion category $\mathcal{D}$ with a braided functor $\mathbf{Rep}(G) \to \mathcal{D}$.

    In particular, this provides a one-to-one correspondence between $G$-crossed braided fusion categories and braided fusion categories with braided embeddings of $\mathbf{Rep}(G)$.
\end{Proposition}

The above ingredients could be simply categorified to construct equivariantization and de-equivariantization for fusion 2-categories. Although the details are missing in the literature, we can still provide a rough sketch here.

Just as the 1-categorical case, it is true by \cite[Example 5.4.5]{D8} that $\mathbf{2Vect}$ is an invertible bimodule 2-category between fusion 2-categories $\mathbf{2Vect}_G$ and $\mathbf{2Rep}(G)$. Let us denote the followings:
\begin{itemize}
    \item $\mathbf{3Rep}(G):= \mathbf{Lmod}(\mathbf{2Vect}_G)$ is the 3-category of finite semisimple 2-categories with $G$-actions, with symmetric monoidal product induced by absolute 2-Deligne tensor product $\boxtimes$;
    
    \item $\Sigma \mathbf{2Rep}(G) := \mathbf{Lmod}(\mathbf{2Rep}(G))$ is the 3-category of finite semisimple $\mathbf{2Rep}(G)$-module 2-categories, with symmetric monoidal product induced by relative 2-Deligne tensor product $\boxtimes_{\mathbf{2Rep}(G)}$, see \cite{D10} for details.
\end{itemize}

We shall have the following equivalence of symmetric fusion 3-categories:
$$\begin{tikzcd}
{\{\text{module 2-categories of }\mathbf{2Vect}_G\text{ with }\boxtimes\text{ as monoidal product}\}}
    \arrow[d,"\text{equivariantization}",shift left=5pt]
\\ {\{\text{module 2-categories of }\mathbf{2Rep}(G)\text{ with }\boxtimes_{\mathbf{2Rep}(G)}\text{ as monoidal product}\}}
    \arrow[u,"\text{de-equivariantization}",shift left=5pt]
\end{tikzcd}$$

By Lemma \ref{lem:BraidedAutoEquivalencesAreTheSameAsMonoidalAuto2EquivalencesOfModularCategories}, one has $\mathbf{Aut}_{\otimes}(\mathbf{Mod}(A)) \simeq \mathbf{Aut}_{br}(\mathcal{A})$, thus a $G$-action on braided fusion category $\mathcal{A}$ is equivalent to a $G$-action on fusion 2-category $\mathbf{Mod}(\mathcal{A})$.

\begin{Conjecture}
    Equivariantization and delooping commute with each other. In other words, for a multifusion category $\mathcal{C}$ with $G$-action and a braided multifusion category $\mathcal{A}$ with $G$-action, we expect that there exists an equivalence of module 2-categories over $\mathbf{2Rep}(G)$: \[ \mathbf{Mod}(\mathcal{C})^G \simeq \mathbf{Mod}(\mathcal{C}^G), \] and an equivalence of $\mathbf{2Rep}(G)$-module multifusion 2-categories: \[ \mathbf{Mod}(\mathcal{A})^G \simeq \mathbf{Mod}(\mathcal{A}^G). \]
\end{Conjecture}

\noindent {\it Crossed Drinfeld Centers.} Equivariantization and Drinfeld center are known to \textit{not commute}. To fix this, we need to replace the Drinfeld center of a multifusion category with $G$-action by a certain \textit{crossed extension}. 

\begin{Proposition}
    Let $\mathcal{C}$ be a fusion category with $G$-action. The Drinfeld center of $\mathcal{C}^G$ contains $\mathbf{Rep}(G)$ as a braided subcategory. The de-equivariantization of $\mathcal{Z}_1(\mathcal{C}^G)$ is a $G$-crossed braided fusion category, whose identity component is $\mathcal{Z}_1(\mathcal{C})$. By \cite{ENO09}, this is equivalent to a morphism $G \to \mathbf{Pic}(\mathcal{Z}_1(\mathcal{C})) \simeq \mathbf{BrPic}(\mathcal{C})$.
\end{Proposition}

This is the so-called \textit{crossed Drinfeld center}. We can take this property as the definition of crossed Drinfeld center.

\begin{Definition} \label{def:CrossedDrinfeld1Center}
    Let $\mathcal{C}$ be a fusion category with $G$-action. Then the crossed Drinfeld center of $\mathcal{C}$ is defined to be $\mathcal{Z}_1(\mathcal{C}^G)_G$.
\end{Definition}

\begin{Remark}
    In \cite{GNN}, the authors define the crossed Drinfeld center for a $G$-graded fusion category by taking the \textit{relative Drinfeld center}, which first appears in \cite{Maj91}. Then from a fusion category $\mathcal{D}$ with $G$-action, we can construct a \textit{crossed product category} $\mathcal{C}:= \mathcal{D} \rtimes G$ \cite{Tam,Nik08}, which is Morita equivalent to the equivariantization $\mathcal{D}^G$. On the other hand, one can take the crossed Drinfeld center $\mathcal{Z}_\mathcal{D}(\mathcal{C})$ of $\mathcal{C}$, which is a $G$-crossed braided fusion category. Finally, by Theorem 3.5 in \cite{GNN}, we have the following equivalence of braided fusion categories containing $\mathbf{Rep}(G)$: \[\mathcal{Z}_1(\mathcal{D}^G) \simeq \mathcal{Z}_\mathcal{D}(\mathcal{C})^G.\]
\end{Remark}

\begin{Remark}
    In \cite{GNN}, the main focus is on the center of $G$-graded fusion categories, rather than fusion categories with $G$-actions. The former can be seen as separable algebras in $\mathbf{2Vect}_G$ while the latter are separable algebras in $\mathbf{2Rep}(G)$. Nevertheless, since the two fusion 2-categories are Morita equivalent, they have the same center $\mathscr{Z}_1(\mathbf{2Vect}_G) \simeq \mathscr{Z}_1(\mathbf{2Rep}(G))$. The notion of crossed Drinfeld center for a $G$-graded fusion category (or respectively for a fusion category with $G$-action) can be understood as taking the categorified \textit{full center} \cite{KR,Dav10,DKR}. Thus the resulting crossed Drinfeld centers are Lagrangian algebras in $\mathscr{Z}_1(\mathbf{2Vect}_G)$ (or respectively in $\mathscr{Z}_1(\mathbf{2Rep}(G))$). In particular, they are $G$-crossed braided fusion categories.
\end{Remark}

Lastly, we can categorify the above construction to fusion 2-categories. 

\begin{Definition} \label{def:CrossedDrinfeld2Center}
    Let $\mathfrak{C}$ be a fusion 2-category with $G$-action, i.e. there is a monoidal 2-functor $G \to \mathbf{Aut}_\otimes(\mathfrak{C})$. Then the crossed Drinfeld center of $\mathfrak{C}$ is defined as the de-equivariantization of the Drinfeld center of the equivariantization of $\mathfrak{C}$, \[\mathscr{Z}^G_1(\mathfrak{C}) := \mathscr{Z}_1(\mathfrak{C}^G)_G \simeq \mathbf{2Vect} \boxtimes_{\mathbf{2Rep}(G)} \mathscr{Z}_1(\mathbf{2Vect} \boxtimes_{\mathbf{2Vect}_G} \mathfrak{C}).\]
\end{Definition}

\begin{Remark}
    By definition, for a fusion 2-category $\mathfrak{C}$ with $G$-action, $\mathscr{Z}^G_1(\mathfrak{C})$ is a $G$-crossed braided fusion 2-category, whose identity component is $\mathscr{Z}_1(\mathfrak{C})$. Extension theory for fusion 2-categories has been developed in \cite{D11}. Therefore, a crossed Drinfeld center of a fusion 2-category $\mathfrak{C}$ with $G$-action is equivalent to a morphism $G \to \mathbf{Pic}(\mathscr{Z}_1(\mathfrak{C})) \simeq \mathbf{BrPic}(\mathfrak{C})$, where the target is now the (Brauer-)Picard 4-group for fusion 2-category $\mathfrak{C}$, with multiplication given by the relative 2-Deligne's tensor product $\boxtimes_\mathfrak{C}$ \cite{D10}.
\end{Remark}

\begin{Proposition} \label{prop:CrossedDrinfeld2CenterOfModA}
    Suppose $\mathfrak{C} = \mathbf{Mod}(\mathcal{A})$ for some braided fusion category $\mathcal{A}$ with $G$-action. By Lemma \ref{lem:BraidedAutoEquivalencesAreTheSameAsMonoidalAuto2EquivalencesOfModularCategories}, fusion 2-category $\mathfrak{C}$ is equipped with a $G$-action. Then its crossed Drinfeld center $\mathscr{Z}^G_1(\mathbf{Mod}(\mathcal{A}))$ consists of the following data:
    \begin{itemize}
        \item An object is given by $\mathcal{M} = \bigoplus_{g \in G} \mathcal{M}_g$ such that for each $g \in G$, $\mathcal{M}_g$ is a finite semisimple right $\mathcal{A}$-module category; moreover, it is equipped with a $G$-\textbf{braided action}, i.e. a $G$-graded natural isomorphism \[(\hbar^\mathcal{M}_g)_{M,X}:M \otimes^\mathcal{M} (g^{-1} \odot X) \to M \otimes^\mathcal{M} X\] for $M$ in $\mathcal{M}_g$ and $X$ in $\mathcal{A}$, satisfying the coherence conditions (\ref{eqn:crossedholonomytriangle}) (\ref{eqn:crossedholonomyassociatorI}) (\ref{eqn:crossedholonomyassociatorII})\footnote{Here we follows Definition \ref{def:crossedholonomy} with $\pi$ trivial and $h=e$, but we drop the $G$-crossed structures on all functors and natural transforms.}.
        
        \item A 1-morphism consists of components $F_g: \mathcal{M}_g \to \mathcal{N}_g$ where $F_g$ is a right $\mathcal{A}$-module functor which preserves the $G$-braided actions.
        
        \item A 2-morphism consists of components $\eta_g: F_g \to G_g$ where $\eta_g$ is a just right $\mathcal{A}$-module natural transform.
    \end{itemize}
\end{Proposition}

\begin{proof}
    Notice that the identity component of $\mathscr{Z}_1^G(\mathbf{Mod}(\mathcal{A}))$ is $\mathscr{Z}_1(\mathbf{Mod}(\mathcal{A}))$, which is equivalent to the 2-category of finite semisimple braided module categories over $\mathcal{A}$ \cite[Theorem 4.11]{DN21}. More generally, the characterization of $G$-braided actions and their coherence conditions can be written down by carefully unpacking the definition of crossed Drinfeld center.
\end{proof}

\begin{Remark} \label{rmk:CrossedDrinfeld2CenterOfModA4Morphism}
    Denote the action on $\mathcal{A}$ as $\gamma:\mathcal{A} \to \mathbf{Aut}_{br}(\mathcal{A})$, then the crossed Drinfeld center $\mathscr{Z}^G_1(\mathbf{Mod}(\mathcal{A}))$ corresponds to the 4-group morphism given by the composition \[G \xrightarrow{\gamma} \mathbf{Aut}_{br}(\mathcal{A}) \xrightarrow{\mathbf{Z}} \mathbf{BrPic}(\mathbf{Mod}(\mathcal{A})),\] where $\mathbf{Z}:\mathbf{Aut}_{br}(\mathcal{A}) \to \mathbf{BrPic}(\mathbf{Mod}(\mathcal{A}))$ is the canonical map from the braided automorphism 2-group of $\mathcal{A}$ to the Brauer-Picard 4-group of $\mathbf{Mod}(\mathcal{A})$ introduced in Remark \ref{rmk:PostnikovTowerOfPicard3Group}.
\end{Remark}

\section{Braided Fusion Extensions of Tannakian Categories} \label{sec:BraidedFusionCatOverRepG}

    Recall that braided multifusion categories are {\'e}tale algebras in $\mathbf{2Vect}$. By Example \ref{exmp:BraidedMultiFusionCatOverRepG}, we know that braided multifusion categories over $\mathbf{Rep}(G)$ are {\'e}tale algebras in $\mathscr{Z}_1(\Sigma \mathbf{Rep}(G))$. In this section, we will classify connected and Lagrangian {\'e}tale algebras in $\mathscr{Z}_1(\Sigma \mathbf{Rep}(G))$.

    \begin{Proposition} \label{prop:ConnectedAndLagrangianEtaleAlgebrasInLocalModules}
        Suppose $A$ is an {\'e}tale algebra in a braided multifusion 2-category $\mathfrak{B}$, then one has the following:
        \begin{enumerate}
            \item An {\'e}tale algebra in $\mathbf{Mod}^{loc}_\mathfrak{B}(A)$ is equivalent to an {\'e}tale algebra $B$ in $\mathfrak{B}$ together with a braided algebra 1-morphism $f:A \to B$. 
            
            \item If $A$ is connected, then $B$ is connected as an {\'e}tale algebra in $\mathfrak{B}$ if and only if it is connected as an {\'e}tale algebra in $\mathbf{Mod}^{loc}_\mathfrak{B}(A)$.
            
            \item There is a braided equivalence of local modules: \[\mathbf{Mod}^{loc}_\mathfrak{B}(B) \simeq \mathbf{Mod}^{loc}_{\mathbf{Mod}^{loc}_\mathfrak{B}(A)}(B).\]
            
            \item In particular, $B$ is Lagrangian as a connected {\'e}tale algebra in $\mathfrak{B}$ if and only if it is Lagrangian as a connected {\'e}tale algebra in $\mathbf{Mod}^{loc}_\mathfrak{B}(A)$.
        \end{enumerate}
    \end{Proposition}

    \begin{proof}
        (1) is proven in \cite[Lemma 3.2.1]{DX}. (3) is proven in \cite[Proposition 3.2.2]{DX}. (4) is a direct consequence of (2) and (3). For (2), we recall that there is a 2-natural adjunction of free 2-functor and forgetful 2-functor between $\mathbf{Mod}_\mathfrak{B}(A)$ and $\mathfrak{B}$, which induces the equivalence of hom categories: \[ \mathbf{Hom}_\mathfrak{B}(I,B) \simeq \mathbf{Hom}_{\mathbf{Mod}_\mathfrak{B}(A)}(A,B).\] Under this equivalence, the unit of $B$ is mapped to $f$. Meanwhile, by Definition \ref{def:LocalModule} the following forgetful 1-functor is fully faithful: \[\mathbf{Hom}_{\mathbf{Mod}^{loc}_\mathfrak{B}(A)}(A,B) \to \mathbf{Hom}_{\mathbf{Mod}_\mathfrak{B}(A)}(A,B).\] Hence, the unit of $B$ is simple in $\mathbf{Hom}_\mathfrak{B}(I,B)$ if and only if $f$ is simple in $\mathbf{Hom}_{\mathbf{Mod}^{loc}_\mathfrak{B}(A)}(A,B)$.
    \end{proof}

    When $\mathfrak{B} = \mathbf{2Vect}$ and $A = \mathbf{Rep}(G)$, we obtain the following characterization of connected and Lagrangian {\'e}tale algebras in $\mathscr{Z}_1(\Sigma \mathbf{Rep}(G))$.

    \begin{Corollary} \label{cor:ConnectedAndLagrangianEtaleAlgebrasInZ(SigmaRepG)}
        An {\'e}tale algebra in $\mathscr{Z}_1(\Sigma \mathbf{Rep}(G))$ consisting of a braided multifusion category $\mathcal{M}$ and a braided functor $F:\mathbf{Rep}(G) \to \mathcal{M}$. Then this {\'e}tale algebra is connected if and only if $\mathcal{M}$ is fusion. 
        
        The braided fusion 2-category of local modules is equivalent to $\mathscr{Z}_1(\mathbf{Mod}(\mathcal{M}))$ in this case. Moreover, this {\'e}tale algebra is Lagrangian if and only if $\mathcal{M}$ is non-degenerate.
    \end{Corollary}

    An arbitrary braided functor $F:\mathbf{Rep}(G) \to \mathcal{M}$ can be uniquely decomposed into a \textit{surjective} (aka dominant) one followed by an \textit{injective} (aka faithful) one. More specifically, one can find some subgroup $H$ in $G$ (determined only up to conjugacy) such that the essential image of $F$ in $\mathcal{M}$ is a symmetric subcategory which is equivalent to $\mathbf{Rep}(H)$. In other word, a connected {\'e}tale algebra in $\mathscr{Z}_1(\Sigma \mathbf{Rep}(G))$ is determined by a subgroup $H$ and a braided fusion extension $\iota:\mathbf{Rep}(H) \hookrightarrow \mathcal{M}$, while a Lagrangian algebra is determined by a subgroup $H$ and a non-degenerate extension $\iota:\mathbf{Rep}(H) \hookrightarrow \mathcal{M}$.

    \begin{Example}
        A non-degenerate braided extension $\iota:\mathbf{Rep}(H) \hookrightarrow \mathcal{M}$ is called \textit{minimal} if its M{\"u}ger centralizer satisfies $\mathcal{Z}_2(\mathbf{Rep}(H);\mathcal{M}) \simeq \mathbf{Rep}(G)$, see \cite{LKW,KLWZZ,JFR} for more details. Moreover, these minimal extensions form an Abelian group under a seemingly-complicated product. They turns out to be equivalent to the braided auto-equivalences on the braided fusion 2-category $\mathscr{Z}_1(\Sigma \mathbf{Rep}(H))$, which are classified by the group cohomology $\mathrm{H}^3(H;\Bbbk^\times)$, see \cite{KLWZZ,JFR}. More explicitly, any group 3-cocycle $\omega$ on $H$ gives rise to a minimal extension $\iota:\mathbf{Rep}(H) \hookrightarrow \mathcal{Z}_1(\mathbf{Vect}^\omega_H)$. They are instances of Lagrangian algebras in $\mathscr{Z}_1(\Sigma \mathbf{Rep}(G))$ by the above corollary.
    \end{Example}

    \begin{Remark} \label{rmk:RightModulesOverBraidedFusionCategoryOverRepG}
        Let's consider the 2-category of right modules over a connected {\'e}tale algebra $B$ in $\mathscr{Z}_1(\Sigma \mathbf{Rep}(G))$.
        \begin{itemize}
            \item When $B$ is provided by a minimal non-degenerate extension \[\mathbf{Rep}(G) \hookrightarrow \mathcal{Z}_1(\mathbf{Vect}^\omega_G)\] for some $\omega \in \mathrm{H}^3(H;\Bbbk^\times)$, the 2-category of right $B$-modules is still equivalent to $\Sigma \mathbf{Rep}(G)$ as a fusion 2-category.
            
            \item When $B$ is provided by a general non-degenerate extension \[\mathbf{Rep}(G) \hookrightarrow \mathcal{M},\] the 2-category of right $B$-modules is \[\mathbf{Mod}_{\mathscr{Z}_1(\Sigma \mathbf{Rep}(G))}(B) \simeq \mathbf{Mod}(\mathcal{Z}_2(\mathbf{Rep}(G);\mathcal{M})).\]
            
            \item When $B$ is a general connected {\'e}tale algebra,\[F: \mathbf{Rep}(G) \twoheadrightarrow \mathbf{Rep}(H) \hookrightarrow \mathcal{D},\] the 2-category of right $B$-modules is equivalent to the Drinfeld centralizer \cite[Definition 3.1.1]{Xu24}: \[\mathbf{Mod}_{\mathscr{Z}_1(\Sigma \mathbf{Rep}(G))}(B) \simeq  \mathscr{Z}_1(\Sigma \mathbf{Rep}(G) \xrightarrow{\Sigma F} \Sigma \mathcal{D}),\] where monoidal 2-functor $\Sigma F$ is induced by braided functor $F$ via the universal property of Karoubi completion \cite{GJF}. Notice that in general this 2-category is not connected, as opposed to the above two special cases.
        \end{itemize}
    \end{Remark}
    
\section{Crossed Braided Multifusion Categories} \label{sec:CrossedBraidedFusionCat}

    In this section, we will classify connected and Lagrangian {\'e}tale algebras in $\mathscr{Z}_1(\mathbf{2Rep}(G))$ via equivariantization and de-equivariantization. We first recall the definition of crossed braided multifusion categories.

    \begin{Definition}
        A $G$-\textbf{crossed braided} multifusion category consists of 
        \begin{enumerate}
            \item An underlying multifusion category $(\mathcal{C},\otimes,I,\alpha,\lambda,\rho)$;
            
            \item A $G$-grading $\mathcal{C} = \bigoplus_{g \in G} \mathcal{C}_g$, which is not assume to be \textit{faithful};
            
            \item This monoidal structure is compatible with the $G$-grading, i.e. for any $g,h,k \in G$, we have $\mathcal{C}_g \otimes \mathcal{C}_h \subseteq \mathcal{C}_{gh}$ and the unit object $I$ in $\mathcal{C}_e$;
        
            \item A $G$-action $\odot:G \to \mathbf{Aut}^\otimes(\mathcal{C})$ such that the grading is compatible with the action: $\forall g, h \in G, g \odot \mathcal{C}_h \subseteq \mathcal{C}_{ghg^{-1}}$. 
            
            \item A $G$-crossed braiding $\beta_{X,Y}:X \otimes Y \to (g \odot Y) \otimes X$ defined for any $g,h \in G$, $X$ in $\mathcal{C}_g$ and $Y$ in $\mathcal{C}_h$ and extended naturally to all objects in $\mathcal{C}$;
            
            \item This $G$-crossed braiding is compatible with the $G$-action such that conditions (81)(82)(83) in \cite[Definition 4.41]{DGNO} are satisfied.
        \end{enumerate}
    \end{Definition}

    Applying the extension theory of Etingof, Nikshych and Ostrik \cite{ENO09}, we can characterize $G$-crossed braided \textit{fusion} categories via the following proposition.

    \begin{Proposition}[{\cite[Theorem 7.12]{ENO09}}] \label{prop:CrossedBraidedExtensionAndPicard3Group}
        A $G$-crossed extension of a braided fusion category $\mathcal{A}$ is equivalent to a 3-group morphism $G \to \mathbf{Pic}(\mathcal{A})$, where $\mathbf{Pic}(\mathcal{A})$ is the Picard 3-group of $\mathcal{A}$ we defined in Definition \ref{def:Picard3Group}.
    \end{Proposition}

    By Example \ref{exmp:BraidedCrossedMultiFusionCat}, we see that {\'e}tale algebras in $\mathscr{Z}_1(\mathbf{2Rep}(G))$ are exactly $G$-crossed braided \textit{multifusion} categories. It is easy to see that $G$-crossed braided \textit{fusion} categories give rise to connected {\'e}tale algebras in $\mathscr{Z}_1(\mathbf{2Rep}(G))$. However, the converse is not true in general. To classify connected and Lagrangian {\'e}tale algebras in $\mathscr{Z}_1(\mathbf{2Rep}(G))$, we may utilize the braided equivalence between $\mathscr{Z}_1(\mathbf{2Rep}(G))$ and $\mathscr{Z}_1(\Sigma \mathbf{Rep}(G))$ via equivariantization (see Proposition \ref{prop:EquivariantizationOfCrossedBraidedFusionCat}).

    \begin{Proposition}[{\cite[Proposition 4.56]{DGNO}}] \label{prop:DeequivariantizationOfBradedFusionCatOverRepG}
        Let $\mathcal{D}$ be a braided fusion category with an embedding $\mathbf{Rep}(G) \hookrightarrow \mathcal{D}$. Its de-equivariantization $\mathcal{C} := \mathcal{D}_G$ is a $G$-crossed braided fusion category, and conversely $\mathcal{D}$ is equivalent to the equivariantization $\mathcal{C}^G$. Let $\mathcal{A}$ be the identity component of $\mathcal{C}$, so it is a braided fusion category. Then we have the follows:
        \begin{enumerate}
            \item $\mathcal{A}^G \simeq \mathcal{Z}_2(\mathbf{Rep}(G);\mathcal{D})$.
            
            \item $\mathcal{D}$ is non-degenerate if and only if $\mathcal{A}$ is non-degenerate and $\mathcal{C}$ is faithfully graded by $G$.
            
            \item $\mathcal{Z}_2(\mathcal{A})^G \simeq \mathcal{Z}_2(\mathcal{D}) \lor \mathbf{Rep}(G)$, where $\mathcal{Z}_2(\mathcal{D}) \lor \mathbf{Rep}(G)$ denotes the smallest fusion subcategory of $\mathcal{D}$ containing both $\mathcal{Z}_2(\mathcal{D})$ and $\mathbf{Rep}(G)$.
            
            \item The support of $\mathcal{C}$ forms a normal subgroup $N$ of $G$, and $\mathbf{Rep}(G/N) \simeq \mathcal{Z}_2(\mathcal{D}) \cap \mathbf{Rep}(G)$.
        \end{enumerate}
    \end{Proposition}

    Combine Proposition \ref{prop:DeequivariantizationOfBradedFusionCatOverRepG} with Proposition \ref{prop:CrossedBraidedExtensionAndPicard3Group}, we obtain the following classification of connected and Lagrangian {\'e}tale algebras in $\mathscr{Z}_1(\mathbf{2Rep}(G))$.

    \begin{Theorem} \label{thm:ClassificationOfConnectedCrossedBraidedMultiFusionCategories}
        A connected {\'e}tale algebra in $\mathscr{Z}_1(\mathbf{2Rep}(G))$ is determined by
        \begin{enumerate}
            \item A subgroup $H$ of $G$ up to conjugacy, 
            
            \item A normal subgroup $N$ of $H$, 
            
            \item A braided fusion category $\mathcal{A}$ with an $H$-action $\rho: H \to \mathbf{Aut}_{br}(\mathcal{A})$,
            
            \item A 4-group morphism $\Gamma: H/N \to \mathbf{BrPic}(\mathbf{Mod}(\mathcal{A}))$.
            
            \item An invertible 2-morphism: \[\begin{tikzcd}
            {H}
                \arrow[d,twoheadrightarrow]
                \arrow[r,"\rho"]
            & {\mathbf{Aut}_{br}(\mathcal{A})}
                \arrow[d,"\mathbf{Z}"]
            \\ {H/N}
                \arrow[r,"\Gamma"']
                \arrow[ur,Rightarrow,shorten >=3ex,shorten <=4ex]
            & {\mathbf{BrPic}(\mathbf{Mod}(\mathcal{A}))}
        \end{tikzcd}, \] where the right arrow is explained in Remark \ref{rmk:PostnikovTowerOfPicard3Group}. The homotopy fiber of this above diagram produces an $N$-crossed extension $\alpha: N \to \mathbf{Pic}(\mathcal{A})$. 
    \end{enumerate} 

        The corresponding $G$-crossed braided multifusion category is given by \[\mathcal{C} := \mathbf{Fun}_H \left( G,\bigoplus_{g \in N} \alpha(g) \right).\] Its 2-category of local modules is equivalent to $\mathscr{Z}_1(\mathbf{Mod}(\mathcal{C}^G))$. Moreover, this {\'e}tale algebra is Lagrangian if and only if $\mathcal{A}$ is non-degenerate and $N=H$.
    \end{Theorem}

    \begin{proof}
        For a connected {\'e}tale algebra $\mathcal{C}$ in $\mathscr{Z}_1(\mathbf{2Rep}(G))$, we first consider its equivariantization $\mathcal{D} := \mathcal{C}^G$, as a connected {\'e}tale algebra in $\mathscr{Z}_1(\Sigma \mathbf{Rep}(G))$. From Section \ref{sec:BraidedFusionCatOverRepG}, we know $\mathcal{D}$ is determined by a subgroup $H$ of $G$ up to conjugacy and a braided fusion extension $\iota:\mathbf{Rep}(H) \hookrightarrow \mathcal{D}$. After equivariantization, we see that $\mathcal{C}$ is induced by an $H$-crossed braided fusion category $\mathcal{D}_H$: \[\mathcal{C} \simeq \mathcal{D}_G \simeq \mathbf{Fun}_H(G,\mathcal{D}_H).\]
        
        Next, by Proposition \ref{prop:DeequivariantizationOfBradedFusionCatOverRepG} we get a complete understanding of the structure of $H$-crossed braided fusion category $\mathcal{D}_H$. The largest fusion subcategory of $\mathcal{D}$ containing both $\mathbf{Rep}(H)$ and $\mathcal{Z}_2(\mathcal{D})$, denoted by $\mathcal{Z}_2(\mathcal{D}) \cap \mathbf{Rep}(H)$, is a fusion subcategory of a Tannakian category, so it is also Tannakian with $\mathcal{Z}_2(\mathcal{D}) \cap \mathbf{Rep}(H) \simeq \mathbf{Rep}(H/N)$ for a normal subgroup $N$ of $H$, which turns out to be the support of $\mathcal{D}_H$ before equivariantization.

        The identity component of $\mathcal{D}_H$ is a braided fusion category $\mathcal{A}$, which is equivalent to the de-equivariantization of $\mathcal{Z}_2(\mathbf{Rep}(H);\mathcal{D})$. By Proposition \ref{prop:CrossedBraidedExtensionAndPicard3Group}, as a faithful $N$-crossed extension of $\mathcal{A}$, $\mathcal{D}_H$ is equivalent to a 3-group morphism $\alpha:N \to \mathbf{Pic}(\mathcal{A})$. Lastly, the $H$-action on $\mathcal{D}_H$ restricts to an $H$-action on $\mathcal{A}$ preserving the braiding, so we obtain a monoidal functor $\rho:H \to \mathbf{Aut}_{br}(\mathcal{A})$. This $N$-crossed extension is compatible with the $H$-action on $\mathcal{A}$, so the following square can be filled by an invertible 2-morphism of 3-groups: \[\begin{tikzcd}
            {N}
                \arrow[d,hook]
                \arrow[r,"\alpha"]
            & {\mathbf{Pic}(\mathcal{A})}
                \arrow[d,"\partial"]
            \\ {H}
                \arrow[r,"\rho"']
            & {\mathbf{Aut}_{br}(\mathcal{A})}
        \end{tikzcd},\] where $\partial$ is defined in Definition \ref{def:InnerAutomorphismsFromPicard3Group}. Moreover, $\alpha$ is a morphism between homotopy fibers (see Remark \ref{rmk:PostnikovTowerOfPicard3Group}), and compatibility implies that itself arises as a homotopy fiber of the following invertible 2-morphism of 4-groups: \[ \begin{tikzcd}
            {H}
                \arrow[d,twoheadrightarrow]
                \arrow[r,"\rho"]
            & {\mathbf{Aut}_{br}(\mathcal{A})}
                \arrow[d,"\mathbf{Z}"]
            \\ {H/N}
                \arrow[r,"\Gamma"']
                \arrow[ur,Rightarrow,shorten >=3ex,shorten <=4ex]
            & {\mathbf{BrPic}(\mathbf{Mod}(\mathcal{A}))}
        \end{tikzcd}, \] where $\Gamma$ corresponds to the 2-category of local modules of $\mathcal{C}$, viewed as an $H/N$-crossed extension of the braided fusion 2-category\footnote{See \cite{D11} for the extension theory of fusion 2-categories.} $\mathscr{Z}_1(\mathbf{Mod}(\mathcal{A}))$.
        
        More explicitly, consider the following fusion subcategories of $\mathcal{D}$: \[\begin{tikzcd}
            {\mathcal{Z}_2(\mathcal{D}) \cap \mathbf{Rep}(H)}
                \arrow[r,hook]
                \arrow[d,hook]
            & {\mathcal{Z}_2(\mathcal{D})}
                \arrow[d,hook]
                \arrow[ddr,hook,bend left=20pt]
            & {}
            \\{\mathbf{Rep}(H)}
                \arrow[r,hook]
                \arrow[drr,hook,bend right=15pt]
            & {\mathcal{Z}_2(\mathcal{D}) \lor \mathbf{Rep}(H)}
                \arrow[dr,hook]
            & {}
            \\ {}
            & {}
            & {\mathcal{Z}_2(\mathbf{Rep}(H);\mathcal{D})}
        \end{tikzcd}.\] Recall that $\mathcal{Z}_2(\mathcal{D}) \cap \mathbf{Rep}(H) \simeq \mathbf{Rep}(H/N)$ and $\mathcal{Z}_2(\mathcal{A})^H \simeq \mathcal{Z}_2(\mathcal{D}) \lor \mathbf{Rep}(H)$. Hence, $\mathcal{D}$ is non-degenerate if and only if $\mathcal{A}$ is non-degenerate and $H=N$. In general, the 2-category of local modules of $\mathcal{C}$ in $\mathscr{Z}_1(\mathbf{2Rep}(G))$ is equivalent to the 2-category of local modules of $\mathcal{D}$ in $\mathscr{Z}_1(\Sigma \mathbf{Rep}(G))$, which is $\mathscr{Z}_1(\mathbf{Mod}(\mathcal{D}))$ by Corollary \ref{cor:ConnectedAndLagrangianEtaleAlgebrasInZ(SigmaRepG)}. Moreover, one can embed $\mathbf{2Rep}(H/N)$ into $\mathscr{Z}_1(\mathbf{Mod}(\mathcal{D}))$, whose de-equivariantization is exactly $\Gamma$.
    \end{proof}

    \begin{Corollary} \label{cor:ConnectedAlgebrasin2RepG}
        A connected {\'e}tale algebra in $\mathbf{2Rep}(G)$ is determined by:
        \begin{enumerate}
            \item A subgroup $H$ of $G$ up to conjugacy,
            
            \item A braided fusion category $\mathcal{A}$ with an $H$-action $\rho: H \to \mathbf{Aut}_{br}(\mathcal{A})$.
        \end{enumerate}

        The corresponding braided multifusion category with $G$-action is given by \[\mathcal{C} := \mathbf{Fun}_H \left( G,\mathcal{A} \right).\]
    \end{Corollary}

    \begin{Remark} \label{rmk:RightModulesOverCrossedBraidedFusionCategory}
        By Remark \ref{rmk:RightModulesOverBraidedFusionCategoryOverRepG} and Theorem \ref{thm:ClassificationOfConnectedCrossedBraidedMultiFusionCategories}, we can implicitly describe the 2-category of right modules over a connected {\'e}tale algebra $\mathcal{C}$ in $\mathscr{Z}_1(\mathbf{2Rep}(G))$. A discussion of the more general case can be found in Remark \ref{rmk:DomainWallsInCondensationOfTwistedCrossedBraidedFusionCategories}.
    \end{Remark}
    
\section{Twisted Crossed Braided Multifusion Categories} \label{sec:TwistedCrossedBraidedFusionCat}

We classify connected and Lagrangian {\'e}tale algebras in $\mathscr{Z}_1(\mathbf{2Rep}(G))$ in the previous section. By Example \ref{exmp:BraidedCrossedMultiFusionCat}, they are exactly connected and Lagrangian {\'e}tale algebras in $\mathscr{Z}_1(\mathbf{2Vect}_G)$. In this section, we would like to turn on a twisting $\pi \in \mathrm{H}^4(G;\Bbbk^\times)$ and generalize the classification of connected and Lagrangian {\'e}tale algebras from $\mathscr{Z}_1(\mathbf{2Rep}(G))$ to $\mathscr{Z}_1(\mathbf{2Vect}^\pi_G)$. We start by providing a concrete description of {\'e}tale algebras in $\mathscr{Z}_1(\mathbf{2Vect}^\pi_G)$.

\begin{Definition} \label{def:TwistedCrossedMonoidalCat}
    A $\pi$-\textbf{twisted} $G$-\textbf{crossed} finite semisimple category consists of:
    \begin{enumerate}
        \item A $G$-graded finite semisimple category $\mathcal{C} = \bigoplus_{g \in G} \mathcal{C}_g$;
        
        \item A projective $G$-action $\odot:G \to \mathbf{Aut}(\mathcal{C})$ such that:
        \begin{enumerate}
            \item $G$-grading is compatible with $G$-action: $\forall g, h \in G, g \odot \mathcal{C}_h \subseteq \mathcal{C}_{ghg^{-1}}$;
            
            \item It restricts to a $Z(g)$-action twisted by 3-cocycle $\tau_g(\pi)$ on $\mathcal{C}_g$ for each $g \in G$. (See Example \ref{exmp:DrinfeldCenterOfTwisted2VectG} for details.)
        \end{enumerate}
    \end{enumerate}
\end{Definition}

\begin{Remark}
    $\pi$-twisted $G$-crossed finite semisimple categories are exactly objects in the 2-category $\mathscr{Z}_1(\mathbf{2Vect}_G^\pi)$. Hence, one can also define $\pi$-twisted $G$-crossed functors and natural transforms as 1- and 2-morphisms within the same ambient 2-category $\mathscr{Z}_1(\mathbf{2Vect}_G^\pi)$.
\end{Remark}

\begin{Definition}
    A $\pi$-\textbf{twisted} $G$-\textbf{crossed monoidal} structure on a $G$-crossed finite semisimple category $\mathcal{C}$ consists of:
    \begin{enumerate}
        \item $\pi$-twisted $G$-crossed functors $\{\otimes_{g,h}: \mathcal{C}_g \boxtimes \mathcal{C}_h \to \mathcal{C}_{gh}\}_{g,h \in G}$ and $I:\mathbf{Vect} \to \mathcal{C}_e$;
        
        \item $\pi$-twisted $G$-crossed natural isomorphisms (given in components where $g,h,k \in G$) \begin{center}
            $\begin{tikzcd}[sep=small]
                {(\mathcal{C}_g \boxtimes \mathcal{C}_h) \boxtimes \mathcal{C}_k} \arrow[rrrr, "\omega_{g,h,k}"] \arrow[dd, "\otimes_{g,h} \boxtimes \mathcal{C}_k"'] &  & {} &  & {\mathcal{C}_g \boxtimes (\mathcal{C}_h \boxtimes \mathcal{C}_k)} \arrow[dd, "\mathcal{C}_g \boxtimes \otimes_{h,k}"] \\
                &  &                           &  &   \\ \mathcal{C}_{gh} \boxtimes \mathcal{C}_k \arrow[rrdd, "\otimes_{gh,k}"'] & {} \arrow[rr, Rightarrow, "\alpha_{g,h,k}", shorten < = 2ex, shorten > = 2ex] & {} & {} & \mathcal{C}_g \boxtimes \mathcal{C}_{hk} \arrow[lldd,"\otimes_{g,hk}"] \\
                                               &  &                           &  &   \\
                                               &  & \mathcal{C}_{ghk}    &  &  
            \end{tikzcd},$
            \end{center}
            
            \begin{center}
            \begin{tabular}{@{}c c@{}}
            $\begin{tikzcd}[sep=small]
            \mathcal{C}_g \arrow[rrrr, equal] \arrow[rrdd, "I \boxtimes \mathcal{C}_g"'] &  & {} \arrow[dd, Rightarrow, "\lambda_g"', near start, shorten > = 1ex] &  & \mathcal{C}_g \\
                                               &  &                           &  &   \\
                                               &  & \mathcal{C}_e \boxtimes \mathcal{C}_g \arrow[rruu, "\otimes_{e,g}"']     &  &  
            \end{tikzcd},$
            
            &
            
            $\begin{tikzcd}[sep=small]
                                              &  & \mathcal{C}_g \boxtimes \mathcal{C}_e \arrow[rrdd, "\otimes_{g,e}"] \arrow[dd, Rightarrow, "\rho_g", shorten > = 1ex, shorten < = 2ex] &  &   \\
                                              &  &                                             &  &   \\
            \mathcal{C}_g \arrow[rruu, "\mathcal{C}_g \boxtimes I"] \arrow[rrrr,equal] &  & {}                                          &  & \mathcal{C}_g
            \end{tikzcd},$

            \end{tabular}
            \end{center}
    \end{enumerate} satisfying the pentagon and triangle conditions:
    
    \begin{enumerate}
        \item [a.] For any $g,h,k,l \in G$, the following equation holds in $\mathbf{Hom}(((\mathcal{C}_g \boxtimes \mathcal{C}_h) \boxtimes \mathcal{C}_k) \boxtimes \mathcal{C}_l,\mathcal{C}_{ghkl})$:
    \end{enumerate}

    \begin{landscape}

    \settoheight{\diagramwidth}{\includegraphics[width=80mm]{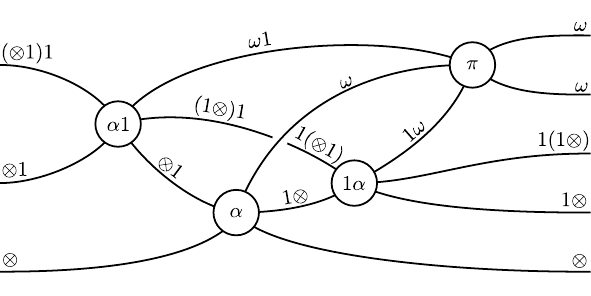}}

    \begin{equation}\label{eqn:crossedpentagon}
    \begin{tabular}{@{}cccc@{}}
    {\includegraphics[width=80mm]{Pictures/Gcrossedbraidedmonoidal/Pentagonleft.pdf}} &
    \raisebox{0.45\diagramwidth}{$=$} &
    \includegraphics[width=80mm]{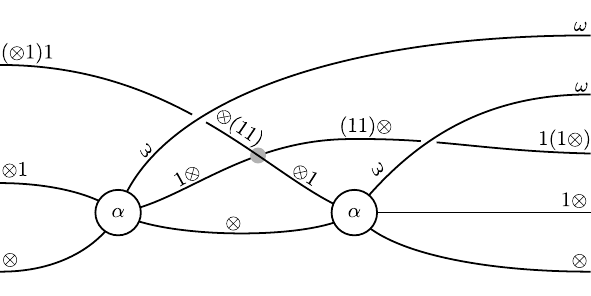} & \raisebox{0.45\diagramwidth}{.}
    \end{tabular}
    \end{equation}

    \vspace*{2cm}

    \begin{enumerate}
        \item [b.] For any $g,h \in G$, the following equation holds in $\mathbf{Hom}(\mathcal{C}_g \boxtimes \mathcal{C}_h,\mathcal{C}_{gh})$:
    \end{enumerate}

    \settoheight{\diagramwidth}{\includegraphics[width=48mm]{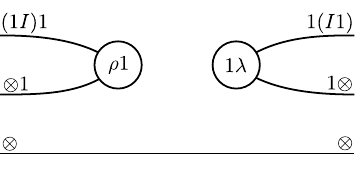}}

    \begin{equation}\label{eqn:crossedtriangle}
        \begin{tabular}{@{}cccc@{}}
        {\includegraphics[width=48mm]{Pictures/Gcrossedbraidedmonoidal/Triangleleft.pdf}} &
        \raisebox{0.45\diagramwidth}{$=$} &
        \includegraphics[width=48mm]{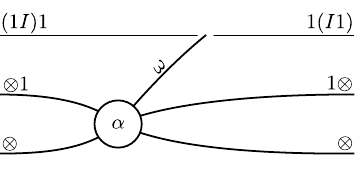} & \raisebox{0.45\diagramwidth}{.}
        \end{tabular}
        \end{equation}
    \end{landscape}
\end{Definition}

\begin{landscape}
\begin{Definition} \label{def:TwistedCrossedBraiding}
    A $\pi$-\textbf{twisted} $G$-\textbf{crossed braiding} is a $\pi$-twisted $G$-crossed natural isomorphism, 
    \[\begin{tikzcd}
        {\mathcal{C}_g \boxtimes \mathcal{C}_h}
            \arrow[rr,"\tau_{g,h}"]
            \arrow[dd,"\otimes_{g,h}"']
        & {}
        & {\mathcal{C}_h \boxtimes \mathcal{C}_g}
            \arrow[dd,"{}^g \odot \boxtimes \mathcal{C}_g"]
        \\ {}
        & {}
        & {}
            \arrow[ll,Rightarrow,shorten < = 8ex, shorten > = 8ex,"\beta_{g,h}"']
        \\ {\mathcal{C}_{gh}}
        & {}
        & {\mathcal{C}_{ghg^{-1}} \boxtimes \mathcal{C}_g}
            \arrow[ll,"\otimes_{ghg^{-1},g}"]
    \end{tikzcd}\] given in components as $\{(\beta_{g,h})_{X,Y}:X \otimes Y \to (g \odot Y) \otimes X\}_{g,h \in G}$ for any $g,h \in G$, $X$ in $\mathcal{C}_g$ and $Y$ in $\mathcal{C}_h$ and extended naturally to all objects in $\mathcal{C}$. Here, $\tau$ is the braiding on the ambient 2-category swapping two entries of Deligne product, and $^g \odot$ is the abbreviation for the action $g \odot -:\mathcal{C}_h \to \mathcal{C}_{ghg^{-1}}$.
    
    It satisfies the following two hexagon equations:
    
    \begin{enumerate}
        \item [a.] For any $g,h,k \in G$, the following equation holds in $\mathbf{Hom}((\mathcal{C}_g \boxtimes \mathcal{C}_h) \boxtimes \mathcal{C}_k,\mathcal{C}_{ghk})$:
    \end{enumerate}

    \settoheight{\diagramwidth}{\includegraphics[width=80mm]{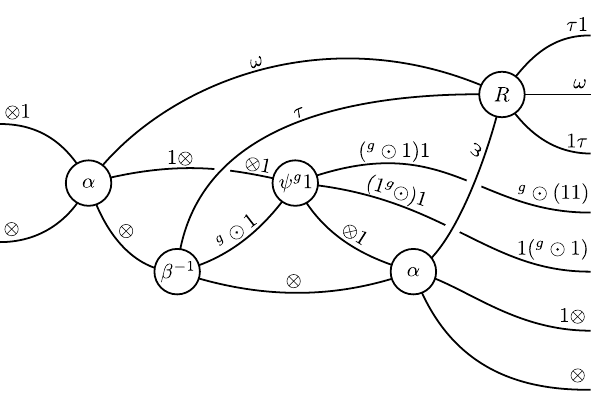}}

    \begin{equation}\label{eqn:crossedhexagon1}
    \begin{tabular}{@{}cccc@{}}
    {\includegraphics[width=80mm]{Pictures/Gcrossedbraidedmonoidal/Hexagon1left.pdf}} &
    \raisebox{0.45\diagramwidth}{$=$} &
    \includegraphics[width=64mm]{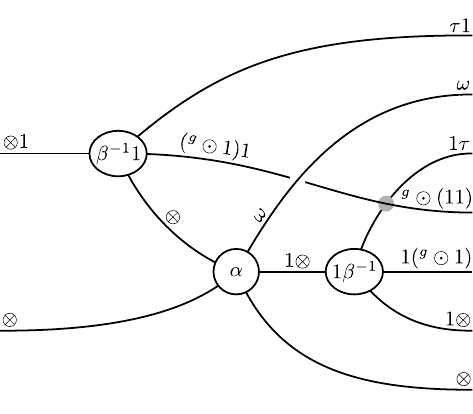} & \raisebox{0.45\diagramwidth}{,}
    \end{tabular}
    \end{equation} where $\psi^g$ is from the $G$-crossed functor structure on $\otimes$.

    \vspace*{2cm}

    \begin{enumerate}
        \item [b.] For any $g,h,k \in G$, the following equation holds in $\mathbf{Hom}(\mathcal{C}_g \boxtimes (\mathcal{C}_h \boxtimes \mathcal{C}_k),\mathcal{C}_{ghk})$:
    \end{enumerate}
    
    \settoheight{\diagramwidth}{\includegraphics[width=80mm]{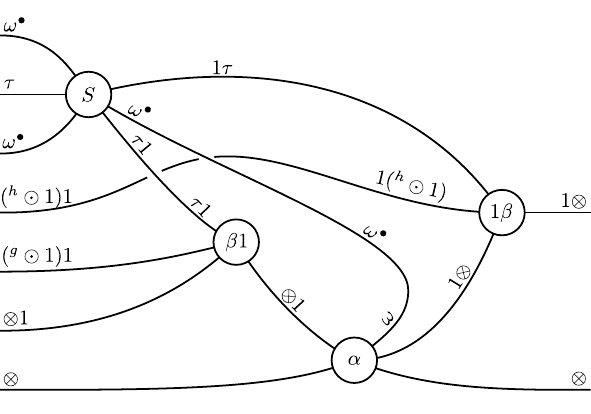}}

    \begin{equation}\label{eqn:crossedhexagon2}
    \begin{tabular}{@{}cccc@{}}
    {\includegraphics[width=80mm]{Pictures/Gcrossedbraidedmonoidal/Hexagon2left.pdf}} &
    \raisebox{0.45\diagramwidth}{$=$} &
    \includegraphics[width=80mm]{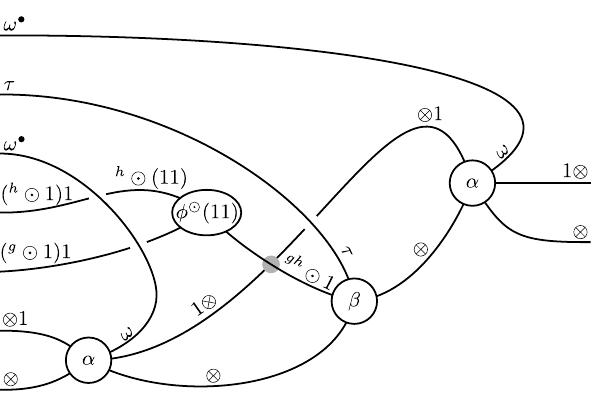} & \raisebox{0.45\diagramwidth}{,}
    \end{tabular}
    \end{equation} where $\phi^\odot$ witnesses from the projective $G$-action structure on $\mathcal{C}$.
\end{Definition}
\end{landscape}

\begin{Definition}
    Let $\mathcal{C}$ be a $\pi$-twisted $G$-crossed finite semisimple category with a $\pi$-twisted $G$-crossed monoidal structure. For $g \in G$ and object $X$ in $\mathcal{C}_g$, a \textbf{left dual} of $X$ consists of:
    \begin{enumerate}
        \item An object $Y$ in $\mathcal{C}_{g^{-1}}$;
        
        \item A morphism $c_{X}:I \to X \otimes Y$ called \textbf{coevaluation} and a morphism $e_{X}:Y \otimes X \to I$ called \textbf{evaluation};
    \end{enumerate} satisfying the two zig-zag equations:
    \[1_X = X \xrightarrow{\lambda_X} I \otimes X \xrightarrow{c_X 1} (X \otimes Y) \otimes X \xrightarrow{\omega_{X,Y,X}} X \otimes (Y \otimes X) \xrightarrow{1 e_X} X \otimes I \xrightarrow{\rho_X} X,\]
    \[1_Y = Y \xrightarrow{\rho^\bullet_Y} Y \otimes I \xrightarrow{1 c_X} Y \otimes (X \otimes Y) \xrightarrow{\omega_{Y,X,Y}^\bullet} (Y \otimes X) \otimes Y \xrightarrow{e_X 1} I \otimes Y \xrightarrow{\lambda^\bullet_Y} Y.\]

    Then, we extend the definition of left duals to any object in $\mathcal{C}$ by additivity. Dually, one can define the right dual of $X$. Notice in the above notation, $(X,c_X,e_X)$ gives a right dual of object $Y$.
\end{Definition}

\begin{Definition}
    A $\pi$-twisted $G$-crossed monoidal category is \textbf{rigid} if every object has a left dual and a right dual. 
\end{Definition}

\begin{Definition} \label{def:TwistedCrossedBraidedFusionCat}
    A $\pi$-twisted $G$-crossed multifusion category is a rigid $\pi$-twisted $G$-crossed finite semisimple category with a $\pi$-twisted $G$-crossed monoidal structure. 
    
    A $\pi$-twisted $G$-crossed braided multifusion category is a $\pi$-twisted $G$-crossed multifusion category with a compatible $\pi$-twisted $G$-crossed braiding.

    A $\pi$-twisted $G$-crossed braided fusion category is a $\pi$-twisted $G$-crossed braided multifusion category with a simple unit. In particular, its identity component $\mathcal{A}$ is a braided fusion category with $G$-action. We can also view this as a $\pi$-twisted $G$-crossed extension for $\mathcal{A}$.
\end{Definition}

\begin{Lemma}
    Separable algebras in $\mathscr{Z}_1(\mathbf{2Vect}^\pi_G)$ are exactly $\pi$-twisted $G$-crossed multifusion categories. {\'E}tale algebras in $\mathscr{Z}_1(\mathbf{2Vect}^\pi_G)$ are exactly $\pi$-twisted $G$-crossed braided multifusion categories.
\end{Lemma}

\begin{Remark} \label{rmk:TwistedCrossedBraidedFunctorAndNaturalTransform}
    The above lemma only requires us to unpack the definitions and carefully check the compatibility conditions. But some readers might be concerned with the precise meaning of the word \textit{exactly}. Indeed, if we would like to interpret it as a statement of equivalences of 2-categories, we need to first provide definitions of $\pi$-twisted $G$-crossed (braided) monoidal functors and natural transforms, and show that they form a 2-category. However, to avoid the verbosity, here we would rather take the alternative approach that, since on the level of objects there exist correspondences between the claimed notations, we choose to \textit{define} the $\pi$-twisted $G$-crossed (braided) monoidal functors and natural transforms as (braided) algebra 1- and 2-morphisms in $\mathscr{Z}_1(\mathbf{2Vect}^\pi_G)$ (where readers can refer to \cite{DY22,DX,Xu24} for the precise definitions). Curious readers are advised to work out the details by themselves.
\end{Remark}

\begin{Proposition} \label{prop:ClassificationOfTwistedCrossedBraidedExtension}
    Let $\mathcal{A}$ be a braided fusion category. Then a $\pi$-twisted $G$-crossed extension of $\mathcal{A}$ is equivalent to a 3-group morphism $\theta: N \rtimes^\pi \mathrm{B}^2 \Bbbk^\times \to \mathbf{Pic}(\mathcal{A})$, where $N$ is the support of this extension and $N \rtimes^\pi \mathrm{B}^2 \Bbbk^\times$ is the 3-group determined by the following fiber sequence: \[\begin{tikzcd}
        {\mathrm{B}^2 \Bbbk^\times}
            \arrow[r]
            \arrow[d]
        & {N \rtimes^\pi \mathrm{B}^2 \Bbbk^\times}
            \arrow[r]
            \arrow[d]
        & {*}
            \arrow[d]
        \\ {*}
            \arrow[r]
        & {N}
            \arrow[r,"\pi|_{N}"]
        & {\mathrm{B}^3 \Bbbk^\times}
    \end{tikzcd}. \]
\end{Proposition}

\begin{proof}
    This is a generalization of Proposition \ref{prop:CrossedBraidedExtensionAndPicard3Group}. Suppose $\mathcal{C} = \bigoplus_{g \in N} \mathcal{C}_g$ is a $\pi$-twisted $G$-crossed extension of $\mathcal{A}$ with support $N$, i.e. $\mathcal{A} \simeq \mathcal{C}_e$ and $N$ is a normal subgroup of $G$. Then the $G$-crossed monoidal structure $\otimes$ induces an $(\mathcal{A},\mathcal{A})$-bimodule structure\footnote{Note that by assumption the twist $\pi$ is \textit{normalized}, i.e. it vanishes if one of the arguments is graded by $e$.} on $\mathcal{C}_g$ for each $g \in N$. After passing to relative Deligne tensor product, we obtain equivalences of $(\mathcal{A},\mathcal{A})$-bimodule categories for $g,h \in G$: \[\mu_{g,h}: \mathcal{C}_g \boxtimes_\mathcal{A} \mathcal{C}_h \simeq \mathcal{C}_{gh}.\] Furthermore, the $G$-crossed braiding $\beta$ induces an identification between the left and right actions of $\mathcal{A}$ on $\mathcal{C}_g$ for each $g \in N$, promoting the above equivalences to right $\mathcal{A}$-module categories. In particular, we obtain a 1-group homomorphism $\pi_0(\theta):N \to \pi_0 \mathbf{Pic}(\mathcal{A})$.

    Next, we would like to extend $\pi_0(\theta)$ upward to a 2-group morphism. This is equivalent to choose the data $\alpha$ for $\mathcal{C}$ such that they still satisfy the $G$-crossed structure and the $\mathcal{A}$-actions, but we do not ask them to satisfy the pentagon condition. It turns out that there is no obstruction to do so, and we can always promote $\pi_0(\theta)$ to a 2-group morphism $\Pi_{\leq 1}(\theta):N \to \Pi_{\leq 1}\mathbf{Pic}(\mathcal{A})$.

    Lastly, to extend $\Pi_{\leq 1}(\theta)$ to a 3-group morphism, we need to impose the pentagon condition on $\alpha$. At this step we encounter the twisting $\pi$, in comparison to the untwisted case from \cite{ENO09}. By definition, $\pi$-twisted $G$-crossed fusion category $\mathcal{C}$ satisfies the pentagon condition (\ref{eqn:crossedpentagon}). This is equivalent to a homotopy between $\pi|_N: N \to \mathrm{B}^3 \Bbbk^\times$ and the obstruction class of $\Pi_{\leq 1}(\theta)$, so it induces the following fiber sequence: \[\begin{tikzcd}
        {N \rtimes^\pi \mathrm{B}^2 \Bbbk^\times}
            \arrow[r,dashed,"\theta"]
            \arrow[d]
        & {\mathbf{Pic}(\mathcal{A})}
            \arrow[d]
            \arrow[r]
        & {*}
            \arrow[d]
        \\ {N}
            \arrow[r,"\Pi_{\leq 1}(\theta)"']
        & {\Pi_{\leq 1}\mathbf{Pic}(\mathcal{A})}
            \arrow[r]
        & {\mathrm{B}^3 \Bbbk^\times}
    \end{tikzcd}, \] where the dashed arrow is the desired 3-group morphism.

    For the other direction, given a 3-group morphism $\theta:N \rtimes^\pi \mathrm{B}^2 \Bbbk^\times \to \mathbf{Pic}(\mathcal{A})$, we can construct a $\pi$-twisted $G$-crossed extension of $\mathcal{A}$ by taking $\mathcal{C}:= \bigoplus_{g \in N} \theta(g)$, with coherence data $(\alpha,\lambda,\rho,\beta)$ induced by the corresponding 3-group morphism structure on $\theta$.
\end{proof}

We are now ready to classify connected {\'e}tale algebras in $\mathscr{Z}_1(\mathbf{2Vect}^\pi_G)$.

\begin{Theorem} \label{thm:ClassificationOfConnectedTwistedCrossedBraidedMultiFusionCategories}
    A connected {\'e}tale algebra in $\mathscr{Z}_1(\mathbf{2Vect}^\pi_G)$ is determined by
    \begin{enumerate}
        \item A subgroup $H$ of $G$ up to conjugacy,
        
        \item A normal subgroup $N$ of $H$,
        
        \item A braided fusion category $\mathcal{A}$ with an $H$-action $\gamma: H \to \mathbf{Aut}_{br}(\mathcal{A})$.
        
        \item A 4-group morphism $H/N \times \mathrm{B}^3 \Bbbk^\times \to \mathbf{BrPic}(\mathbf{Mod}(\mathcal{A}))$, where
        \begin{itemize}
            \item The first component \[\Gamma: H/N \to \mathbf{BrPic}(\mathbf{Mod}(\mathcal{A})) \simeq \mathbf{Pic}(\mathscr{Z}_1(\mathbf{Mod}(\mathcal{A})))\] is an $H/N$-crossed extension of $\mathscr{Z}_1(\mathbf{Mod}(\mathcal{A}))$,
            
            \item The second component is the \textbf{pointwise inverse} of \[\iota: \mathrm{B}^3 \Bbbk^\times \to \mathbf{BrPic}(\mathbf{Mod}(\mathcal{A})),\] the 3-truncation in the Postnikov tower of $\mathbf{BrPic}(\mathbf{Mod}(\mathcal{A}))$, i.e. it induces the identity on the top homotopy group.
        \end{itemize}
        
        \item An invertible 2-morphism: \[
            \begin{tikzcd}
                {H}
                    \arrow[d,"{(q,\pi|_H)}"']
                    \arrow[r,"\gamma"]
                & {\mathbf{Aut}_{br}(\mathcal{A})}
                    \arrow[d,"\mathbf{Z}"]
                \\ {H/N \times \mathrm{B}^3 \Bbbk^\times}
                    \arrow[r,"{(\Gamma,{-}\iota)}"']
                    \arrow[ur,Rightarrow,shorten >=3ex,shorten <=4ex]
                & {\mathbf{BrPic}(\mathbf{Mod}(\mathcal{A}))}
            \end{tikzcd}, \]
        where the left arrow is a combination of the canonical quotient morphism $q:H \twoheadrightarrow H/N$ and the 4-cocycle $\pi|_H: H \to \mathrm{B}^3 \Bbbk^\times$, and the right arrow is explained in Remark \ref{rmk:PostnikovTowerOfPicard3Group}. The homotopy fiber of the above square produces a $\pi$-twisted $N$-crossed extension: \[\begin{tikzcd}
        {N \rtimes^\pi \mathrm{B}^2 \Bbbk^\times}
            \arrow[d]
            \arrow[r,"\theta"]
        & {\mathbf{Pic}(\mathcal{A})}
            \arrow[d,"\partial"]
        \\ {H}
            \arrow[r,"\gamma"']
            \arrow[ur,Rightarrow,shorten >=4ex,shorten <=4ex]
        & {\mathbf{Aut}_{br}(\mathcal{A})}
    \end{tikzcd}. \]
\end{enumerate}

    The corresponding $\pi$-twisted $G$-crossed braided multifusion category is given by \[\mathcal{C} := \mathbf{Fun}_H \left( G,\bigoplus_{g \in N} \theta(g) \right).\]
\end{Theorem}

\begin{proof}
    Suppose we are given a connected {\'e}tale algebra $\mathcal{C}$ in $\mathscr{Z}_1(\mathbf{2Vect}^\pi_G)$. We recognize that $\mathcal{C}$ is a $\pi$-twisted $G$-crossed braided multifusion category. After restricted to its identity component, it produces a braided multifusion category with $G$-action, which has to be connected as an {\'e}tale algebra in $\mathbf{2Rep}(G)$. By Corollary \ref{cor:ConnectedAlgebrasin2RepG}, we know that it is determined by a subgroup $H$ of $G$ up to conjugacy and a braided fusion category $\mathcal{A}$ with an $H$-action $\gamma:H \to \mathbf{Aut}_{br}(\mathcal{A})$. Moreover, we have an embedding $\mathbf{Fun}_H(G,\mathcal{A}) \hookrightarrow \mathcal{C}$ as the identity component. Furthermore, the support of its unit forms a braided embedding \[\mathcal{E}:=\mathbf{Fun}_H(G,\mathbf{Vect}) \hookrightarrow \mathbf{Fun}_H(G,\mathcal{A}) \hookrightarrow \mathcal{C}.\]

    By Proposition \ref{prop:ConnectedAndLagrangianEtaleAlgebrasInLocalModules}, we can first take the local modules of $\mathcal{E}$ in $\mathscr{Z}_1(\mathbf{2Vect}^\pi_G)$ and then $\mathcal{C}$ becomes a connected {\'e}tale algebra in this new braided fusion 2-category. Notice that \[\mathbf{Mod}^{loc}_{\mathscr{Z}_1(\mathbf{2Vect}^\pi_G)}(\mathcal{E}) \simeq \mathscr{Z}_1(\mathbf{2Vect}^{\pi|_H}_H)\] and then $\mathcal{C} \simeq \mathbf{Fun}_H(G,\mathcal{F})$, for some $\pi|_H$-twisted $H$-crossed braided \textit{fusion} category $\mathcal{F}$ extending the identity component $\mathcal{A}$. Then by Proposition \ref{prop:ClassificationOfTwistedCrossedBraidedExtension}, this data is equivalent to \begin{enumerate}
        \item A normal subgroup $N$ of $H$,
        
        \item A 3-group morphism $\theta:N \rtimes^\pi \mathrm{B}^2 \Bbbk^\times \to \mathbf{Pic}(\mathcal{A})$.
    \end{enumerate} together with its compatibility with the $H$-action on $\mathcal{A}$, which gives us an invertible 2-morphism filling the following diagram: \[\begin{tikzcd}
        {N \rtimes^\pi \mathrm{B}^2 \Bbbk^\times}
            \arrow[d]
            \arrow[r,"\theta"]
        & {\mathbf{Pic}(\mathcal{A})}
            \arrow[d,"\partial"]
        \\ {H}
            \arrow[r,"\gamma"']
        & {\mathbf{Aut}_{br}(\mathcal{A})}
    \end{tikzcd}. \] Finally, we observe that the above square is a homotopy fiber of the following square of 4-groups: \[\begin{tikzcd}
        {H}
            \arrow[d,"{(q,\pi|_H)}"']
            \arrow[r,"\gamma"]
        & {\mathbf{Aut}_{br}(\mathcal{A})}
            \arrow[d,"\mathbf{Z}"]
        \\ {H/N \times \mathrm{B}^3 \Bbbk^\times}
            \arrow[r,"{(\Gamma,{-}\iota)}"']
            \arrow[ur,Rightarrow,shorten >=3ex,shorten <=4ex]
        & {\mathbf{BrPic}(\mathbf{Mod}(\mathcal{A}))}
    \end{tikzcd}, \] where the equivariantization of $\Gamma$ gives the 2-category of local modules of $\mathcal{C}$ in $\mathscr{Z}_1(\mathbf{2Vect}_G^\pi)$. We save the details to the remaining part of this section (see in particular Theorem \ref{thm:LocalModulesOverTwistedCrossedBraidedFusionCategories}).
\end{proof}

\begin{Corollary} \label{cor:ClassificationOfLagrangianTwistedCrossedBraidedFusionCategories}
    A connected {\'e}tale algebra in $\mathscr{Z}_1(\mathbf{2Vect}^\pi_G)$ (determined by the data stated in the theorem) is Lagrangian if and only if $\mathcal{A}$ is non-degenerate and $N=H$. In other word, a Lagrangian algebra in $\mathscr{Z}_1(\mathbf{2Vect}^\pi_G)$ is equivalent to a sequence $(H,\mathcal{A},\gamma,\theta,\phi)$ where \begin{enumerate}
        \item $H$ is a subgroup of $G$ determined up to conjugation,
        
        \item $\mathcal{A}$ is a non-degenerate braided fusion category with $H$-action $\gamma$,
        
        \item 3-group morphism $\theta: H \rtimes^\pi \mathrm{B}^2 \Bbbk^\times \to \mathbf{Pic}(\mathcal{A})$ corresponds to a $\pi|_H$-twisted $H$-crossed extension of $\mathcal{A}$,
        
        \item $\phi$ is an equivalence of $H$-crossed braided fusion 2-categories: \[\mathbf{2Vect}^{\pi|_H}_H \boxtimes_{\mathbf{Fun}(H,\mathbf{2Vect})} \mathscr{Z}^H_1(\mathbf{Mod}(\mathcal{A})) \simeq \mathbf{2Vect}_H \boxtimes \mathscr{Z}_1(\mathbf{Mod}(\mathcal{A}))\] by Theorem \ref{thm:LocalModulesOverTwistedCrossedBraidedFusionCategories}, or equivalently the cell filling the following diagram:
    \end{enumerate} \[\begin{tikzcd}
        {H}
            \arrow[d,"{\pi|_H}"']
            \arrow[r,"\gamma"]
        & {\mathbf{Aut}_{br}(\mathcal{A})}
            \arrow[d,"\mathbf{Z}"]
        \\ {\mathrm{B}^3 \Bbbk^\times}
            \arrow[r,"{-}\iota"']
            \arrow[ur,Rightarrow,shorten >=3ex,shorten <=4ex]
        & {\mathbf{BrPic}(\mathbf{Mod}(\mathcal{A}))} 
    \end{tikzcd}. \]
\end{Corollary}

We would like to introduce the explicit descriptions for modules and local modules over $\pi$-twisted $G$-crossed braided multifusion categories. From now on, let $(\mathcal{C},\otimes^\mathcal{C},I,\alpha^\mathcal{C},\lambda^\mathcal{C},\rho^\mathcal{C})$ be a $\pi$-twisted $G$-crossed multifusion category.

\begin{Definition}
    A $\pi$-\textbf{twisted} $G$-\textbf{crossed finite semisimple right module category} over $\mathcal{C}$ consists of:
    \begin{enumerate}
        \item A $\pi$-twisted $G$-crossed finite semisimple category $\mathcal{M} = \bigoplus_{g \in G} \mathcal{M}_g$;
        
        \item A $\pi$-twisted $G$-crossed functor $\{\otimes^\mathcal{M}_{g,h}:\mathcal{M}_g \boxtimes \mathcal{C}_h \to \mathcal{M}_{gh} \}_{g,h \in G}$;
        
        \item $\pi$-twisted $G$-crossed natural isomorphisms in components for $g,h,k,l \in G$: \begin{center}
            $\begin{tikzcd}[sep=small]
                {(\mathcal{M}_g \boxtimes \mathcal{C}_h) \boxtimes \mathcal{C}_k} \arrow[rrrr, "\omega_{g,h,k}"] \arrow[dd, "\otimes^\mathcal{M}_{g,h} \boxtimes \mathcal{C}_k"'] &  & {} &  & {\mathcal{M}_g \boxtimes (\mathcal{C}_h \boxtimes \mathcal{C}_k)} \arrow[dd, "\mathcal{M}_g \boxtimes \otimes_{h,k}"] \\
                &  &                           &  &   \\ \mathcal{M}_{gh} \boxtimes \mathcal{C}_k \arrow[rrdd, "\otimes^\mathcal{M}_{gh,k}"'] & {} \arrow[rr, Rightarrow, "\alpha^\mathcal{M}_{g,h,k}", shorten < = 2ex, shorten > = 2ex] & {} & {} & \mathcal{M}_g \boxtimes \mathcal{C}_{hk} \arrow[lldd,"\otimes^\mathcal{M}_{g,hk}"] \\
                                               &  &                           &  &   \\
                                               &  & \mathcal{M}_{ghk}    &  &  
            \end{tikzcd},$
            \end{center}
            
            \begin{center}
            $\begin{tikzcd}[sep=small]
                                              &  & \mathcal{M}_g \boxtimes \mathcal{C}_e \arrow[rrdd, "\otimes^\mathcal{M}_{g,e}"] \arrow[dd, Rightarrow, "\rho^\mathcal{M}_g", shorten > = 1ex, shorten < = 2ex] &  &   \\
                                              &  &                                             &  &   \\
            \mathcal{M}_g \arrow[rruu, "\mathcal{M}_g \boxtimes I"] \arrow[rrrr,equal] &  & {}                                          &  & \mathcal{M}_g
            \end{tikzcd},$
            \end{center}
    \end{enumerate} satisfying the pentagon and triangle conditions:

    \begin{enumerate}
        \item [a.] For any $g,h,k,l \in G$, the following equation holds in $\mathbf{Hom}(((\mathcal{M}_g \boxtimes \mathcal{C}_h) \boxtimes \mathcal{C}_k) \boxtimes \mathcal{C}_l,\mathcal{M}_{ghkl})$:
    \end{enumerate}

    \begin{landscape}

    \settoheight{\diagramwidth}{\includegraphics[width=80mm]{Pictures/Gcrossedbraidedmonoidal/Pentagonleft.pdf}}

    \begin{equation}\label{eqn:modulecrossedpentagon}
    \begin{tabular}{@{}cccc@{}}
    {\includegraphics[width=80mm]{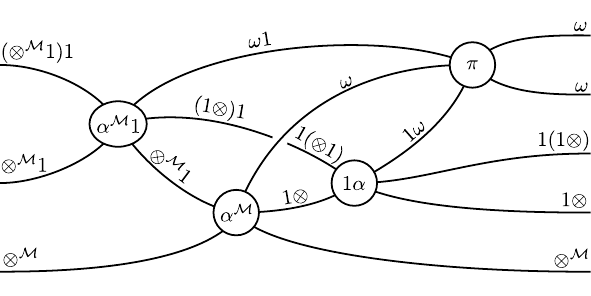}} &
    \raisebox{0.45\diagramwidth}{$=$} &
    \includegraphics[width=80mm]{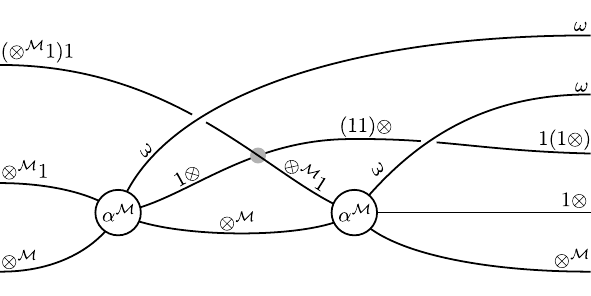} & \raisebox{0.45\diagramwidth}{.}
    \end{tabular}
    \end{equation}

    \vspace*{2cm}

    \begin{enumerate}
        \item [b.] For any $g,h \in G$, the following equation holds in $\mathbf{Hom}(\mathcal{M}_g \boxtimes \mathcal{C}_h,\mathcal{M}_{gh})$:
    \end{enumerate}

    \settoheight{\diagramwidth}{\includegraphics[width=48mm]{Pictures/Gcrossedbraidedmonoidal/Triangleleft.pdf}}

    \begin{equation}\label{eqn:modulecrossedtriangle}
        \begin{tabular}{@{}cccc@{}}
        {\includegraphics[width=48mm]{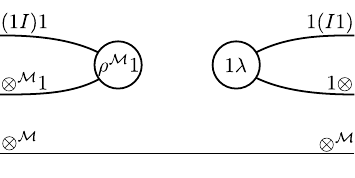}} &
        \raisebox{0.45\diagramwidth}{$=$} &
        \includegraphics[width=48mm]{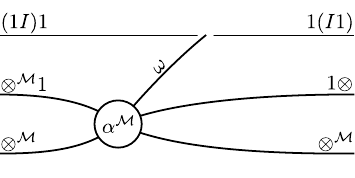} & \raisebox{0.45\diagramwidth}{.}
        \end{tabular}
        \end{equation}
    \end{landscape}
\end{Definition}

\begin{landscape}
    \begin{Definition} \label{def:crossedholonomy}
        A $\pi$-\textbf{twisted} $G$-\textbf{crossed holonomy} on a $\pi$-twisted $G$-crossed finite semisimple right module category $\mathcal{M}$ over $\mathcal{C}$ is a $\pi$-twisted $G$-crossed natural isomorphism given in components for $g,h \in G$: \[ \begin{tikzcd}
            {\mathcal{M}_g \boxtimes \mathcal{C}_h}
                \arrow[r,"\tau_{g,h}"]
                \arrow[dd,"\otimes_{g,h}"']
            & {\mathcal{C}_h \boxtimes \mathcal{M}_g}
                \arrow[r,"{}^g \odot \, \boxtimes \, \mathcal{M}_g"]
            & {\mathcal{C}_{ghg^{-1}} \boxtimes \mathcal{M}_g}
                \arrow[d,"\tau_{ghg^{-1},g}"]
            \\ {}
            & {}
            & {\mathcal{M}_g \boxtimes \mathcal{C}_{ghg^{-1}}}
                \arrow[d,"{}^{ghg^{-1}} \odot \, \boxtimes \, \mathcal{C}_{ghg^{-1}}"]
                \arrow[ll,Rightarrow,shorten < = 12ex, shorten > = 12ex,"h^\mathcal{M}_{g,h}"']
            \\ {\mathcal{M}_{gh}}
            & {}
            & {\mathcal{M}_{g h gh^{-1}g^{-1}} \boxtimes \mathcal{C}_{ghg^{-1}}}
                \arrow[ll,"\otimes_{g h gh^{-1}g^{-1},ghg^{-1}}"]
        \end{tikzcd},  \] or equivalently, a family of isomorphisms $(\hbar^\mathcal{M}_{g,h})_{M,X}: M \otimes^\mathcal{M} (g^{-1} \odot X) \to (h \odot M) \otimes^\mathcal{M} X $ natural for $g,h \in G$, $M$ in $\mathcal{M}_g$ and $X$ in $\mathcal{C}_h$, such that it is $\pi$-twisted $G$-crossed and satisfies the following coherence equations:

        \begin{enumerate}
            \item [a.] For any $g,h,k \in G$, object $X$ in $\mathcal{M}_g$, we have the triangle commutes 
            \begin{equation} \label{eqn:crossedholonomytriangle} 
                \begin{tikzcd}
                    {M \otimes^\mathcal{M} (g^{-1} \odot I)}
                        \arrow[rr,"{(\hbar^\mathcal{M}_{g,e})_{M,I}}"]
                        \arrow[dr]
                    & {}
                    & {(e \odot M) \otimes^\mathcal{M} I}
                        \arrow[dl]
                    \\ {}
                    & {M \otimes^\mathcal{M} I}
                    & {}
                \end{tikzcd}, 
            \end{equation} where the left leg witnesses the $G$-crossed functor structure on the unit $I$ and the right leg witnesses the unitality of projective $G$-action structure on $\mathcal{M}$.
        \end{enumerate}
        
        \begin{enumerate}
            \item [b.] For any $g,h,k \in G$, object $X$ in $\mathcal{M}_g$, $X$ in $\mathcal{C}_h$ and $Y$ in $\mathcal{C}_k$, the following diagram commutes:
        \end{enumerate}
        
        \begin{equation} \label{eqn:crossedholonomyassociatorI}
            \begin{tikzcd}[row sep=30pt,column sep=70pt]
                {(M \otimes^\mathcal{M} (g^{-1} \odot X)) \otimes^\mathcal{M} (g^{-1} h^{-1} \odot Y)}
                    \arrow[d,"\hbar^\mathcal{M}_{hg,k}"']
                    \arrow[r,"\alpha^\mathcal{M}_{g,g^{-1} h g,g^{-1} h^{-1} k h g}"]
                & {M \otimes^\mathcal{M} ((g^{-1} \odot X) \otimes (g^{-1} h^{-1} \odot Y))}
                    \arrow[d]
                \\ {(k \odot (M \otimes^\mathcal{M} (g^{-1} \odot X))) \otimes^\mathcal{M} Y}
                    \arrow[d]
                & {M \otimes^\mathcal{M} (g^{-1} \odot (X \otimes (h^{-1} \odot Y)))}
                    \arrow[d,"1 \beta_{h,h^{-1} k h}"]
                \\ {((k \odot M) \otimes^\mathcal{M} (k g^{-1} \odot X )) \otimes^\mathcal{M} Y}
                    \arrow[d,"\alpha^\mathcal{M}_{kgk^{-1}, k g^{-1} h g k^{-1},k}"']
                & {M \otimes^\mathcal{M} (g^{-1} \odot ((h \odot (h^{-1} \otimes Y)) \otimes X))}
                    \arrow[d]
                \\ {(k \odot M) \otimes^\mathcal{M} ((k g^{-1} \odot X ) \otimes Y)}
                    \arrow[d]
                & {M \otimes^\mathcal{M} ((g^{-1} \odot Y) \otimes (g^{-1} \odot X))}
                    \arrow[d,"(\alpha^\mathcal{M}_{g,g^{-1} k g,g^{-1} h g})^\bullet"]
                \\ {(k \odot M) \otimes^\mathcal{M} ((k \odot (g^{-1} \odot X) ) \otimes Y)}
                    \arrow[d,"1 (\beta^\mathcal{M}_{k,g^{-1} h g})^\bullet"']
                & {(M \otimes^\mathcal{M} (g^{-1} \odot Y)) \otimes^\mathcal{M} (g^{-1} \odot X)}
                    \arrow[d,"\hbar^\mathcal{M}_{g,k} 1"]
                \\ {(k \odot M) \otimes^\mathcal{M} (Y \otimes (g^{-1} \odot X))}
                    \arrow[r,"(\alpha^\mathcal{M}_{k g k^{-1},k,g^{-1} h g})^\bullet"']
                & {((k \odot M) \otimes^\mathcal{M} Y) \otimes^\mathcal{M} (g^{-1} \odot X)}
            \end{tikzcd} 
        \end{equation} where the unlabeled arrows are compositions of isomorphisms from various $G$-crossed functor structures and projective $G$-action structure on $\mathcal{M}$.
        
        \vspace*{1cm}
        \begin{enumerate}
            \item [c.] For any $g,h,k \in G$, object $X$ in $\mathcal{M}_g$, $X$ in $\mathcal{C}_h$ and $Y$ in $\mathcal{C}_k$, the following diagram commutes:
        \end{enumerate}

        \begin{equation} \label{eqn:crossedholonomyassociatorII} 
            \begin{tikzcd}[row sep=30pt,column sep=70pt]
                {(M \otimes^\mathcal{M} (g^{-1} \odot X)) \otimes^\mathcal{M} (g^{-1} h^{-1} \odot Y)}
                    \arrow[d,"\hbar^\mathcal{M}_{g,h} 1"']
                    \arrow[r,"\alpha^\mathcal{M}_{g,g^{-1} h g,g^{-1} h^{-1} k h g}"]
                & {M \otimes^\mathcal{M} ((g^{-1} \odot X) \otimes (g^{-1} h^{-1} \odot Y))}
                    \arrow[d]
                \\ {((h \odot M) \otimes^\mathcal{M} X) \otimes^\mathcal{M} (g^{-1} h^{-1} \odot Y)}
                    \arrow[d,"\hbar^\mathcal{M}_{hg,k} "']
                & {M \otimes^\mathcal{M} (g^{-1} \odot (X \otimes (h^{-1} \odot Y)))}
                    \arrow[d,"1 \beta_{h,h^{-1} k h}"]
                \\ {(k \odot ((h \odot M) \otimes^\mathcal{M} X)) \otimes^\mathcal{M} Y}
                    \arrow[d]
                & {M \otimes^\mathcal{M} (g^{-1} \odot ((h \odot (h^{-1} \odot Y)) \otimes X))}
                    \arrow[d]
                \\ {((kh \odot M) \otimes^\mathcal{M} (k \odot X)) \otimes^\mathcal{M} Y}
                    \arrow[d,"\alpha^\mathcal{M}_{k h g h^{-1} k^{-1}, k h k^{-1}, k}"']
                & {M \otimes^\mathcal{M} (g^{-1} \odot (Y \otimes X))}
                    \arrow[d,"1 \beta_{k,h}"]
                \\ {(kh \odot M) \otimes^\mathcal{M} ((k \odot X) \otimes Y)}
                & {M \otimes^\mathcal{M} (g^{-1} \odot ((k \odot X) \otimes Y))}
                    \arrow[l,"\hbar^\mathcal{M}_{g,kh}"]
            \end{tikzcd} 
        \end{equation} where the unlabeled arrows are also compositions of isomorphisms from various $G$-crossed functor structures and projective $G$-action structure on $\mathcal{M}$.

        \vspace*{0.5cm}
        A $\pi$-twisted $G$-crossed finite semisimple right module category over $\mathcal{C}$ together with a $\pi$-twisted $G$-crossed holonomy is called a $\pi$-\textbf{twisted} $G$-\textbf{crossed finite semisimple local module category} over $\mathcal{C}$.
    \end{Definition} 
\end{landscape}

\begin{Lemma}
    Right modules over $\mathcal{C}$ in $\mathscr{Z}_1(\mathbf{2Vect}^\pi_G)$ are exactly $\pi$-twisted $G$-crossed finite semisimple right module categories over $\mathcal{C}$. 
    
    Local modules over $\mathcal{C}$ in $\mathscr{Z}_1(\mathbf{2Vect}^\pi_G)$ are exactly $\pi$-twisted $G$-crossed finite semisimple local module categories over $\mathcal{C}$.
\end{Lemma}

\begin{Remark}
    Just as Remark \ref{rmk:TwistedCrossedBraidedFunctorAndNaturalTransform}, the above lemma can be interpreted as a statement of equivalences of 2-categories, where the definitions of $\pi$-twisted $G$-crossed (local) module functors and natural transforms are given as (local) module 1- and 2-morphisms in $\mathscr{Z}_1(\mathbf{2Vect}^\pi_G)$. Curious readers are encouraged to work out the details by themselves.
\end{Remark}

Let $\mathcal{A}$ be a braided fusion category with a $G$-action $\gamma: G \to \mathbf{Aut}_{br}(\mathcal{A})$. Clearly, $\mathcal{A}$ is a connected étale algebra in $\Sigma \mathbf{Rep}(G)$. Since the identity component of $\mathscr{Z}_1(\mathbf{2Vect}^\pi_G)$ can be identified with $\Sigma \mathbf{Rep}(G)$, $\mathcal{A}$ could also be viewed as a connected étale algebra in $\mathscr{Z}_1(\mathbf{2Vect}^\pi_G)$.

\begin{Theorem} \label{thm:ModulesOverBraidedFusionCategoryWithGAction}
    $\pi$-twisted $G$-crossed finite semisimple right module categories over $\mathcal{A}$ form a 2-category equivalent to $\mathbf{2Vect}^\pi_{G} \boxtimes_G \mathbf{Mod}(\mathcal{A})$, where $G$ acts on $\mathbf{2Vect}^\pi_G$ via \textbf{conjugation} and on $\mathbf{Mod}(\mathcal{A})$ by Lemma \ref{lem:BraidedAutoEquivalencesAreTheSameAsMonoidalAuto2EquivalencesOfModularCategories}.
\end{Theorem}

\begin{proof}
    Let us consider the following 2-functor \[\mathbf{B}: \mathbf{2Vect}^\pi_G \boxtimes \mathbf{Mod}(\mathcal{A}) \to \mathbf{Mod}_{\mathscr{Z}_1(\mathbf{2Vect}^\pi_G)}(\mathcal{A}),\]
    \begin{itemize}
        \item Given a $\pi$-twisted $G$-graded finite semisimple category $\mathcal{M} = \bigoplus_{g \in G} \mathcal{M}_g$ and a right $\mathcal{A}$-module category $\mathcal{N}$, this 2-functor assigns \[\mathbf{B}(\mathcal{M},\mathcal{N}) := \bigoplus_{x,g \in G} \mathcal{M}_{xgx^{-1}} \boxtimes \mathcal{N}_{\langle x \rangle},\] where $\mathcal{N}_{\langle x \rangle}$ has the same underlying category as $\mathcal{N}$ and its right $\mathcal{A}$-action is twisted by the braided auto-equivalence $x \odot -: \mathcal{A} \to \mathcal{A}$.
         
        Each component is equipped with the right $\mathcal{A}$-action. On the other hand, the $\pi$-twisted $G$-crossed structure on $\mathbf{B}(\mathcal{M},\mathcal{N})$ is defined as follows:
        \begin{itemize}
            \item For $g \in G$, the $G$-graded component is $\bigoplus_{x \in G} \mathcal{M}_{xgx^{-1}} \boxtimes \mathcal{N}_{\langle x \rangle}$.
            
            \item For $h \in G$, the $G$-action sends $G$-graded components $\mathbf{B}(\mathcal{M},\mathcal{N})_g$ to $\mathbf{B}(\mathcal{M},\mathcal{N})_{hgh^{-1}}$ via relabeling the index $x$ by $xh$, and sending each $\mathcal{N}_{\langle xh \rangle}$ to $\mathcal{N}_{\langle x \rangle}$ using the braided auto-equivalence $h \odot -: \mathcal{A} \to \mathcal{A}$.
            
            \item The twisted $G$-action is inherited from the twisted $G$-grading on $\mathcal{M}$.
        \end{itemize}
        
        \item There is a balanced 2-functor structure $\alpha^\mathbf{B}$, i.e. a $\pi$-twisted $G$-crossed isomorphism $\alpha^\mathbf{B}_{\mathcal{M},y,\mathcal{N}}: \mathbf{B}(\mathcal{M} \cdot y,\mathcal{N}) \simeq \mathbf{B}(\mathcal{M},y \cdot \mathcal{N})$, natural with respect to any $\mathcal{M} = \bigoplus_{g \in G} \mathcal{M}_g$ in $\mathbf{2Vect}^\pi_G$, $y \in G$ and $\mathcal{N}$ in $\mathbf{Mod}(\mathcal{A})$, defined via relabeling the index $x \mapsto x' = x y^{-1}$: \[\bigoplus_{x,g \in G} \mathcal{M}_{x y^{-1}g y x^{-1}} \boxtimes \mathcal{N}_{\langle x \rangle} \equiv \bigoplus_{x',g \in G} \mathcal{M}_{x' g (x')^{-1}} \boxtimes \mathcal{N}_{\langle x' y \rangle}. \] There are canonical invertible modifications $\omega^\mathbf{B}$ and $\gamma^\mathbf{B}$:
        \[\begin{tikzcd}
            {}
            & {\mathbf{B}\left(\mathcal{M} \cdot y,z \cdot \mathcal{N}\right)}
                \arrow[dd,"\alpha^\mathbf{B}_{\mathcal{M},y,z\mathcal{N}}"]
            \\ {\mathbf{B}\left(\mathcal{M} \cdot yz,\mathcal{N}\right)}
                \arrow[ru,"\alpha^\mathbf{B}_{\mathcal{M}y,z,\mathcal{N}}"]
                \arrow[rd,"\alpha^\mathbf{B}_{\mathcal{M},yz,\mathcal{N}}"']
                \arrow[r,Rightarrow,shorten <=10pt, shorten >=10pt,"\omega^\mathbf{B}_{\mathcal{M},y,z,\mathcal{N}}"]
            & {}
            \\ {}
            & {\mathbf{B}\left(\mathcal{M},yz \cdot \mathcal{N} \right)}
        \end{tikzcd}, \quad \begin{tikzcd}
            {}
            & {\mathbf{B}\left(\mathcal{M},\mathcal{N}\right)}
                \arrow[dd,equal]
            \\ {\mathbf{B}\left(\mathcal{M},\mathcal{N}\right)}
                \arrow[ru,equal]
                \arrow[rd,"\alpha^\mathbf{B}_{\mathcal{M},e,\mathcal{N}}"']
                \arrow[r,Rightarrow,shorten <=10pt, shorten >=10pt,"\gamma^\mathbf{B}_{\mathcal{M},\mathcal{N}}"]
            & {}
            \\ {}
            & {\mathbf{B}\left(\mathcal{M},\mathcal{N}\right)}
        \end{tikzcd},\] satisfying the coherence conditions in \cite[Definition 2.1.1]{D10}.
    \end{itemize}

    Thus, by 3-universal property of relative 2-Deligne tensor product \cite[Definition 2.2.1]{D10}, this induces an essentially unique 2-functor \[\mathbf{2Vect}^\pi_G \boxtimes_G \mathbf{Mod}(\mathcal{A}) \to \mathbf{Mod}_{\mathscr{Z}_1(\mathbf{2Vect}^\pi_G)}(\mathcal{A}).\] Furthermore, it actually is an equivalence, as the above balanced 2-functor $(\mathbf{B},\alpha^\mathbf{B},\omega^\alpha,\gamma^\alpha)$ turns out to be universal among all such balanced 2-functors. Details are omitted here for brevity.
\end{proof}

\begin{Corollary}
    If $\mathcal{A}$ has trivial $G$-action, then one has \[\mathbf{Mod}_{\mathscr{Z}_1(\mathbf{2Vect}^\pi_G)}(\mathcal{A}) \simeq  \mathscr{Z}_1(\mathbf{2Vect}^\pi_G) \boxtimes \mathbf{Mod}(\mathcal{A}). \]
\end{Corollary}

\begin{Remark}
    One can extend $\mathbf{B}$ to a 3-condensation, where \[\mathbf{Mod}_{\mathscr{Z}_1(\mathbf{2Vect}^\pi_G)}(\mathcal{A}) \to \mathbf{2Vect}^\pi_G \boxtimes \mathbf{Mod}(\mathcal{A})\] forgets a $\pi$-twisted $G$-crossed finite semisimple right module category $\mathcal{M}$ over $\mathcal{A}$ to its underlying $\pi$-twisted $G$-crossed finite semisimple category $\bigoplus_{g \in G} \mathcal{M}_g$ with the $\mathcal{A}$-action on each component.
\end{Remark}

\begin{Theorem} \label{thm:LocalModulesOverBraidedFusionCategoryWithGAction}
    $\pi$-twisted $G$-crossed finite semisimple local module categories over $\mathcal{A}$ form a 2-category equivalent to the equivariantization of \[\mathbf{2Vect}^\pi_{G} \boxtimes_{\mathbf{Fun}(G,\mathbf{2Vect})} \mathscr{Z}^G_1(\mathbf{Mod}(\mathcal{A})),\] where $\mathscr{Z}^G_1(\mathbf{Mod}(\mathcal{A}))$ is the $G$-crossed Drinfeld center introduced in Proposition \ref{prop:CrossedDrinfeld2CenterOfModA}, and $\mathbf{Fun}(G,\mathbf{2Vect})$ actions on $\mathbf{2Vect}^\pi_{G}$ and $\mathscr{Z}^G_1(\mathbf{Mod}(\mathcal{A}))$ are both induced by viewing them as $G$-graded 2-categories.
\end{Theorem}

\begin{proof}
    Let's first consider $\mathbf{2Vect}^\pi_{G}$ and $\mathscr{Z}^G_1(\mathbf{Mod}(\mathcal{A}))$ as $G$-graded 2-categories. Equivalently, applying the extension theory \cite{D11}, one has:
    \begin{enumerate}
        \item $\mathbf{2Vect}^\pi_{G}$ is determined by the 4-cocycle $\pi$ viewed as a 4-group morphism $G \to \mathrm{B}^3 \Bbbk^\times$. There is a further embedding of $\mathrm{B}^3 \Bbbk^\times$ into $\mathbf{BrPic}(\mathbf{Mod}(\mathcal{A}))$ as the top piece of Postnikov tower. This further induces a 4-group morphism $G \to \mathbf{BrPic}(\mathbf{Mod}(\mathcal{A}))$.
        
        \item $\mathscr{Z}^G_1(\mathbf{Mod}(\mathcal{A}))$ is a $G$-crossed extension of $\mathscr{Z}_1(\mathbf{Mod}(\mathcal{A}))$, and it corresponds to a 4-group morphism $G \xrightarrow{\gamma} \mathbf{Aut}_{br}(\mathcal{A}) \xrightarrow{\mathbf{Z}} \mathbf{BrPic}(\mathbf{Mod}(\mathcal{A}))$, as mentioned in Remark \ref{rmk:CrossedDrinfeld2CenterOfModA4Morphism}.
    \end{enumerate}

    On the other hand, the 2-category $\mathbf{2Vect}^\pi_{G} \boxtimes_{\mathbf{Fun}(G,\mathbf{2Vect})} \mathscr{Z}^G_1(\mathbf{Mod}(\mathcal{A}))$ is simply the components in the absolute 2-Deligne's tensor product $\mathbf{2Vect}^\pi_{G} \boxtimes \mathscr{Z}^G_1(\mathbf{Mod}(\mathcal{A}))$ which corresponds to the \textit{pointwise} tensor product with respect to the $G$-gradings. Therefore, this 2-category is characterized by the 4-group morphism $G \to \mathbf{BrPic}(\mathbf{Mod}(\mathcal{A}))$ which is the \textit{pointwise} product of the two 4-group morphisms above. In particular, this induces a $G$-crossed braided fusion 2-category structure on $\mathbf{2Vect}^\pi_{G} \boxtimes_{\mathbf{Fun}(G,\mathbf{2Vect})} \mathscr{Z}^G_1(\mathbf{Mod}(\mathcal{A}))$.

    Finally, the identification of $\mathbf{Mod}^{loc}_{\mathscr{Z}_1(\mathbf{2Vect}_G^\pi)}(\mathcal{A})$ with the equivariantization of $\mathbf{2Vect}^\pi_{G} \boxtimes_{\mathbf{Fun}(G,\mathbf{2Vect})} \mathscr{Z}^G_1(\mathbf{Mod}(\mathcal{A}))$ follows almost by definition: in Proposition \ref{prop:CrossedDrinfeld2CenterOfModA}, we have identified the crossed Drinfeld center $\mathscr{Z}^G_1(\mathbf{Mod}(\mathcal{A}))$ as the 2-category of \textit{untwisted} $G$-\textit{graded} finite semisimple braided $\mathcal{A}$-module categories. Therefore, after tensoring with $\mathbf{2Vect}^\pi_{G}$ pointwise with respect to the $G$-gradings, we obtain the 2-category of $\pi$-twisted $G$-\textit{graded} finite semisimple local $\mathcal{A}$-module categories. Lastly, $\pi$-twisted $G$-\textit{crossed} finite semisimple local $\mathcal{A}$-module categories are given by equivariantization of the above $G$-crossed braided fusion 2-category.
    
    In summary, the de-equivariantization of $\mathbf{Mod}^{loc}_{\mathscr{Z}_1(\mathbf{2Vect}_G^\pi)}(\mathcal{A})$ is a $G$-crossed extension of $\mathscr{Z}_1(\mathbf{Mod}(\mathcal{A}))$ produced by pointwise tensor product of $\mathbf{2Vect}^\pi_{G}$ and $\mathscr{Z}^G_1(\mathbf{Mod}(\mathcal{A}))$.
\end{proof}

\begin{Remark}
    Equivalently, one has a commutative diagram: \[\begin{tikzcd}
        {G}
            \arrow[d,"{(\mathrm{id},\pi)}"']
            \arrow[r,"\gamma"]
        & {\mathbf{Aut}_{br}(\mathcal{C}_e)}
            \arrow[d,"\mathbf{Z}"]
        \\ {G \times \mathrm{B}^3 \Bbbk^\times}
            \arrow[r,"{(\Gamma,{-}\iota)}"']
        & {\mathbf{BrPic}(\mathbf{Mod}(\mathcal{C}_e))}
    \end{tikzcd},\] where $\iota:\mathrm{B}^3 \Bbbk^\times \to \mathbf{BrPic}(\mathbf{Mod}(\mathcal{C}_e))$ is the top layer of the Postnikov tower, and $\Gamma$ is the 4-group morphism corresponding to the de-equivariantization of $\mathbf{Mod}^{loc}_{\mathscr{Z}_1(\mathbf{2Vect}_G^\pi)}(\mathcal{C}_e)$ as a $G$-crossed extension of $\mathscr{Z}_1(\mathbf{Mod}(\mathcal{A}))$.
\end{Remark}

\begin{Corollary}
    If $\mathcal{A}$ has the trivial $G$-action, then one has \[\mathbf{Mod}^{loc}_{\mathscr{Z}_1(\mathbf{2Vect}^\pi_G)}(\mathcal{A}) \simeq  \mathscr{Z}_1(\mathbf{2Vect}^\pi_G) \boxtimes \mathscr{Z}_1(\mathbf{Mod}(\mathcal{A})). \]
\end{Corollary}

\begin{Theorem} \label{thm:LocalModulesOverTwistedCrossedBraidedFusionCategories}
    Let $\mathcal{C}$ be a $\pi$-twisted $G$-crossed braided fusion category, or equivalently, a $\pi$-twisted $G$-crossed extension of a braided fusion category $\mathcal{C}_e$ with a $G$-action $\gamma: G \to \mathbf{Aut}_{br}(\mathcal{C}_e)$. Then the 2-category of $\pi$-twisted $G$-crossed finite semisimple local module categories over $\mathcal{C}$ corresponds to a 4-group morphism \[ \Gamma: G/N \to \mathbf{BrPic}(\mathbf{Mod}(\mathcal{C}_e)),\] where $N$ is the support of $\mathcal{C}$. Moreover, we have an equivalence of $G$-crossed braided fusion 2-categories: \[\bigoplus_{g \in G} \Gamma(gN) \simeq \mathbf{2Vect}^\pi_{G} \boxtimes_{\mathbf{Fun}(G,\mathbf{2Vect})} \mathscr{Z}^G_1(\mathbf{Mod}(\mathcal{C}_e)).\]
\end{Theorem}

\begin{proof}
    The support of $\mathcal{C}$ forms a normal subgroup $N$ of $G$, and $\mathcal{C}$ gives rise to a 3-group morphism $\theta:N \rtimes^\pi \mathrm{B}^2 \Bbbk^\times \to \mathbf{Pic}(\mathcal{C}_e)$ by Proposition \ref{prop:ClassificationOfTwistedCrossedBraidedExtension}, together with a homotopy filling the following square: \[\begin{tikzcd}
        {N \rtimes^\pi \mathrm{B}^2 \Bbbk^\times}
            \arrow[r,"\theta"]
            \arrow[d]
        & {\mathbf{Pic}(\mathcal{C}_e)}
            \arrow[d,"\partial"]
        \\ {G}
            \arrow[r,"\gamma"']
        & {\mathbf{Aut}(\mathcal{C}_e)}
    \end{tikzcd}.\]

    Meanwhile, notice that one has the following exact sequence of 4-groups, i.e. both a fiber sequence and a cofiber sequence: \[\begin{tikzcd}
        {N \rtimes^\pi \mathrm{B}^2 \Bbbk^\times}
            \arrow[r]
            \arrow[d]
        & {G}   
            \arrow[d,"{(q,\pi)}"]
        \\ {*}
            \arrow[r]
        & {G/N \times \mathrm{B}^3 \Bbbk^\times}
    \end{tikzcd},\] where $q:G \twoheadrightarrow G/N$ is the canonical quotient homomorphism. Therefore, it induces a 4-group morphism filling the following fiber sequence of 4-group morphisms: \[\begin{tikzcd}
        {N \rtimes^\pi \mathrm{B}^2 \Bbbk^\times}
            \arrow[r,"\theta"]
            \arrow[d]
        & {\mathbf{Pic}(\mathcal{C}_e)}
            \arrow[d,"\partial"]
        \\ {G}
            \arrow[r,"\gamma"]
            \arrow[d,"{(q,\pi)}"']
        & {\mathbf{Aut}(\mathcal{C}_e)}
            \arrow[d,"\mathbf{Z}"]
        \\ {G/N \times \mathrm{B}^3 \Bbbk^\times}
            \arrow[r,dashed,"{(\Gamma,{-}\iota)}"']
        & {\mathbf{BrPic}(\mathbf{Mod}(\mathcal{C}_e))}
    \end{tikzcd}.\] Notice that the second component $\mathrm{B}^3 \Bbbk^\times \to \mathbf{BrPic}(\mathbf{Mod}(\mathcal{C}_e))$ has to coincide with $-\iota$ since its fiber $\theta$ also induces the identity map between the top layer of Postnikov systems by Proposition \ref{prop:ClassificationOfTwistedCrossedBraidedExtension}. Hence we now focus on the 4-group morphism $\Gamma: G/N \to \mathbf{BrPic}(\mathbf{Mod}(\mathcal{C}_e))$, which is equivalent to an $H/N$-crossed extension of $\mathscr{Z}_1(\mathbf{Mod}(\mathcal{C}_e))$ by the extension theory of fusion 2-categories developed in \cite{D11}.
    
    By Theorem \ref{thm:LocalModulesOverBraidedFusionCategoryWithGAction}, the de-equivariantization of $\mathbf{Mod}^{loc}_{\mathscr{Z}_1(\mathbf{2Vect}_G^\pi)}(\mathcal{C}_e)$ (with respect to $G$) is equivalent to the pointwise tensor product of $\mathbf{2Vect}^\pi_{G}$ and $\mathscr{Z}^G_1(\mathbf{Mod}(\mathcal{C}_e))$. So it factors through the square: \[\begin{tikzcd}
        {G}
            \arrow[r,"\gamma"]
            \arrow[d,"{(\mathrm{id},\pi)}"']
        & {\mathbf{Aut}(\mathcal{C}_e)}
            \arrow[dd,"\mathbf{Z}"]
        \\ {G \times \mathrm{B}^3 \Bbbk^\times}
            \arrow[d,"{(q,\mathrm{id})}"']
            \arrow[dr,dashed]
        & {}
        \\ {G/N \times \mathrm{B}^3 \Bbbk^\times}
            \arrow[r,"{(\Gamma,{-}\iota)}"']
        & {\mathbf{BrPic}(\mathbf{Mod}(\mathcal{C}_e))}
    \end{tikzcd},\] which can be interpreted as the equivalence stated: \[\bigoplus_{g \in G} \Gamma(gN) \simeq \mathbf{2Vect}^\pi_{G} \boxtimes_{\mathbf{Fun}(G,\mathbf{2Vect})} \mathscr{Z}^G_1(\mathbf{Mod}(\mathcal{C}_e)).\]
    
    The final thing to check is that $\Gamma$ is the de-equivariantization of the 2-category of local modules over $\mathcal{C}$ with respect to $G/N$. To condense a $\pi$-twisted $G$-crossed braided fusion category $\mathcal{C}$, by Proposition \ref{prop:ConnectedAndLagrangianEtaleAlgebrasInLocalModules}, we decompose this into two steps: first, we condense its identity component $\mathcal{C}_e$, then we view $\mathcal{C}$ as a connected {\'e}tale algebra in local modules over $\mathcal{C}_e$ and further condense the entire $\mathcal{C}$, i.e.

    \[\mathbf{Mod}^{loc}_{\mathscr{Z}_1(\mathbf{2Vect}^\pi_G)}(\mathcal{C}) \simeq \mathbf{Mod}^{loc}_{\mathbf{Mod}^{loc}_{\mathscr{Z}_1(\mathbf{2Vect}^\pi_G)}(\mathcal{C}_e)}(\mathcal{C}).\]

    Take any local $\mathcal{C}$-module $\mathcal{W}$ in $\mathbf{Mod}^{loc}_{\mathscr{Z}_1(\mathbf{2Vect}^\pi_G)}(\mathcal{C}_e)$, it has an underlying $\pi$-twisted $G$-crossed finite semisimple local $\mathcal{C}_e$-module category. If $\mathcal{M}$ is a simple local $\mathcal{C}_e$-module category lying in $\mathcal{W}$, then for any $g \in N$, one also has $\mathcal{M} \boxtimes_{\mathcal{C}_e} \mathcal{C}_g$ lying in $\mathcal{W}$ as another simple $\mathcal{C}_e$-module category, with $G$-braiding translated. Omitting details here, we observe that $\pi$-twisted $G$-gradings of local modules over $\mathcal{C}$ necessarily factor through the quotient $q: G \twoheadrightarrow G/N$, such that for any $g \in G$ and $n \in N$, there is a canonical identification of $G$-graded components $\mathcal{W}_{gn} \simeq \mathcal{W}_g \boxtimes_{\mathcal{C}_e} \mathcal{C}_n$ induced by the right $\mathcal{C}$-action on $\mathcal{W}$. This functorial process implies that $\Gamma$ construct above has to agree with the de-equivariantization of local modules over $\mathcal{C}$ with respect to $G/N$.
\end{proof}

\begin{Remark}
    If the $G$-action on $\mathcal{C}_e$ is trivial, then the 4-cocycle $\pi$ descends onto $G/N$, or equivalently there exists a group 3-cochain $\phi$ on $N$ such that $\mathrm{d} \phi = \pi|_N$. \[\begin{tikzcd}
        {}
        & {N}
            \arrow[rd,"\pi|_N"]
            \arrow[d,hook]
            \arrow[ld]
        & {}
        \\ {*}
            \arrow[rd]
        & {G}
            \arrow[r,"\pi"]
            \arrow[d,twoheadrightarrow]
        & {\mathrm{B}^3 \Bbbk^\times}
        \\ {}
        & {G/N}
            \arrow[ru,dashed,"\pi|_{G/N}"']
        & {}
    \end{tikzcd} \quad \Leftrightarrow \quad \exists\phi\text{ fills the outer square.}\] In this case, the local modules over $\mathcal{C}$ are given by \[\mathbf{Mod}^{loc}_{\mathscr{Z}_1(\mathbf{2Vect}^\pi_G)}(\mathcal{C}) \simeq  \mathscr{Z}_1(\mathbf{2Vect}^{\pi|_{G/N}}_{G/N}) \boxtimes \mathscr{Z}_1(\mathbf{Mod}(\mathcal{C}_e)). \]
\end{Remark}

\begin{Remark} \label{rmk:DomainWallsInCondensationOfTwistedCrossedBraidedFusionCategories}
    In the proof of Theorem \ref{thm:LocalModulesOverTwistedCrossedBraidedFusionCategories}, we describe the braided fusion 2-categories involved in the two-step condensation process. On the other hand, it is also important to understand the domain walls involved in the condensation process, which are mathematically described by two fusion 2-categories with monoidal bimodule 2-categorical structures. Furthermore, these two domain walls could be fused together via relative 2-Deligne tensor product to produce a new fusion 2-category, corresponding to the entire condensation process done in one step. We leave the details for future work. 
\end{Remark}

\section{Application to Bosonic Fusion 2-Categories} \label{sec:ClassificationOfBosonicFusion2Categories}

Our previous results, combined with {Décoppet}'s classification result of fusion 2-categories up to Morita equivalence \cite{D9}, can be applied to the classification of bosonic fusion 2-categories. See also the recent work \cite{DHJFNPPRY}.

\begin{Definition}
    A \textbf{bosonic fusion 2-category} is a fusion 2-category $\mathfrak{C}$ with $\mathcal{Z}_2(\Omega \mathfrak{C})$ Tannakian, i.e. there is a finite group $H$ such that $\mathcal{Z}_2(\Omega \mathfrak{C}) \simeq \mathbf{Rep}(H)$ as symmetric fusion categories, where braided fusion category $\Omega \mathfrak{C} := \mathbf{End}_\mathfrak{C}(I)$ is the looping of $\mathfrak{C}$.
\end{Definition}

\begin{Theorem}[{\cite[Theorem 4.1.6]{D9}}] \label{thm:CenterOfBosonicFusion2Categories}
    A bosonic fusion 2-category $\mathfrak{C}$ is always Morita equivalent to $\mathbf{2Vect}_G^\pi \boxtimes \mathbf{Mod}(\mathcal{A})$, where $G$ is a finite group $G$, $\pi \in \mathrm{H}^4(G;\Bbbk^\times)$ and $\mathcal{A}$ is a non-degenerate braided fusion category. In particular, one obtains an equivalence of braided fusion 2-categories $\mathscr{Z}_1(\mathfrak{C}) \simeq \mathscr{Z}_1(\mathbf{2Vect}^\pi_G)$.
\end{Theorem}

Recall that $\mathfrak{C}$ is enriched over its Drinfeld center $\mathscr{Z}_1(\mathfrak{C})$ \cite{D4}. Moreover, this enrichment respects the monoidal structures, hence for any algebra $A$ in $\mathfrak{C}$, the enriched hom $[I,A]$ provides an algebra in $\mathscr{Z}_1(\mathfrak{C})$.

\begin{Proposition} \label{prop:FullCenterOfUnit}
    The enriched endo-hom of the monoidal unit $[I,I]$ is endowed with a Lagrangian algebra structure in $\mathscr{Z}_1(\mathbf{2Vect}^\pi_G)$ such that \[\mathfrak{C} \simeq \mathbf{Mod}_{\mathscr{Z}_1(\mathfrak{C})}([I,I]).\]
\end{Proposition}

\begin{proof}
    Décoppet has prove \cite[Theorem 4.2.2]{D4} that any finite semisimple module 2-category $\mathfrak{M}$ over a fusion 2-category $\mathfrak{D}$ is canonically enriched in $\mathfrak{D}$ in the sense of Garner and Schulman \cite[Section 3.1]{GS}. In particular, since the forgetful functor $\mathscr{Z}_1(\mathfrak{C}) \to \mathfrak{C}$ is dominant \cite[Corollary 5.3.6]{D9}, one has \[\mathfrak{C} \simeq \mathbf{Mod}_{\mathscr{Z}_1(\mathfrak{C})}([I,I])\] as finite semisimple module 2-categories over $\mathscr{Z}_1(\mathfrak{C})$. The enriched endo-hom $[I,I]$ is a rigid algebra by construction. By \cite[Theorem 3.2.4]{D7}, it is also separable.

    We would like to promote the above equivalence to an equivalence of fusion 2-categories. This requires that the enrichment of $\mathfrak{C}$ in $\mathscr{Z}_1(\mathfrak{C})$ needs to have certain compatibility with the monoidal structures. This is a reminiscent of the 1-categorical case done by Morrison and Penneys \cite{MP17}. In Appendix \ref{sec:M2CatEnrichedInB2Cat}, we provide a detailed proof of this fact. Therefore, $[I,I]$ is indeed an étale algebra in $\mathscr{Z}_1(\mathbf{2Vect}^\pi_G)$, and this promotes desired equivalence to an equivalence of fusion 2-categories.

    By \cite[Theorem 4.4.5]{Xu24}, local modules over $[I,I]$ is characterized as the braided centralizer: \[\mathbf{Mod}^{loc}_{\mathscr{Z}_1(\mathfrak{C})}([I,I]) \simeq \mathscr{Z}_2(\mathscr{Z}_1(\mathfrak{C}) \to \mathscr{Z}_1(\mathbf{Mod}_{\mathscr{Z}_1(\mathfrak{C})}([I,I]))) \] \[\simeq \mathscr{Z}_2(\mathscr{Z}_1(\mathfrak{C}) \xrightarrow{\mathrm{Id}} \mathscr{Z}_1(\mathfrak{C})) \simeq \mathbf{2Vect}. \] In other word, étale algebra $[I,I]$ is Lagrangian in $\mathscr{Z}_1(\mathfrak{C})$.
\end{proof}

\begin{Theorem} \label{thm:ClassificationOfBosonicFusion2Categories}
    Combining the above theorem and proposition, one finds that every bosonic fusion 2-category corresponds to a Lagrangian algebra in the Drinfeld center $\mathscr{Z}_1(\mathbf{2Vect}^\pi_G)$ for some finite group $G$ and $\pi \in \mathrm{H}^4(G;\Bbbk^\times)$. 
    
    By Corollary \ref{cor:ClassificationOfLagrangianTwistedCrossedBraidedFusionCategories}, this Lagrangian algebra is determined by a subgroup $H$ of $G$, a non-degenerate braided fusion category $\mathcal{A}$ with $H$-action $\gamma$ and a faithful $\pi$-twisted $H$-crossed extension, which is the fiber of the following square of 4-group morphisms \[\begin{tikzcd}
        {H \rtimes^\pi \mathrm{B}^2 \Bbbk^\times}
            \arrow[r,"\theta"]
            \arrow[d]
        & {\mathbf{Pic}(\mathcal{A})}
            \arrow[d,"\partial"]
        \\ {H}
            \arrow[d,"{\pi|_H}"']
            \arrow[r,"\gamma"]
            \arrow[ur,Rightarrow,shorten >=4ex,shorten <=5ex]
        & {\mathbf{Aut}_{br}(\mathcal{A})}
            \arrow[d,"\mathbf{Z}"]
        \\ {\mathrm{B}^3 \Bbbk^\times}
            \arrow[r,"{-}\iota"']
            \arrow[ur,Rightarrow,shorten >=3ex,shorten <=4ex]
        & {\mathbf{BrPic}(\mathbf{Mod}(\mathcal{A}))} 
    \end{tikzcd}. \]
\end{Theorem}

\begin{Remark}
    One can directly extract the above data from a bosonic fusion 2-category $\mathfrak{C}$ as follows. Given any braided fusion 2-category $\mathfrak{C}$, one can extract a finite group $G$ and a 4-cocycle $\pi \in \mathrm{H}^4(G;\Bbbk^\times)$ via Theorem \ref{thm:CenterOfBosonicFusion2Categories}. In particular, one has an equivalence of symmetric fusion categories: $\Omega \mathscr{Z}_1(\mathfrak{C}) \simeq \mathbf{Rep}(G)$. Meanwhile, since $\mathfrak{C}$ is bosonic, one has another symmetric fusion category $\mathcal{Z}_2(\Omega \mathfrak{C}) \simeq \mathbf{Rep}(H)$, from which one can reconstruct a finite group $H$. 
    
    Notice that the canonical forgetful 2-functor $\mathscr{Z}_1(\mathfrak{C}) \to \mathfrak{C}$ induces a braided functor of the endo-categories: $\Omega \mathscr{Z}_1(\mathfrak{C}) \to \Omega \mathfrak{C}$, which has to factor through the Müger center: $\Omega \mathscr{Z}_1(\mathfrak{C}) \twoheadrightarrow \mathcal{Z}_2(\Omega \mathfrak{C}) \hookrightarrow \Omega \mathfrak{C}$. This identifies $H$ as a subgroup of $G$ up to conjugation.

    Now we get a braided embedding $\mathbf{Rep}(H) \simeq \mathcal{Z}_2(\Omega \mathfrak{C}) \hookrightarrow \Omega \mathfrak{C}$, so the de-equivariantization of $\Omega \mathfrak{C}$ is a braided fusion category $\mathcal{A}$ with an $H$-action $\gamma$. Finally, the $\pi$-twisted $H$-crossed extension is obtained by taking the fiber of the square of 4-group morphisms: \[\begin{tikzcd}
        {H}
            \arrow[d,"{\pi|_H}"']
            \arrow[r,"\gamma"]
        & {\mathbf{Aut}_{br}(\mathcal{A})}
            \arrow[d,"\mathbf{Z}"]
        \\ {\mathrm{B}^3 \Bbbk^\times}
            \arrow[r,"{{-}\iota}"']
            \arrow[ur,Rightarrow,shorten >=3ex,shorten <=4ex]
        & {\mathbf{BrPic}(\mathbf{Mod}(\mathcal{A}))}
    \end{tikzcd}.\]
\end{Remark}

\begin{Remark}
    The above parametrization is not optimal. Suppose the 4-cocycle $\pi$ is trivial, for any 3-cocycle $\omega$ on $G$, viewed as an $G$-crossed extension of trivial algebra $\mathbf{Vect}$, one can verify that the 2-category of right modules over this Lagrangian algebra is always given by \textit{the same fusion 2-category}: \[\mathbf{Mod}_{\mathscr{Z}_1(\mathbf{2Vect}_G)}(\mathbf{Vect}^\omega_G) \simeq \mathbf{2Rep}(G). \]
    
    One could have two \textit{non-equivalent} Lagrangian algebras in $\mathscr{Z}_1(\mathbf{2Vect}^\pi_G)$ whose 2-categories of right modules are \textit{equivalent} as fusion 2-categories. This subtlety is recorded in the collection of all possible twistings of the canonical forgetful 2-functor $\mathscr{Z}_1(\mathbf{2Vect}^\pi_G) \simeq \mathscr{Z}_1(\mathfrak{C}) \to \mathfrak{C}$. The author suspects that this phenomenon is related to the existence of $G$-symmetry protected topological orders in (1+1)D \cite{KLWZZ}. It would be interesting to investigate this in the future work from a physical perspective.
\end{Remark}

\appendix
\appendixpage
\addappheadtotoc

\section{Enriched Monoidal 2-Categories} \label{sec:M2CatEnrichedInB2Cat}

Given a semistrict\footnote{In particular we use \textit{Gray monoids} as the underlying models of semistrict multifusion 2-categories, see \cite{GPS}.} monoidal 2-category $(\mathfrak{C},\Box,I)$, we denote the \textit{interchanger} and its inverse using the string diagram notation as follows\footnote{We use the string diagram notation for 2-categories as in \cite{GS,DX}: regions are labelled by objects (which are usually omitted), lines are labelled by 1-morphisms whose composition goes from top to bottom, and coupons are labelled by 2-morphisms whose composition goes from left to right. We often omit the product $\Box$ and write $1$ for identity morphisms.}: 

\newlength{\diagramlength}
\settoheight{\diagramlength}{\includegraphics[width=20mm]{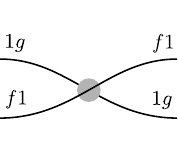}}

\[\begin{tabular}{c c c c}
\includegraphics[width=20mm]{Pictures/Appendex/inter.pdf}\quad\raisebox{0.45\diagramlength}{,}\ \ \ \   & \ \ \ \  \includegraphics[width=20mm]{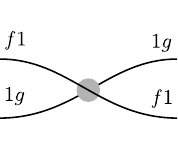}\quad\raisebox{0.45\diagramlength}{,}
\end{tabular}\] where $f:x_0 \to y_0$ and $g:x_1 \to y_1$ are 1-morphisms in $\mathfrak{C}$.

Let $F,G:\mathfrak{A}\rightarrow \mathfrak{B}$ be two 2-functors, and let $\tau:F \to G$ be a 2-natural transform. Its \textit{naturality} provides a 2-isomorphism between $\tau_y \circ F(f)$ and $G(f) \circ \tau_x$ for any 1-morphism $f:x \to y$ in $\mathfrak{A}$. We denote this 2-isomorphism and its inverse using the following string diagrams: 

\settoheight{\diagramlength}{\includegraphics[width=20mm]{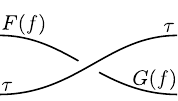}}
\[ \begin{tabular}{c c c c}
\includegraphics[width=20mm]{Pictures/Appendex/naturality.pdf}\quad\raisebox{0.45\diagramlength}{,}\ \ \ \   & \ \ \ \  \includegraphics[width=20mm]{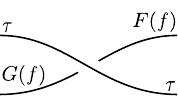}\quad\raisebox{0.45\diagramlength}{.}
\end{tabular}\]

A \textit{braiding} on a semistrict monoidal 2-category $(\mathfrak{C},\Box,I)$, by Gurski \cite{Gur2}, consists of:

\begin{enumerate}
    \item An adjoint 2-natural isomorphism $b_{x,y}: x \, \Box \, y \to y \, \Box \, x$ for any objects $x,y$ in $\mathfrak{C}$.
    
    \item Two invertible modifications $R$ and $S$, given on the objects $x,y,z$ in $\mathfrak{C}$ as
    \begin{center}
    \begin{tabular}{@{}c c@{}}
    
    $\begin{tikzcd}
    xyz \arrow[rr, "b"] \arrow[rd, "b1"'] & {} \arrow[d, Rightarrow, "R"]          & yzx \\
                                          & yxz \arrow[ru, "1b"'] &    
    \end{tikzcd},$
    
    &
    
    $\begin{tikzcd}
    xyz \arrow[rr, "b_2"] \arrow[rd, "1b"'] & {} \arrow[d, Rightarrow, "S"]     & zxy \\
                                            & xzy \arrow[ru, "b1"'] &    
    \end{tikzcd},$
    \end{tabular}
    \end{center}
    
    \noindent where the subscript in $b_2$ records were the braiding occurs. We will write a $b$ instead of a $b_1$ as this can be confused with $b1$.
\end{enumerate} subject to the following relations, which are taken from of \cite[Section 2.1.1]{DY22}:

\begin{enumerate}
    \item [a.] For every objects $x,y,z,w$ in $\mathfrak{C}$, we have
\end{enumerate}

\settoheight{\diagramlength}{\includegraphics[width=37.5mm]{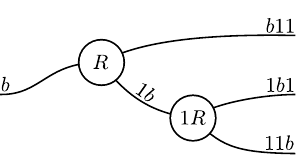}}

\begin{equation}\label{eqn:braidingaxiom1}
\begin{tabular}{@{}ccc@{}}

\includegraphics[width=37.5mm]{Pictures/Appendex/braidingaxiom1.pdf} & \raisebox{0.45\diagramlength}{$=$} &
\includegraphics[width=37.5mm]{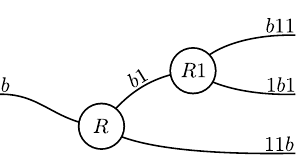}

\end{tabular}
\end{equation}

\begin{enumerate}
\item [] in $\mathbf{Hom}_{\mathfrak{C}}(x \, \Box \, y \, \Box \, z \, \Box \, w, y \, \Box \, z \, \Box \, w \, \Box \, x)$,
\item [b.] For every objects $x,y,z,w$ in $\mathfrak{C}$, we have
\end{enumerate}

\settoheight{\diagramlength}{\includegraphics[width=37.5mm]{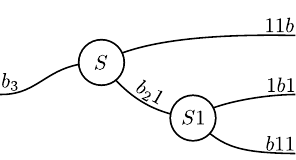}}

\begin{equation}\label{eqn:braidingaxiom2}
\begin{tabular}{@{}ccc@{}}

\includegraphics[width=37.5mm]{Pictures/Appendex/braidingaxiom3.pdf} & \raisebox{0.45\diagramlength}{$=$} &
\includegraphics[width=37.5mm]{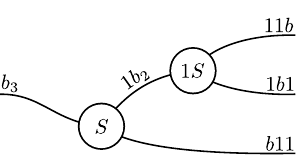}

\end{tabular}
\end{equation}

\begin{enumerate}
\item [] in $\mathbf{Hom}_{\mathfrak{C}}(x \, \Box \, y \, \Box \, z \, \Box \, w, w \, \Box \, x \, \Box \, y \, \Box \, z)$,
\item [c.] For every objects $x,y,z,w$ in $\mathfrak{C}$, we have
\end{enumerate}

\settoheight{\diagramlength}{\includegraphics[width=37.5mm]{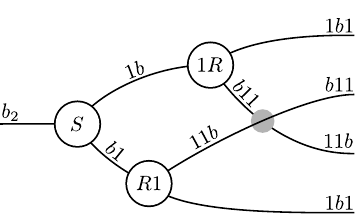}}

\begin{equation}\label{eqn:braidingaxiom3}
\begin{tabular}{@{}ccc@{}}

\includegraphics[width=45mm]{Pictures/Appendex/braidingaxiom5.pdf} & \raisebox{0.45\diagramlength}{$=$} &
\includegraphics[width=37.5mm]{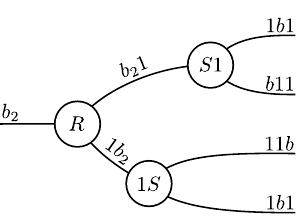}

\end{tabular}
\end{equation}

\begin{enumerate}
\item [] in $\mathbf{Hom}_{\mathfrak{C}}(x \, \Box \, y \, \Box \, z \, \Box \, w,z \, \Box \, w \, \Box \, x \, \Box \, y)$,
\item [d.] For every objects $x,y,z$ in $\mathfrak{C}$, we have
\end{enumerate}

\settoheight{\diagramlength}{\includegraphics[width=45mm]{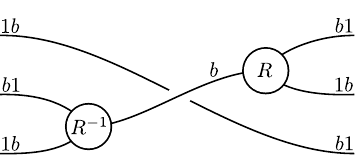}}

\begin{equation}\label{eqn:braidingaxiom4}
\begin{tabular}{@{}ccc@{}}

\includegraphics[width=45mm]{Pictures/Appendex/braidingaxiom7.pdf} & \raisebox{0.45\diagramlength}{$=$} &
\includegraphics[width=45mm]{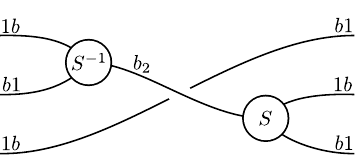}

\end{tabular}
\end{equation}

\begin{enumerate}
\item [] in $\mathbf{Hom}_{\mathfrak{C}}(x \, \Box \, y \, \Box \, z, z \, \Box \, y \, \Box \, x)$,
\item [e.] For every object $x$ in $\mathfrak{C}$, the adjoint 2-natural isomorphisms $$b_{x,I}:x \equiv x \, \Box \, I\rightarrow I \, \Box \, x \equiv x \textrm{ and } b_{I,x}:x \equiv I \, \Box \, x \rightarrow x \, \Box \, I \equiv x$$ are the identity adjoint 2-natural isomorphisms,

\item [f.] For every objects $x,y,z$ in $\mathfrak{C}$, the modifications $R_{x,y,z}$ and $S_{x,y,z}$ are the identity 2-isomorphisms if any of $x$, $y$, or $z$ is equal to $I$.
\end{enumerate}

\begin{Definition} \label{def:Enriched2Category}
    An enriched 2-category over a semistrict monoidal 2-category $\mathfrak{C}$ consists of:
    \begin{enumerate}
        \item A set\footnote{We assume smallness in general and ignore the size issue.} of objects $Ob(\mathfrak{M})$;

        \item For objects $x,y$ in $Ob(\mathfrak{M})$, an enriched hom object $[x,y]^\mathfrak{M}_\mathfrak{C}$ in $\mathfrak{C}$;

        \item For object $x$ in $Ob(\mathfrak{M})$, an enriched unit $j_x:I \to [x,y]^\mathfrak{M}_\mathfrak{C}$;

        \item For objects $x,y,z$ in $Ob(\mathfrak{M})$, an enriched composition \[m_{x,y,z}:[y,z]^\mathfrak{M}_\mathfrak{C} \,\Box \, [x,y]^\mathfrak{M}_\mathfrak{C} \to [x,z]^\mathfrak{M}_\mathfrak{C};\]

        \item For objects $x,y$ in $Ob(\mathfrak{M})$, enriched unitors in $\mathfrak{C}$ \[\begin{tikzcd}[sep=30pt]
            {I \, \Box \, [x,y]^\mathfrak{M}_\mathfrak{C}}
                \arrow[r,equal]
                \arrow[d,"j_y 1"']
            & {I \, \Box \, [x,y]^\mathfrak{M}_\mathfrak{C}}
                \arrow[d,equal]
            \\ {[y,y]^\mathfrak{M}_\mathfrak{C} \, \Box \, [x,y]^\mathfrak{M}_\mathfrak{C}}
                \arrow[r,"m_{x,y,y}"']
                \arrow[ur,Rightarrow,shorten <=20pt, shorten >=20pt,"\lambda_{x,y}"]
            & {[x,y]^\mathfrak{M}_\mathfrak{C}}
        \end{tikzcd},\] \[\begin{tikzcd}[sep=30pt]
            {[x,y]^\mathfrak{M}_\mathfrak{C} \, \Box \, I}
                \arrow[r,equal]
                \arrow[d,"1 j_x"']
            & {[x,y]^\mathfrak{M}_\mathfrak{C} \, \Box \, I}
                \arrow[d,equal]
            \\ {[x,y]^\mathfrak{M}_\mathfrak{C} \, \Box \, [x,x]^\mathfrak{M}_\mathfrak{C}}
                \arrow[r,"m_{x,x,y}"']
                \arrow[ur,Rightarrow,shorten <=20pt, shorten >=20pt,"\rho_{x,y}"]
            & {[x,y]^\mathfrak{M}_\mathfrak{C}}
        \end{tikzcd};\]

        \item For objects $x,y,z,w$ in $Ob(\mathfrak{M})$, an enriched associator in $\mathfrak{C}$ \[\begin{tikzcd}[sep=30pt]
            {([z,w]^\mathfrak{M}_\mathfrak{C} \, \Box \, [y,z]^\mathfrak{M}_\mathfrak{C}) \, \Box \, [x,y]^\mathfrak{M}_\mathfrak{C}}
                \arrow[r,"m_{y,z,w} 1"]
                \arrow[d,equal]
            & {[y,w]^\mathfrak{M}_\mathfrak{C} \, \Box \, [x,y]^\mathfrak{M}_\mathfrak{C}}
                \arrow[dd,"m_{x,y,w}"]
            \\{[z,w]^\mathfrak{M}_\mathfrak{C} \, \Box \,( [y,z]^\mathfrak{M}_\mathfrak{C} \, \Box \, [x,y]^\mathfrak{M}_\mathfrak{C})}
                \arrow[d,"1 m_{x,y,z}"']
            & {}
                \arrow[l,Rightarrow,shorten <=10pt, shorten >=10pt,"\pi_{x,y,z,w}"']
            \\ {[z,w]^\mathfrak{M}_\mathfrak{C} \, \Box \, [x,z]^\mathfrak{M}_\mathfrak{C}}
                \arrow[r,"m_{x,z,w}"']
            & {[x,w]^\mathfrak{M}_\mathfrak{C}}
        \end{tikzcd};\]
    \end{enumerate}
    subject to the following simplified conditions in \cite[Section 3.1]{GS}:
    \begin{enumerate}
        \item [a.] For objects $x,y,z,u,v$ in $Ob(\mathfrak{M})$, we have
    \end{enumerate}
    
    \settoheight{\diagramlength}{\includegraphics[height=20mm]{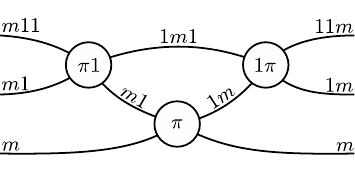}}
    
    \begin{equation}\label{eqn:Enriched2CategoryAssociativity}
    \begin{tabular}{@{}ccc@{}}
    
    \includegraphics[height=20mm]{Pictures/Appendex/EnrichedAssociativityleft.pdf} & \raisebox{0.45\diagramlength}{$=$} &
    \includegraphics[height=20mm]{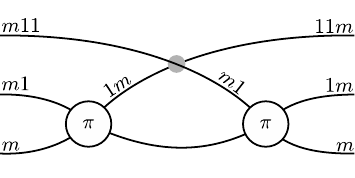}
    
    \end{tabular}
    \end{equation}

    \begin{enumerate}
        \item [] in $\mathbf{Hom}_{\mathfrak{C}}([u,v] \, \Box \, [z,u] \, \Box \, [y,z] \, \Box \, [x,y], [x,v])$,
        \item [b.] For objects $x,y,z$ in $Ob(\mathfrak{M})$, we have
    \end{enumerate}

    \settoheight{\diagramlength}{\includegraphics[height=20mm]{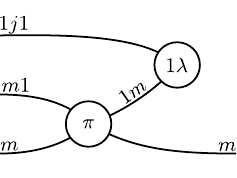}}
    
    \begin{equation}\label{eqn:Enriched2CategoryUnitality}
    \begin{tabular}{@{}ccc@{}}
    
    \includegraphics[height=20mm]{Pictures/Appendex/EnrichedUnitalityleft.pdf} & \raisebox{0.45\diagramlength}{$=$} &
    \includegraphics[height=20mm]{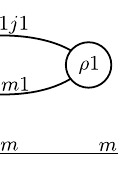}
    
    \end{tabular}
    \end{equation}

    \begin{enumerate}
        \item [] in $\mathbf{Hom}_{\mathfrak{C}}([y,z] \, \Box \, [x,y], [x,z])$.
    \end{enumerate}
\end{Definition}
    
\begin{Remark} \label{rmk:Underlying2CategoryOfEnriched2Category}
    Post-composing the enriched hom object with the forgetful 2-functor $\mathfrak{C} \to \mathbf{Cat}; \, x \mapsto \mathbf{Hom}_\mathfrak{C}(I,x)$, we obtain a $\mathbf{Cat}$-enriched 2-category structure on $\mathfrak{M}$. This is nothing else but the underlying 2-category structure on $\mathfrak{M}$, with $Ob(\mathfrak{M})$ as its genuine set of objects now.
\end{Remark}

\begin{Example}[{\cite[Proposition 4.1.1]{D4}}] \label{exmp:FiniteSemisimpleModule2CategoryIsEnriched} 
    Let $\mathfrak{C}$ be a fusion 2-category. Then every finite semisimple $\mathfrak{C}$-module 2-category $(\mathfrak{M},\Diamond)$ has a canonical enriched 2-category structure over $\mathfrak{C}$, determined by the tensor-hom 2-adjunction: for any object $z$ in $\mathfrak{C}$ and objects $x,y$ in $\mathfrak{M}$, one has the following 2-natural equivalence of categories: \[\Phi_{z,x;y}: \mathbf{Hom}_{\mathfrak{M}}(z \, \Diamond \, x,y) \simeq \mathbf{Hom}_{\mathfrak{C}}(z,[x,y]).\]
\end{Example}

\begin{Definition} \label{def:EnrichedMonoidal2Category}
    Let $\mathfrak{B}$ be a semistrict braided monoidal 2-category. We define a semistrict enriched monoidal 2-category over $\mathfrak{B}$ consisting of:
    \begin{enumerate}
        \item An underlying enriched 2-category $(\mathfrak{C},[-,-],m,j,\pi,\lambda,\rho)$ over $\mathfrak{B}$;
        
        \item For objects $x,y$ in $Ob(\mathfrak{C})$, a unique object $xy$ in $Ob(\mathfrak{C})$ called their monoidal product;
        
        \item A distinguished object $I_\mathfrak{C}$ in $Ob(\mathfrak{C})$ called the monoidal unit;
        
        \item For objects $x,y,z,w$ in $Ob(\mathfrak{C})$, an enriched product in $\mathfrak{B}$ \[\mu_{\left(\substack{x,y \\ z,w}\right)}: [x,y] \, \Box \, [z,w] \to [xz,yw];\]
        
        \item For objects $x,y,z,u,v,w$ in $Ob(\mathfrak{C})$, an enriched monoidal associator in $\mathfrak{B}$ \[\begin{tikzcd}[sep=40pt]
            {[x,y] \, \Box \, [z,u] \, \Box \, [v,w]}
                \arrow[r,"1 \mu"]
                \arrow[d,"\mu 1"']
            & {[x,y] \, \Box \, [zv,uw]}
                \arrow[d,"\mu"]
            \\ {[xz,yu] \, \Box \, [v,w]}
                \arrow[r,"\mu"']
                \arrow[ur,Rightarrow,shorten <=30pt, shorten >=30pt,"\alpha_{\left(\substack{x,y \\ z,u \\ v,w}\right)}"']
            & {[xzv,yuw]}
        \end{tikzcd};\]

        \item For objects $x,y$ in $Ob(\mathfrak{C})$, enriched monoidal unitors in $\mathfrak{B}$ \[\begin{tikzcd}
            {I_\mathfrak{B} \, \Box \, [x,y]}
                \arrow[d,equal]
                \arrow[r,"j_{I} 1"]
            & {[I_\mathfrak{C},I_\mathfrak{C}] \, \Box \, [x,y]}
                \arrow[d,"\mu"]
                \arrow[dl,Rightarrow,shorten <=10pt, shorten >=15pt,"l_{x,y}"]
            \\ {[x,y]}
                \arrow[r,equal]
            & {[x,y]}
        \end{tikzcd},\] \[\begin{tikzcd}
            {[x,y] \, \Box \, I_\mathfrak{B}}
                \arrow[d,equal]
                \arrow[r,"1 j_{I}"]
            & {[x,y] \, \Box \, [I_\mathfrak{C},I_\mathfrak{C}]}
                \arrow[d,"\mu"]
                \arrow[dl,Rightarrow,shorten <=10pt, shorten >=15pt,"r_{x,y}"]
            \\ {[x,y]}
                \arrow[r,equal]
            & {[x,y]}
        \end{tikzcd};\] 

        \item For objects $x,y,z,u,v,w$ in $Ob(\mathfrak{C})$, an enriched monoidal interchanger \[\begin{tikzcd}[sep=40pt]
            {[y,z] \, \Box \, [v,w] \, \Box \, [x,y] \, \Box \, [u,v]}
                \arrow[r,"\mu \mu"]
                \arrow[dd,"1 b 1"']
            & {[yv,zw] \, \Box \, [xu,yv]}
                \arrow[d,"m"]
            \\ {}
            & {[yvz,zwu]}
            \\ {[y,z] \, \Box \, [x,y] \, \Box \, [v,w] \, \Box \, [u,v]}
                \arrow[r,"m m"']
                \arrow[uur,Rightarrow,shorten <=60pt, shorten >=60pt,"\varepsilon_{\left(\substack{x,y,z \\ u,v,w}\right)}"']
            & {[x,z] \, \Box \, [u,w]}
                \arrow[u,"\mu"']
        \end{tikzcd};\]

        \item For objects $x,y$ in $Ob(\mathfrak{C})$, a product of enriched units in $\mathfrak{B}$ \[\begin{tikzcd}
            {I_\mathfrak{B} \, \Box \, I_\mathfrak{B} }
                \arrow[d,equal]
                \arrow[r,"j_x j_y"]
            & {[x,x] \, \Box \, [y,y]}
                \arrow[d,"\mu"]
                \arrow[dl,Rightarrow,shorten <=10pt, shorten >=15pt,"\iota_{x,y}"]
            \\ {I_\mathfrak{B} }
                \arrow[r,"j_{xy}"']
            & {[xy,xy]}
        \end{tikzcd};\]
    \end{enumerate}
    satisfying the following conditions:
    \begin{enumerate}
        \item[a.] For object $x$ in $Ob(\mathfrak{C})$, one always has $I_\mathfrak{C} \; x \equiv x \equiv x \; I_\mathfrak{C}$,
        
        \item[b.] For objects $x,y,z$ in $Ob(\mathfrak{C})$, one always has $(xy)z \equiv xyz \equiv x(yz)$,
        
        \item[c.] For objects $x,y,z,u,v,w,p,q$ in $Ob(\mathfrak{C})$, we have
    \end{enumerate}
    
    \settoheight{\diagramlength}{\includegraphics[height=20mm]{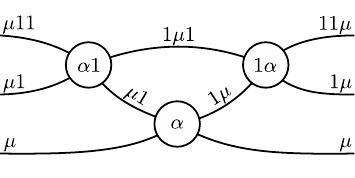}}
    
    \begin{equation}\label{eqn:EnrichedMonoidalAssociativity}
    \begin{tabular}{@{}ccc@{}}
    
    \includegraphics[height=20mm]{Pictures/Appendex/EnrichedMonoidalAssociativityleft.pdf} & \raisebox{0.45\diagramlength}{$=$} &
    \includegraphics[height=20mm]{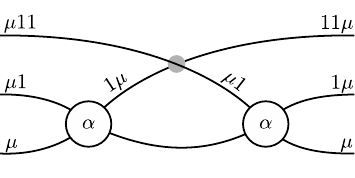}
    
    \end{tabular}
    \end{equation}

    \begin{enumerate}
        \item [] in $\mathbf{Hom}_{\mathfrak{B}}([x,y] \, \Box \, [z,u] \, \Box \, [v,w] \, \Box \, [p,q], [xzvp,yuwq])$,
        \item [d.] For objects $x,y,z,w$ in $Ob(\mathfrak{C})$, we have
    \end{enumerate}

    \settoheight{\diagramlength}{\includegraphics[height=20mm]{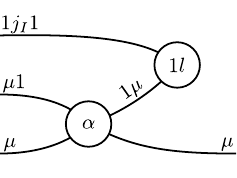}}
    
    \begin{equation}\label{eqn:EnrichedMonoidalUnitality}
    \begin{tabular}{@{}ccc@{}}
    
    \includegraphics[height=20mm]{Pictures/Appendex/EnrichedMonoidalUnitalityleft.pdf} & \raisebox{0.45\diagramlength}{$=$} &
    \includegraphics[height=20mm]{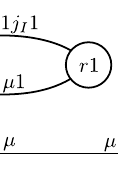}
    
    \end{tabular}
    \end{equation}

    \begin{enumerate}
        \item [] in $\mathbf{Hom}_{\mathfrak{B}}([x,y] \, \Box \, [z,w], [xz,yw])$,
        \item [e.] For objects $x,y,z$ in $Ob(\mathfrak{C})$, we have
    \end{enumerate}

    \settoheight{\diagramlength}{\includegraphics[height=33mm]{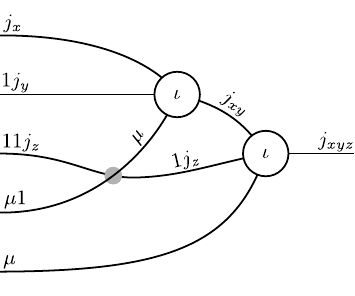}}
    
    \begin{equation}\label{eqn:ProductEnrichedUnitsAssociativity}
    \begin{tabular}{@{}ccc@{}}
    
    \includegraphics[height=33mm]{Pictures/Appendex/ProductOfEnrichedUnitsAssociatorLeft.pdf} & \raisebox{0.45\diagramlength}{$=$} &
    \includegraphics[height=33mm]{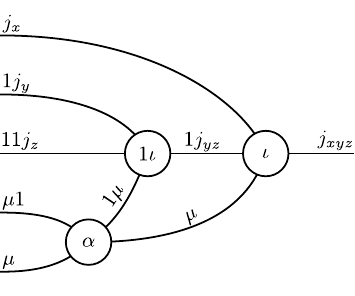}
    \end{tabular}
    \end{equation}

    \begin{enumerate}
        \item [] in $\mathbf{Hom}_{\mathfrak{B}}(I_\mathfrak{B}, [xyz,xyz])$,
        \item [f.] For object $x$ in $Ob(\mathfrak{C})$, we have
    \end{enumerate}

    \settoheight{\diagramlength}{\includegraphics[height=20mm]{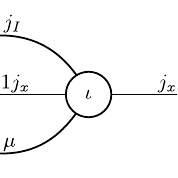}}
    
    \begin{equation}\label{eqn:ProductEnrichedUnitsUnitality1}
    \begin{tabular}{@{}ccc@{}}
    
    \includegraphics[height=20mm]{Pictures/Appendex/ProductOfEnrichedUnitsLeftUnitorLeft.pdf} & \raisebox{0.45\diagramlength}{$=$} &
    \includegraphics[height=20mm]{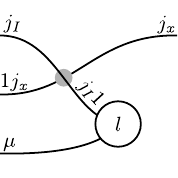}
    \end{tabular}
    \end{equation}

    \begin{enumerate}
        \item [] in $\mathbf{Hom}_{\mathfrak{B}}(I_\mathfrak{B}, [x,x])$, and
    \end{enumerate}

    \settoheight{\diagramlength}{\includegraphics[height=20mm]{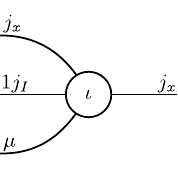}}
    
    \begin{equation}\label{eqn:ProductEnrichedUnitsUnitality2}
    \begin{tabular}{@{}ccc@{}}
    
    \includegraphics[height=20mm]{Pictures/Appendex/ProductOfEnrichedUnitsRightUnitorLeft.pdf} & \raisebox{0.45\diagramlength}{$=$} &
    \includegraphics[height=20mm]{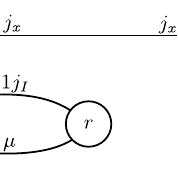}
    \end{tabular}
    \end{equation}

    \begin{enumerate}
        \item [] in $\mathbf{Hom}_{\mathfrak{B}}(I_\mathfrak{B}, [x,x])$,
        \item [g.] For objects $x,y,u,v$ in $Ob(\mathfrak{C})$, we have
    \end{enumerate}

    \settoheight{\diagramlength}{\includegraphics[height=40mm]{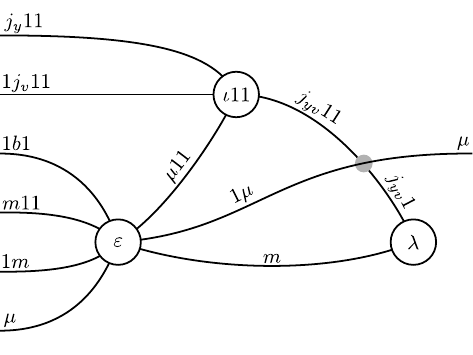}}
    
    \begin{equation}\label{eqn:InterchangerUnitality1}
    \begin{tabular}{@{}ccc@{}}
    
    \includegraphics[height=40mm]{Pictures/Appendex/InterchangerUnitality1Left.pdf} & \raisebox{0.45\diagramlength}{$=$} &
    \includegraphics[height=40mm]{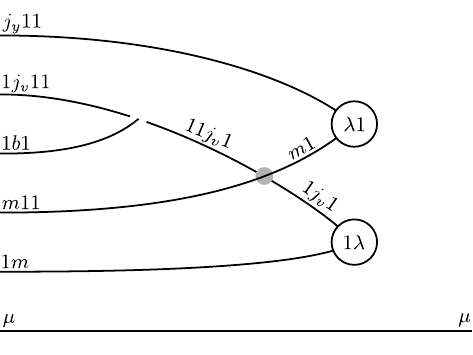}
    \end{tabular}
    \end{equation}

    \begin{enumerate}
        \item [] in $\mathbf{Hom}_{\mathfrak{B}}([y,y] \, \Box \, [v,v] \, \Box \, [x,y] \, \Box \, [u,v], [xu,yv])$, and
    \end{enumerate}

    \settoheight{\diagramlength}{\includegraphics[height=40mm]{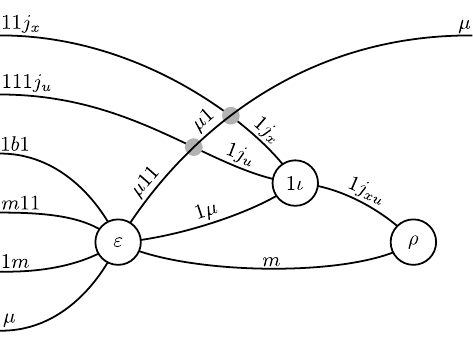}}

    \begin{equation}\label{eqn:InterchangerUnitality2}
        \begin{tabular}{@{}ccc@{}}
        
        \includegraphics[height=40mm]{Pictures/Appendex/InterchangerUnitality2Left.pdf} & \raisebox{0.45\diagramlength}{$=$} &
        \includegraphics[height=40mm]{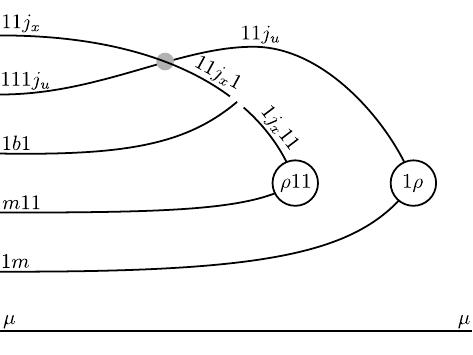}
        \end{tabular}
    \end{equation}

    \begin{enumerate}
        \item [] in $\mathbf{Hom}_{\mathfrak{B}}([x,y] \, \Box \, [u,v] \, \Box \, [x,x] \, \Box \, [u,u], [xu,yv])$,
        \item [h.] For objects $x,y,z$ in $Ob(\mathfrak{C})$, we have
    \end{enumerate}

    \settoheight{\diagramlength}{\includegraphics[height=40mm]{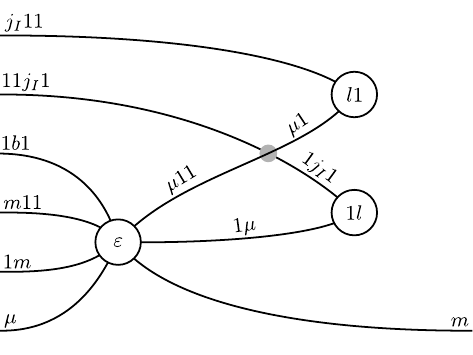}}
    
    \begin{equation}\label{eqn:InterchangerUnitality3}
    \begin{tabular}{@{}ccc@{}}
    
    \includegraphics[height=40mm]{Pictures/Appendex/InterchangerUnitality3Left.pdf} & \raisebox{0.45\diagramlength}{$=$} &
    \includegraphics[height=40mm]{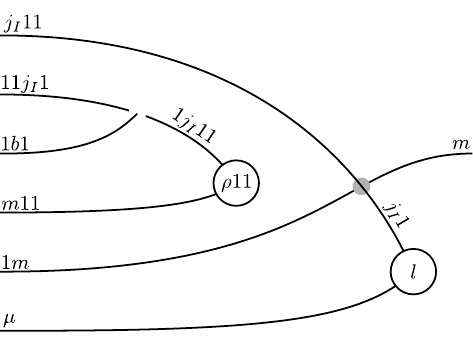}
    \end{tabular}
    \end{equation}

    \begin{enumerate}
        \item [] in $\mathbf{Hom}_{\mathfrak{B}}([I,I] \, \Box \, [y,z] \, \Box \, [I,I] \, \Box \, [x,y], [x,z])$, and
    \end{enumerate}

    \settoheight{\diagramlength}{\includegraphics[height=40mm]{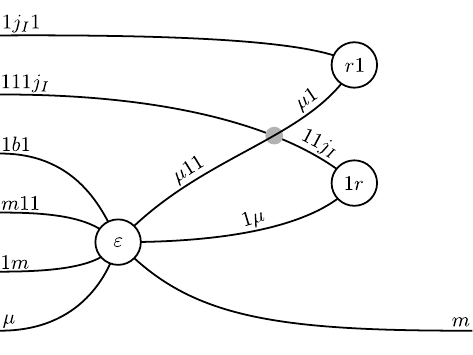}}

    \begin{equation}\label{eqn:InterchangerUnitality4}
        \begin{tabular}{@{}ccc@{}}
        
        \includegraphics[height=40mm]{Pictures/Appendex/InterchangerUnitality4Left.pdf} & \raisebox{0.45\diagramlength}{$=$} &
        \includegraphics[height=40mm]{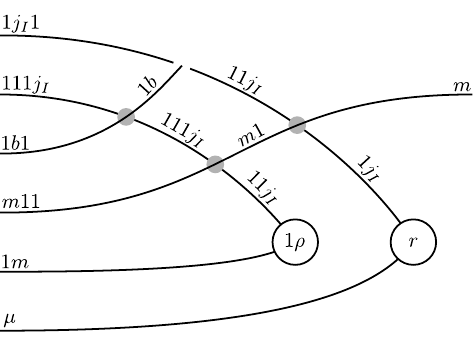}
        \end{tabular}
    \end{equation}

    \begin{enumerate}
        \item [] in $\mathbf{Hom}_{\mathfrak{B}}([y,z] \, \Box \, [I,I] \, \Box \, [x,y] \, \Box \, [I,I], [x,z])$,
        \item [i.] For objects $x,y,z,p,u,v,w,q$ in $Ob(\mathfrak{C})$, we have Equation (\ref{eqn:InterchangerAssociativity1}) in \[\mathbf{Hom}_{\mathfrak{B}}([z,p] \, \Box \, [w,q] \, \Box \, [y,z] \, \Box \, [v,w] \, \Box \, [x,y] \, \Box \, [u,v], [xu,pq]),\]
        
        \item [j.] For objects $x,y,z,u,v,w,p,q,s$ in $Ob(\mathfrak{C})$, we have Equation (\ref{eqn:InterchangerAssociativity2}) in \[\mathbf{Hom}_{\mathfrak{B}}([y,z] \, \Box \, [v,w] \, \Box \, [q,s] \, \Box \, [x,y] \, \Box \, [u,v] \, \Box \, [p,q], [xup,zws]).\]
    \end{enumerate}
\end{Definition}

\begin{landscape}
    \phantom{text}
    \begin{equation}\label{eqn:InterchangerAssociativity1}
    \begin{tabular}{c}
    \includegraphics[height=50mm]{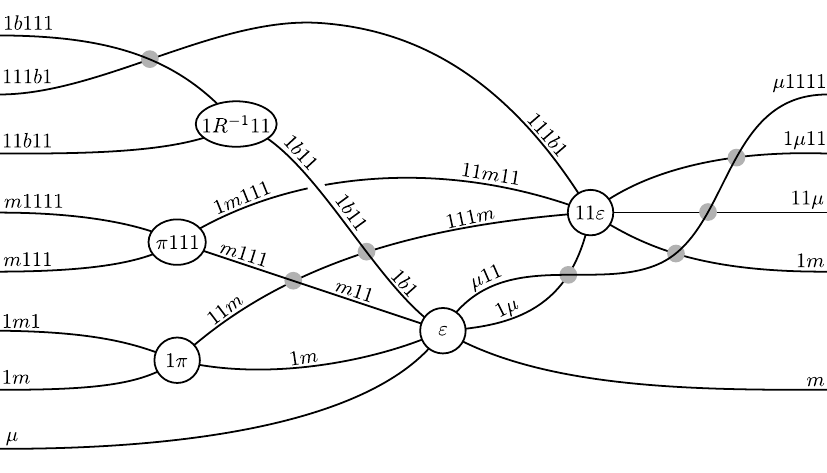} \\ \rotatebox{90}{$=$} \\
    \includegraphics[height=50mm]{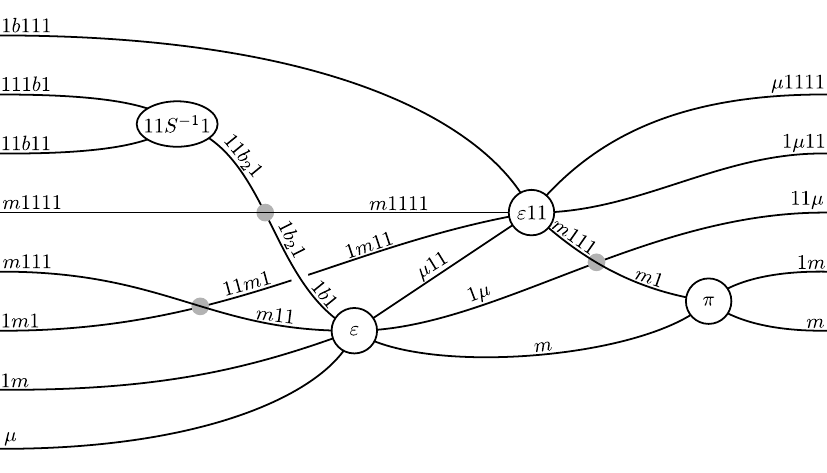}
    \end{tabular}
    \end{equation}

    \begin{equation}\label{eqn:InterchangerAssociativity2}
        \begin{tabular}{@{}c@{}}
        
        \includegraphics[height=50mm]{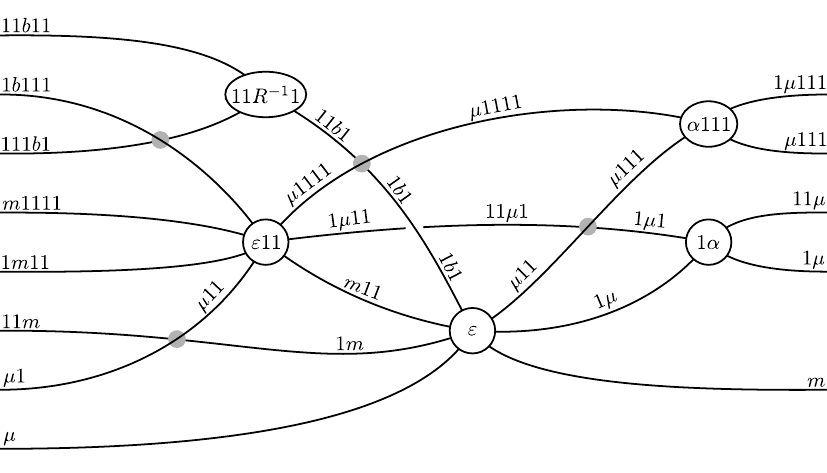} \\ \rotatebox{90}{$=$} \\
        \includegraphics[height=50mm]{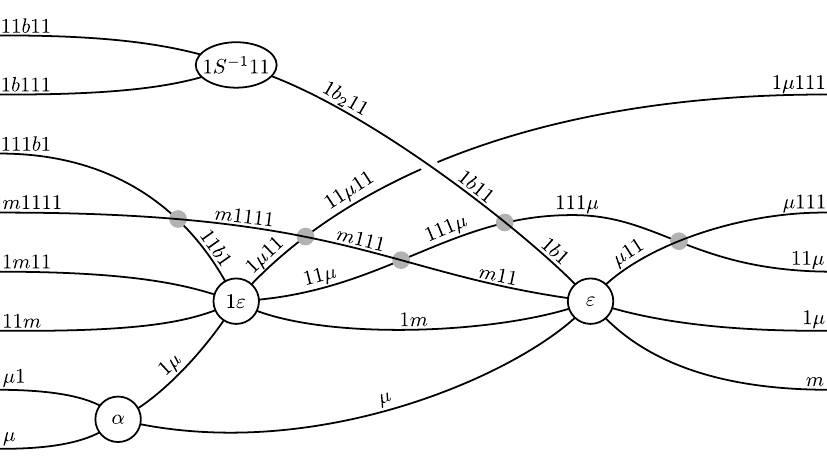}
        \end{tabular}
        \end{equation}
\end{landscape}

\begin{Corollary} \label{cor:EnrichedEndoHomOfUnitIsBraidedAlgebra}
    In particular, $m$ and $\mu$ endows two compatible algebra structures on $[I_\mathfrak{C},I_\mathfrak{C}]$, which turns it into a braided algebra in $\mathfrak{B}$ via Eckmann-Hilton argument.
\end{Corollary}

\begin{Example}
    If we take $\mathfrak{B} = \mathbf{Cat}$ with its Cartesian product as the symmetric monoidal structure, then a $\mathbf{Cat}$-enriched semistrict monoidal 2-category $\mathfrak{C}$ is equivalent to a semistrict monoidal 2-category in the normal sense. Firstly, the underlying 2-category of $\mathfrak{C}$ is given by the $\mathfrak{B}$-enriched 2-category structure on $\mathfrak{C}$ by Remark \ref{rmk:Underlying2CategoryOfEnriched2Category}. Secondly, the monoidal product and unit are given on the objects by definition. The interchanger is now induced by the followed 2-isomorphism:
    \settoheight{\diagramlength}{\includegraphics[height=40mm]{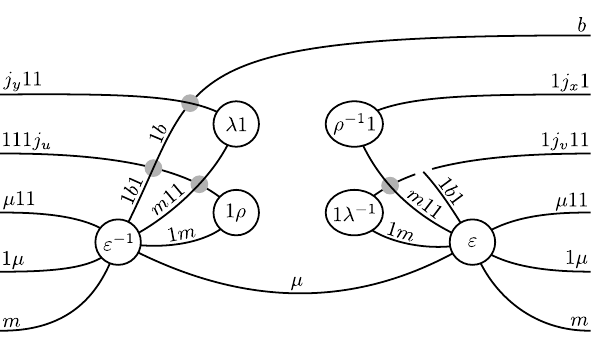}}

    \begin{equation}\label{eqn:UnderlyingInterchanger}
        \begin{tabular}{@{}ccc@{}}
        
        \includegraphics[height=40mm]{Pictures/Appendex/UnderlyingInterchanger.pdf}
        \end{tabular}
    \end{equation}
    in $\mathbf{Hom}_\mathfrak{B}([x,y] \, \Box \, [u,v],[xu,yv])$, i.e. for any 1-morphism $f:x \to y$ and $g:u \to v$ in $\mathfrak{C}$, it provides a 2-isomorphism filling the following square: \[\begin{tikzcd}
        {x \, \Box \, u}
            \arrow[r,"f 1"]
            \arrow[d,"1 g"']
        & {y \, \Box \, u}
            \arrow[d,"1 g"]
            \arrow[dl,Rightarrow,shorten <=10pt, shorten >=10pt]
        \\ {x \, \Box \, v}
            \arrow[r,"f 1"']
        & {y \, \Box \, v}
    \end{tikzcd}.\]
\end{Example}

\begin{Example}
    We would like to show that every fusion 2-category $\mathfrak{C}$ (assumed to be semistrict) is enriched as a semistrict monoidal 2-category over its Drinfeld center\footnote{See \cite[Section 2.1]{D9} for the definition of semistrict Drinfeld centers.} $\mathscr{Z}_1(\mathfrak{C})$. Let's denote the forgetful 2-functor by $\mathbf{U}:\mathscr{Z}_1(\mathfrak{C}) \to \mathfrak{C}$.

    Firstly, $\mathfrak{C}$ is a finite semisimple module 2-category over $\mathscr{Z}_1(\mathfrak{C})$ via the forgetful 2-functor preserving the monoidal structures: \[\mathscr{Z}_1(\mathfrak{C}) \to \mathbf{End}(\mathfrak{C}); \quad (x,b^x_{-},R^x_{-,-}) \mapsto x \, \Box \, -. \] Therefore, by Example \ref{exmp:FiniteSemisimpleModule2CategoryIsEnriched}, $\mathfrak{C}$ is enriched over $\mathscr{Z}_1(\mathfrak{C})$.

    Moreover, this enrichment is compatible with monoidal structure on $\mathfrak{C}$ and braided monoidal structure on $\mathscr{Z}_1(\mathfrak{C})$. \begin{enumerate}
        \item On objects, the monoidal product and unit are obviously compatible.
        
        \item On morphisms, the additional structures $(\mu,\alpha,l,r,\varepsilon,\iota)$ are induced by the universal property of the tensor-hom 2-adjunction.
    \end{enumerate}

    More specifically, for objects $x,y,u,v$ in $\mathfrak{C}$, the enriched product \[\mu_{\left(\substack{x,y \\ u,v}\right)}: [x,y] \, \Box \, [u,v] \to [xu,yv]\] lies in \[\mathbf{Hom}_{\mathscr{Z}_1(\mathfrak{C})}([x,y] \, \Box \, [u,v],[xu,yv]) \simeq \mathbf{Hom}_{\mathfrak{C}}(\mathbf{U}[x,y] \, \Box \, \mathbf{U}[u,v] \, \Box \, x \, \Box \, u,y \, \Box \, v)\] via tensor-hom 2-adjunction. Then we can apply the half-braiding on $[u,v]$ to exchange the positions of $\mathbf{U}[u,v]$ and $x$, which now lands in the hom category \[\mathbf{Hom}_{\mathfrak{C}}(\mathbf{U}[x,y] \, \Box \, x \, \Box \, \mathbf{U}[u,v] \, \Box \, u,y \, \Box \, v).\] Again, tensor hom 2-adjunction provides canonical evaluations $\mathbf{U}[x,y] \, \Box \, x \to y$ and $\mathbf{U}[u,v] \, \Box \, u \to v$, which induce the functor \[\mathbf{Hom}_{\mathfrak{C}}(y \, \Box \, v,y \, \Box \, v) \to \mathbf{Hom}_{\mathfrak{C}}(\mathbf{U}[x,y] \, \Box \, x \, \Box \, \mathbf{U}[u,v] \, \Box \, u,y \, \Box \, v).\] Lastly, the enriched product $\mu$ is provided by the identity 1-morphism $y \, \Box \, v \to y \, \Box \, v$ via the above sequence of functors.

    Curious readers are encouraged to work out the details of the remaining structures $(\alpha,l,r,\varepsilon,\iota)$ and check the compatibility conditions.
\end{Example}

Finally, this reaches the conclusion of the following proposition, which is the technical step skipped in the proof of Proposition \ref{prop:FullCenterOfUnit}.

\begin{Proposition}
    For any fusion 2-category $\mathfrak{C}$, there is an equivalence of enriched monoidal 2-categories over $\mathscr{Z}_1(\mathfrak{C})$ between $\mathfrak{C}$ and $\mathbf{Mod}_{\mathscr{Z}_1(\mathfrak{C})}([I,I])$.
\end{Proposition}

\bibliography{bibliography.bib}

\end{document}